\documentclass{amsart}

\usepackage{amsmath}
\usepackage[latin1]{inputenc}
\usepackage{amsfonts}
\usepackage{amssymb}
\usepackage{graphicx,epsfig}
\usepackage{amsthm}
\usepackage{amscd}
\usepackage[all]{xy}
\usepackage[latin1]{inputenc}

\usepackage[T1]{fontenc}

\usepackage{amsmath}

\usepackage{amsfonts}

\usepackage{amssymb}

\usepackage{amsthm}

\usepackage{graphics}

\newcommand{\mk}{\medskip}

\newcommand{\ZZ}{\mathbb{Z}}

\newcommand{\CC}{\mathbb{C}}

\newcommand{\NN}{\mathbb{N}}

\newcommand{\QQ}{\mathbb{Q}}

\newcommand{\Glie}{\mathfrak{g}}

\newcommand{\Yim}{\mathcal{Y}}

\newcommand{\Hlie}{\mathfrak{h}}

\newcommand{\demo}{\noindent {\it \small Proof:}\quad}

\renewcommand{\NN}{\ensuremath{\mathbb{N}}}

\renewcommand{\CC}{\ensuremath{\mathbb{C}}}

\renewcommand{\QQ}{\ensuremath{\mathbb{Q}}}

\newcommand{\U}{\mathcal{U}}

\newtheorem{thm}{Theorem}[section]

\newtheorem{defi}[thm]{Definition}

\newtheorem{cor}[thm]{Corollary}

\newtheorem{prop}[thm]{Proposition}

\newtheorem{lem}[thm]{Lemma}

\newtheorem{conj}[thm]{Conjecture}

\author{David Hernandez}

\thanks{\'Ecole Normale Sup\'erieure - DMA, 45, Rue d'Ulm F-75230 PARIS, Cedex 05  FRANCE
\newline email: David.Hernandez@ens.fr, URL: http://www.dma.ens.fr/$\sim$dhernand}

\usepackage[all]{xy}
\title{The $t$-Analogs of $q$-Characters at Roots of Unity for Quantum Affine algebras and Beyond}

\begin{document}

\begin{abstract} The $q$-characters were introduced by Frenkel and Reshetikhin \cite{Fre} to study finite dimensional representations of the untwisted quantum affine algebra $\U_q(\hat{\Glie})$ for $q$ generic. The $\epsilon$-characters at roots of unity were constructed by Frenkel and Mukhin \cite{Fre3} to study finite dimensional representations of various specializations of $\U_q(\hat{\Glie})$ at $q^s=1$. In the finite simply laced case Nakajima \cite{Naa}\cite{Nab} defined deformations of $q$-characters called $q,t$-characters for $q$ generic and also at roots of unity. The definition is combinatorial but the proof of the existence uses the geometric theory of quiver varieties which holds only in the simply laced case. In \cite{her} we proposed an algebraic general (non necessarily simply laced) new approach to $q,t$-characters for $q$ generic. In this paper we treat the root of unity case. Moreover we construct $q$-characters and $q,t$-characters for a large class of generalized Cartan matrices (including finite and affine cases except $A_1^{(1)}$, $A_2^{(2)}$) by extending the approach of \cite{her}. In particular we generalize the construction of analogs of Kazhdan-Lusztig polynomials at roots of unity of \cite{Nab} to those cases. We also study properties of various objects used in this article : deformed screening operators at roots of unity, $t$-deformed polynomial algebras, bicharacters arising from symmetrizable Cartan matrices, deformation of the Frenkel-Mukhin's algorithm. \end{abstract}

\maketitle

\tableofcontents

\pagestyle{myheadings}

\markboth{DAVID HERNANDEZ}{THE $t$-ANALOGS OF $q$-CHARACTERS AT ROOTS OF UNITY}

\section{Introduction} 

V.G. Drinfel'd \cite{Dri1} and M. Jimbo \cite{jim} associated, independently, to any symmetrizable Kac-Moody algebra $\Glie$ and any complex number $q\in\CC^*$ a Hopf algebra $\U_q(\Glie)$ called quantum group or quantum Kac-Moody algebra.

\noindent First we suppose that $q\in\CC^*$ is not a root of unity. In the case of a semi-simple Lie algebra $\Glie$ of rank $n$, the structure of the Grothendieck ring $\text{Rep}(\U_q(\Glie))$ of finite dimensional representations of the quantum finite algebra $\U_q(\Glie)$ is well understood. It is analogous to the classical case $q=1$. In particular we have ring isomorphisms:
$$\text{Rep}(\U_q(\Glie))\simeq \text{Rep}(\Glie)\simeq \ZZ [\Lambda]^W\simeq \ZZ[T_1,...,T_n]$$ 
deduced from the injective homomorphism of characters $\chi$:
$$\chi(V)=\underset{\lambda\in\Lambda}{\sum}\text{dim}(V_{\lambda})\lambda$$
where $V_{\lambda}$ are weight spaces of a representation $V$ and $\Lambda$ is the weight lattice.

\noindent For the general case of Kac-Moody algebras the picture is less clear. The representation theory of the quantum affine algebra $\U_q(\hat{\Glie})$ is of particular interest (see \cite{Cha}, \cite{Cha2}). In this case there is a crucial property of $\U_q(\hat{\Glie})$: it has two realizations, the usual Drinfel'd-Jimbo realization and a new realization (see \cite{Dri2} and \cite{bec}) as a quantum affinization of the quantum finite algebra $\U_q(\Glie)$.

\noindent To study the finite dimensional representations of $\U_q(\hat{\Glie})$ Frenkel and Reshetikhin \cite{Fre} introduced $q$-characters which encode the (pseudo)-eigenvalues of some commuting elements in the Cartan subalgebra $\U_q(\hat{\Hlie})\subset\U_q(\hat{\Glie})$ (see also \cite{kn}). The morphism of $q$-characters is an injective ring homomorphism:
$$\chi_q:\text{Rep}(\U_q(\hat{\Glie}))\rightarrow\ZZ[Y_{i,a}^{\pm}]_{i\in I,a\in\CC^*}$$
where $\text{Rep}(\U_q(\hat{\Glie}))$ is the Grothendieck ring of finite dimensional (type 1)-representations of $\U_q(\hat{\Glie})$ and $I=\{1,...,n\}$. In particular $\text{Rep}(\U_q(\hat{\Glie}))$ is commutative and isomorphic to $\ZZ[X_{i,a}]_{i\in I,a\in\CC^*}$.

\noindent The morphism of $q$-characters has a symmetry property analogous to the classical action of the Weyl group $\text{Im}(\chi)=\ZZ[\Lambda]^W$: Frenkel and Reshetikhin \cite{Fre} defined $n$ screening operators $S_i$ and showed that $\text{Im}(\chi_q)=\underset{i\in I}{\bigcap}\text{Ker}(S_i)$ for $\Glie=sl_2$. The result was proved by Frenkel and Mukhin for all finite $\Glie$ in \cite{Fre2}.

\noindent In the simply laced case Nakajima \cite{Naa}\cite{Nab} introduced $t$-analogs of $q$-characters. The motivations are the study of filtrations induced on representations by (pseudo)-Jordan decompositions, the study of the decomposition in irreducible modules of tensorial products and the study of cohomologies of certain quiver varieties. The morphism of $q,t$-characters is a $\ZZ[t^{\pm}]$-linear map 
$$\chi_{q,t}:\text{Rep}(\U_q(\hat{\Glie}))\rightarrow\ZZ[Y_{i,a}^{\pm},t^{\pm}]_{i\in I,a\in\CC^*}$$ 
which is a deformation of $\chi_q$ and multiplicative in a certain sense. A combinatorial axiomatic definition of $q,t$-characters is given. But the existence is non-trivial and is proved with the geometric theory of quiver varieties which holds only in the simply laced case. 

\noindent In \cite{her} we defined and constructed $q,t$-characters in the general (non necessarily simply laced) case with a new approach motivated by the non-commutative structure of $\U_q(\hat{\Hlie})\subset\U_q(\hat{\Glie})$, the study of screening currents of \cite{Freb} and of deformed screening operators $S_{i,t}$ of \cite{her01}. In particular we have a symmetry property: the image of $\chi_{q,t}$ is a completion of $\underset{i\in I}{\bigcap}\text{Ker}(S_{i,t})$.

\noindent The representation theory of the quantum affine algebras $\U_q(\hat{\Glie})$ depends crucially whether $q$ is a root of unity or not (see \cite{Cha3}). Frenkel and Mukhin \cite{Fre2} generalized $q$-characters at roots of unity : if $\epsilon$\label{epsilon} is a $s^{th}$-primitive root of unity the morphism of $\epsilon$-characters is\label{y}:
$$\chi_{\epsilon}:\text{Rep}(\U_{\epsilon}^{\text{res}}(\hat{\Glie}))\rightarrow\ZZ[Y_{i,a}^{\pm}]_{i\in I,a\in\CC^*}$$
where $\text{Rep}(\U_{\epsilon}^{\text{res}}(\hat{\Glie}))$ is the Grothendieck ring of finite dimensional (type 1)-representations of the restricted specialization $\U_{\epsilon}^{\text{res}}(\hat{\Glie})$ of $\U_q(\hat{\Glie})$ at $q=\epsilon$. In particular $\text{Rep}(\U_{\epsilon}^{\text{res}}(\hat{\Glie}))$ is commutative and isomorphic to $\ZZ[X_{i,a}]_{i\in I,a\in\CC^*}$.

\noindent Moreover $\chi_{\epsilon}$ can be characterized by
$$\chi_{\epsilon}(\underset{i\in I,l\in\ZZ/s\ZZ}{\prod}X_{i,\epsilon^l}^{x_{i,l}})=\tau_s(\chi_q(\underset{i\in I, 0\leq l\leq s-1}{\prod}X_{i,q^l}^{x_{i,[l]}}))$$
where $\tau_s:\ZZ[Y_{i,q^l}^{\pm}]_{i\in I,l\in\ZZ}\rightarrow\ZZ[Y_{i,\epsilon^l}^{\pm}]_{i\in I,l\in\ZZ/s\ZZ}$ is the ring homomorphism such that $\tau_s(Y_{i,q^l}^{\pm})=Y_{i,\epsilon^{[l]}}^{\pm}$ (for $l\in\ZZ$ we denote by $[l]$ its image in $\ZZ/s\ZZ$\label{lc}).

\noindent In the simply laced case Nakajima generalized the theory of $q,t$-characters at roots of unity with the help of quiver varieties \cite{Nab}.

In this paper we construct $q,t$-characters at roots of unity in the general (non necessarily simply laced) case by extending the approach of \cite{her}. As an application we construct analogs of Kazhdan-Lusztig polynomials at roots of unity in the same spirit as Nakajima did for the simply laced case. We also study properties of various objects used in this paper: deformed screening operators at roots of unity, $t$-deformed polynomial algebras, bicharacters arising from general symmetrizable Cartan matrices, deformation of the Frenkel-Mukhin's algorithm.

\noindent The construction is also extended beyond the case of a quantum affine algebra, that is to say by replacing the finite Cartan matrix by a generalized symmetrizable Cartan matrix: the construction of $q$-characters as well as $q,t$-characters (generic and roots of unity cases) is explained in this paper for (non necessarily finite) Cartan matrices such that $i\neq j\Rightarrow C_{i,j}C_{j,i}\leq 3$ (it includes finite and affine types except $A_1^{(1)}$, $A_2^{(2)}$). The notion of a quantum affinization is more general than the construction of a quantum affine algebra from a quantum finite algebra: it can be extended to any general symmetrizable Cartan matrix (see \cite{Naams}). For example for an affine Cartan matrix one gets a quantum toroidal algebra (see \cite{var}). In general a quantum affinization is not a quantum Kac-Moody algebra and few is known about the representation theory outside the quantum affine algebra case. However for an integrable representation one can define $q$-characters as Frenkel-Reshetikhin did for quantum affine algebras. So the $q$-characters constructed in this paper for some generalized symmetrizable Cartan matrix are to be linked with representation theory of the associated quantum affinization. We will address further developments on this point in a separate publication.

\noindent This paper is organized as follows: after some backgrounds in section \ref{bckg}, we generalize in section \ref{un} the construction of $t$-deformed polynomial algebras of \cite{her} to the root of unity case. We give a ``concrete'' construction using Heisenberg algebras. We show that this twisted multiplication can also be ``abstractly'' defined with two bicharacters $d_1$, $d_2$ as Nakajima did for the simply laced case (for which there is only one bicharacter $d_1=d_2$). 

\noindent In section \ref{quatre} we remind how $q,t$-characters are constructed for $q$ generic and $C$ finite in \cite{her}. We extend the construction of $q$-characters and of $q,t$-characters to symmetrizable (non necessarily finite) Cartan matrices such that $i\neq j\Rightarrow C_{i,j}C_{j,i}\leq 3$, in particular for affine Cartan matrices (except $A_1^{(1)}$ and $A_2^{(2)}$). The $q,t$-characters can be computed by the algorithm described in \cite{her} which is a deformation of the algorithm of Frenkel-Mukhin \cite{Fre2}.

\noindent In section \ref{rofcase} we construct $q,t$-characters at roots of unity. Let us explain the crucial technical point of this section: we can not use directly a $t$-deformation of the definition of Frenkel-Mukhin because there is no analog of $\tau_s$ which is an algebra homomorphism for the $t$-deformed structures. But we can construct $\tau_{s,t}$ which is multiplicative for some ordered products (see section \ref{deficonst}). In particular $\tau_{s,t}$ has nice properties and we can define $\chi_{\epsilon,t}$ such that ``$\chi_{\epsilon,t}=\tau_{s,t}\circ \chi_{q,t}$''. We give properties of $\chi_{\epsilon,t}$ analogous to the property of $\chi_{\epsilon}$ (proposition \ref{axiomes}, theorems \ref{qtint} and \ref{axquat}). In particular in the $ADE$-case we get a formula which is Axiom 4 of \cite{Nab}, and so the construction coincides with the construction of \cite{Nab} for the $ADE$-case.

\noindent In section \ref{app} we give some applications about Kazhdan-Lusztig polynomials and quantization of the Grothendieck ring. If $C$ is finite the technical point in the root of unity case is to show that the algorithm produces a finite number of dominant monomials. We give a conjecture about the multiplicity of an irreducible module in a standard module at roots of unity. For the $ADE$-case it is a result of Nakajima \cite{Nab}. An analogous conjecture was given in \cite{her} for $q$ generic. We also study the non finite cases.

\noindent In section \ref{complements} we give some complements: first we discuss the finiteness of the algorithm; at $t=1$ it stops if $C$ is finite and it does not stop if $C$ is affine. We relate the structure of the deformed ring in the affine $A_r^{(1)}$-case to the structure of quantum toroidal algebras. We study some combinatorial properties of the Cartan matrices which are related to the bicharacters $d_1$ and $d_2$ (propositions \ref{facile}, \ref{faciledeux}, \ref{qsym} and theorem \ref{result}).

\noindent For convenience of the reader we give at the end of this article an index of notations defined in the main body of the text.

\noindent In the course of writing this paper we were informed by H. Nakajima that the $t$-analogs of $q$-characters for some quantum toroidal algebras are also mentioned in the remark 6.9 of \cite{Nad}. This incited us to add the construction of analogs of Kazdhan-Lusztig polynomials at roots of unity also in the non finite cases (section \ref{nonfini}). 

\noindent {\bf Acknowledgments.} The author would like to thank H. Nakajima for useful comments on a previous version of this paper, N. Reshetikhin and M. Rosso for encouraging him to study the root of unity cases, and M. Varagnolo for indications on quantum toroidal algebras.

\section{Background}\label{bckg}

\subsection{Cartan matrices}\label{backcartan}

A generalized Cartan matrix is $C=(C_{i,j})_{1\leq i,j\leq n}$ such that $C_{i,j}\in\ZZ$ and\label{carmat}:
$$C_{i,i}=2$$
$$i\neq j\Rightarrow C_{i,j}\leq 0$$
$$C_{i,j}=0\Leftrightarrow C_{j,i}=0$$
Let $I=\{1,...,n\}$.

\noindent $C$ is said to be decomposable if it can be written in the form $C=P\begin{pmatrix}A&0\\0&B\end{pmatrix}P^{-1}$ where $P$ is a permutation matrix, $A$ and $B$ are square matrices. Otherwise $C$ is said to be indecomposable.

\noindent $C$ is said to be symmetrizable if there is a matrix $D=\text{diag}(r_1,...,r_n)$ ($r_i\in\NN^*$)\label{ri} such that $B=DC$\label{symcar} is symmetric. In particular if $C$ is symmetric then it is symmetrizable with $D=I_n$.

\noindent If $C$ is indecomposable and symmetrizable then there is a unique choice of $r_1,...,r_n>0$ such that $r_1\wedge ...\wedge r_n=1$: indeed if $C_{j,i}\neq 0$ we have the relation $r_i=\frac{C_{j,i}}{C_{i,j}}r_j$.

\noindent In the following $C$ is a symmetrizable and indecomposable generalized Cartan matrix. For example:

\noindent $C$ is said to be of finite type if all its principal minors are positive (see \cite{bou} for a classification).

\noindent $C$ is said to be of affine type if all its proper principal minor are positive and $\text{det}(C)=0$ (see \cite{kac} for a classification).

\noindent Let $r^{\vee}=\text{max}\{r_i-1-C_{j,i})/i\neq j\}\cup\{1\}$\label{rvee}. If $C$ is finite we have $r^{\vee}=\text{max}\{r_i/i\in I\}=\text{max}\{-C_{i,j}/i\neq j\}$. In particular if $C$ is of type $ADE$ we have $r^{\vee}=1$, if $C$ is of type $B_lC_lF_4$ ($l\geq 2$) we have $r^{\vee}=2$, if $C$ of type $G_2$ we have $r^{\vee}=3$.

\noindent Let $z$\label{z} be an indeterminate and $z_i=z^{r_i}$. The matrix $C(z)=(C_{i,j}(z))_{1\leq i,j\leq n}$\label{cz} with coefficients in $\ZZ[z^{\pm}]$ is defined by $C_{i,i}(z)=[2]_{z_i}=z_i+z_i^{-1}$ and $C_{i,j}=[C_{i,j}]_z$ for $i\neq j$ where for $l\in\ZZ$, we use the notation:
$$[l]_z=\frac{z^l-z^{-l}}{z-z^{-1}}\text{ ($=z^{-l+1}+z^{-l+3}+...+z^{l-1}$ for $l\geq 1$)}$$ 
Let $B(z)=D(z)C(z)$\label{bz} where $D(z)$ is the diagonal matrix $D_{i,j}(z)=\delta_{i,j}[r_i]_z$, that is to say $B_{i,j}(z)=[r_i]_zC_{i,j}(z)$.

\noindent In particular, the coefficients of $C(z)$ and $B(z)$ are symmetric Laurent polynomials (invariant under $z\mapsto z^{-1}$).

\noindent In the following we suppose that $\text{det}(C(z))\neq 0$. It includes finite and affine Cartan matrices (if $C$ is of type $A_1^{(1)}$ we set $r_1=r_2=2$) and also the matrices such that $i\neq j\Rightarrow C_{i,j}C_{j,i}\leq 3$ which will appear later (see lemma \ref{condcinv} and section \ref{cartancompl} for complements).

\subsection{Quantum affine algebras}\label{qaa} In the following $q$ is a complex number $q\in\CC^*$\label{qz}. If $q$ is not a root of unity we set $s=0$ and we say that $q$ is generic. Otherwise $s\geq 1$ is set such that $q$ is a $s^{th}$ primitive root of unity.\label{s}

\noindent We suppose in this section that $C$ is finite. We refer to \cite{Fre3} for the definition of the untwisted quantum affine algebra $\U_q(\hat{\Glie})$\label{qadefi} associated to $C$ (for $q$ generic) and of the restricted specialization $\U_{\epsilon}^{\text{res}}(\hat{\Glie})$ of $\U_q(\hat{\Glie})$ at $q=\epsilon$ (for $\epsilon$ root of unity).

\noindent We briefly describe the construction of $\U_{\epsilon}^{\text{res}}(\hat{\Glie})$ from $\U_q(\hat{\Glie})$: we consider a $\ZZ[q,q^{-1}]$-subalgebra of $\U_q(\hat{\Glie})$ containing the $(x_i^{\pm})^{(r)}=\frac{(x_i^{\pm})^s}{[r]_{q_i}!}$ (where $[r]_q!=[r]_q[r-1]_q...[1]_q$) for some generators $x_i^{\pm}$, and we set $q=\epsilon$.

\noindent One can define a Hopf algebra structure on $\U_q(\hat{\Glie})$ and $\U_{\epsilon}^{\text{res}}(\hat{\Glie})$, and so we consider the Grothendieck ring of finite dimensional (type 1)-representations: $\text{Rep}(\U_q(\hat{\Glie}))$ and $\text{Rep}(\U_{\epsilon}^{\text{res}}(\hat{\Glie}))$. 

\noindent The morphism of $q$-characters $\chi_q$\label{chiqdefi} (Frenkel-Reshetikhin \cite{Fre}) and the morphism of $\epsilon$-character $\chi_{\epsilon}$ 
\\(Frenkel-Mukhin \cite{Fre3}) are injective ring homomorphisms:
$$\chi_q:\text{Rep}(\U_q(\hat{\Glie}))\rightarrow\ZZ[Y_{i,a}^{\pm}]_{i\in I,a\in\CC^*}\text{ , }\chi_{\epsilon}:\text{Rep}(\U_{\epsilon}^{\text{res}}(\hat{\Glie}))\rightarrow\ZZ[Y_{i,a}^{\pm}]_{i\in I,a\in\CC^*}$$
In particular $\text{Rep}(\U_q(\hat{\Glie}))$ and $\text{Rep}(\U_{\epsilon}^{\text{res}}(\hat{\Glie}))$ are commutative and isomorphic to $\ZZ[X_{i,a}]_{i\in I,a\in\CC^*}$.

\noindent Frenkel and Mukhin \cite{Fre2}\cite{Fre3} have proven that for $i\in I$, $a\in\CC^*$:
$$\chi_q(X_{i,a})\in\ZZ[Y_{i,aq^m}^{\pm}]_{i\in I, m\in\ZZ}\text{ and }\chi_{\epsilon}(X_{i,a})\in\ZZ[Y_{i,a\epsilon^m}^{\pm}]_{i\in I, m\in\ZZ}$$

\noindent Indeed it suffices to study (see \cite{her} for details)\label{eps}:
$$\chi_q:\text{Rep}=\ZZ[X_{i,l}]_{i\in I,l\in\ZZ}\rightarrow \Yim=\ZZ[Y_{i,l}^{\pm}]_{i\in I,l\in\ZZ}$$
\noindent (where $X_{i,l}=X_{i,q^l}$, $Y_{i,l}^{\pm}=Y_{i,q^l}^{\pm}$)\label{rep}, and:
$$\chi_{\epsilon}^s:\text{Rep}^s=\ZZ[X_{i,l}]_{i\in I,l\in\ZZ/s\ZZ}\rightarrow\Yim^s=\ZZ[Y_{i,l}^{\pm}]_{i\in I,l\in\ZZ/s\ZZ}$$
\noindent (where $X_{i,l}=X_{i,\epsilon^l}$\label{reps}, $Y_{i,l}^{\pm}=Y_{i,\epsilon^l}^{\pm}$).

\section{$t$-deformed polynomial algebras}\label{un}

\subsection{The $t$-deformed algebra $\hat{\Yim}_t^s$}\label{trois} In this section we generalize at roots of unity the construction of \cite{her} of $t$-deformed polynomial algebras. 

\subsubsection{Construction}\label{defihcal} 

In this section we suppose that $B(z)$ is symmetric.

\begin{defi}\label{hatcalh} $\hat{\mathcal{H}}$\label{zq} is the $\CC$-algebra defined by generators $a_i[m],y_i[m]$\label{aim} ($i\in I, m\in\ZZ-\{0\}$)\label{yim}, central elements $c_r$\label{cr} ($r>0$) and relations ($i,j\in I,m,r\in\ZZ-\{0\}$):
\begin{equation}\label{aa}[a_i[m],a_j[r]]=\delta_{m,-r}(q^m-q^{-m})B_{i,j}(q^m)c_{|m|}\end{equation}
\begin{equation}\label{ay}[a_i[m],y_j[r]]=(q^{mr_i}-q^{-r_im})\delta_{m,-r}\delta_{i,j}c_{|m|}\end{equation}
\begin{equation}\label{yy}[y_i[m],y_j[r]]=0\end{equation}\end{defi}

\noindent Let $\hat{\mathcal{H}}_h=\hat{\mathcal{H}}[[h]]$\label{zqh}. For $i\in I$, $l\in\ZZ/s\ZZ$ we define $\tilde{Y}_{i,l}^{\pm},\tilde{A}_{i,l}^{\pm},t_l^{\pm}\in\hat{\mathcal{H}}_h$\label{tail} such that:
$$\tilde{Y}_{i,l}=\text{exp}(\underset{m>0}{\sum}h^m q^{lm}y_i[m])\text{exp}(\underset{m>0}{\sum}h^m q^{-lm}y_i[-m])$$
$$\tilde{A}_{i,l}=\text{exp}(\underset{m>0}{\sum}h^m q^{lm}a_{i}[m])\text{exp}(\underset{m>0}{\sum}h^m q^{-lm}a_i[-m])$$
$$t_l=\text{exp}(\underset{m>0}{\sum}h^{2m}q^{lm}c_m)$$
and for $R=\underset{l\in\ZZ}{\sum}R_lz^l\in\ZZ[z^{\pm}]$:
$$t_{R}=\underset{l\in\ZZ}{\prod}t_l^{R_l}=\text{exp}(\underset{m>0}{\sum}h^{2m}R(q^{m})c_m)\in\hat{\mathcal{H}}_h$$\label{tr}
Note that the root of unity condition, that is to say $s\geq 1$, is a periodic condition ($\tilde{Y}_{i,l+s}=\tilde{Y}_{i,l}$). 

\begin{lem}\label{reltr}(\cite{her}) We have the following relations in $\hat{\mathcal{H}}_h$:
$$\tilde{A}_{i,l}\tilde{Y}_{j,k}\tilde{A}_{i,l}^{-1}\tilde{Y}_{j,k}^{-1}
=t_{\delta_{i,j}(z^{-r_i}-z^{r_i})(-z^{(l-k)}+z^{(k-l)})}$$
$$\tilde{A}_{i,l}\tilde{A}_{j,k}\tilde{A}_{i,l}^{-1}\tilde{A}_{j,k}^{-1}=t_{B_{i,j}(z)(z^{-1}-z)(-z^{(l-k)}+z^{(k-l)})}$$
\end{lem}

\begin{defi}\label{yu} $\hat{\Yim}_u^s$ is the $\ZZ$-subalgebra of $\hat{\mathcal{H}}_h$ generated by the $\tilde{Y}_{i,l},\tilde{A}_{i,l}^{-1},t_l$ ($i\in I,l\in\ZZ/z\ZZ$).\end{defi}

\noindent Note that if $s\geq 1$, the elements $\tilde{A}_{i,0}^{-1}\tilde{A}_{i,1}^{-1}...\tilde{A}_{i,s-1}^{-1}$ and $\tilde{Y}_{i,0}\tilde{Y}_{i,1}...\tilde{Y}_{i,s-1}$ are central in $\hat{\Yim}_u$.

\begin{defi}\label{tyt} $\hat{\Yim}_t^s$ is the quotient-algebra of $\hat{\Yim}_u^s$ by relations $t_l=1$ if $l\in\ZZ/s\ZZ-\{0\}$. \end{defi}

\noindent We keep the notations $\tilde{Y}_{i,l},\tilde{A}_{i,l}^{-1}$ for their image in $\hat{\Yim}_t^s$. We denote by $t$\label{t} the image of $t_0=\text{exp}(\underset{m>0}{\sum}h^{2m}c_m)$ in $\hat{\Yim}_t^s$. In particular the image of $t_R$ is $t^{R_0}$. We denote by $\hat{\Yim}_t=\hat{\Yim}_t^0$ the algebra in the generic case.

\subsubsection{Structure} For $a,b\in\ZZ/s\ZZ$, let $\delta_{a,b}=1$ if $a=b$ and $\delta_{a,b}=0$ if $a\neq b$.

\noindent The following theorem gives the structure of $\hat{\Yim}_t^s$: 

\begin{thm}\label{dessus} The algebra $\hat{\Yim}_t^s$ is defined by generators $\tilde{Y}_{i,l},\tilde{A}_{i,l}^{-1},t^{\pm}$ $(i\in I,l\in\ZZ/s\ZZ)$ and relations ($i,j\in I, k,l\in\ZZ/s\ZZ$):
$$\tilde{Y}_{i,l}\tilde{Y}_{j,k}=\tilde{Y}_{j,k}\tilde{Y}_{i,l}$$
$$\tilde{A}_{i,l}^{-1}\tilde{A}_{j,k}^{-1}=t^{\alpha(i,l,j,k)}\tilde{A}_{j,k}^{-1}\tilde{A}_{i,l}^{-1}$$
$$\tilde{Y}_{j,k}\tilde{A}_{i,l}^{-1}=t^{\beta(i,l,j,k)}\tilde{A}_{i,l}^{-1}\tilde{Y}_{j,k}$$
where $\alpha,\beta : (I\times\ZZ/s\ZZ)^2\rightarrow\ZZ$\label{alpha} are given by ($l,k\in\ZZ/s\ZZ$, $i,j\in I$)\label{beta}:
$$\alpha(i,l,i,k)=2(\delta_{l-k,-2r_i}-\delta_{l-k,2r_i})$$
$$\alpha(i,l,j,k)=2\underset{r=C_{i,j}+1,C_{i,j}+3,...,-C_{i,j}-1}{\sum}(\delta_{l-k,r+r_i}-\delta_{l-k,r-r_i})\text{ (if $i\neq j$)}$$
$$\beta(i,l,j,k)=2\delta_{i,j}(-\delta_{l-k,r_i}+\delta_{l-k,-r_i})$$
\end{thm}

\noindent Note that for $i,j\in I$ and $l,k\in\ZZ/s\ZZ$ we have $\alpha(i,l,j,k)=-\alpha(j,k,i,l)$ and $\beta(i,l,j,k)=-\beta(j,k,i,l)$.

\noindent This theorem is a generalization of theorem 3.11 of \cite{her}. It is proved in the same way except for lemma 3.7 of \cite{her} whose proof is changed at roots of unity: for $N\geq 1$ we denote by $\ZZ_N[z]\subset\ZZ[z]$ the subset of polynomials of degree lower that $N$. The following lemma is a generalization of lemma 3.7 of \cite{her} at roots of unity :

\begin{lem}\label{indgene} We suppose that $s\geq 1$. Let $J=\{1,...,r\}$ be a finite set of cardinal $r$ and $\Lambda$ be the polynomial commutative algebra 
\\$\Lambda=\CC[\lambda_{j,m}]_{j\in J,m\geq 0}$. For $R=(R_1,...,R_r)\in {\ZZ_{s-1}[z]}^{r}$, consider:
$$\Lambda_R=\text{exp}(\underset{j\in I,m>0}{\sum}h^mR_j(q^m)\lambda_{j,m})\in\Lambda[[h]]$$
Then the $(\Lambda_R)_{R\in {\ZZ_{s-1}[z]}^r}$ are $\CC$-linearly independent. In particular the $\Lambda_{j,l}=\Lambda_{(0,...,0,z^l,0,...,0)}$ ($j\in I$, $0\leq l\leq s-1$) are $\CC$-algebraically independent.\end{lem}

\demo Suppose we have a linear combination ($\mu_R\in\CC$, only a finite number of $\mu_R\neq 0$):
$$\underset{R\in {\ZZ_{s-1}[z]}^r}{\sum}\mu_R\Lambda_R=0$$
In the proof of lemma 3.7 of \cite{her} we saw that for $N\geq 0$, $j_1,...,j_N\in J$, $l_1,...,l_N>0$, $\alpha_1,...,\alpha_N\in\CC$ we have:
$$\underset{R\in {\ZZ_{s-1}[z]}^r/R_{j_1}(q^{l_1})=\alpha_1,...,R_{j_{N}}(q^{l_{R}})=\alpha_{N}}{\sum}\mu_{R}=0$$
We set $N=sr$ and 
$$((j_1,l_1),...,(j_N,l_N))=((1,1),(1,2),...,(1,s),(2,1),...,(2,s),(3,1),...,(r,s))$$
We get for all $\alpha_{j,l}\in\CC$ ($j\in J,1\leq l\leq L$):
$$\underset{R\in\ZZ_{s-1}[z]^r/\forall j\in J,1\leq l\leq s, R_{j}(q^{l})=\alpha_{j,l}}{\sum}\mu_{R}=0$$
It suffices to show that there is at most one term is this sum. But consider $P,Q\in\ZZ_{s-1}[z]$ such that for all $1\leq l\leq s$, $P(q^l)=P'(q^l)$. As $q$ is primitive the $q^l$ are different and so $P-P'=0$.\qed

\subsection{Bicharacters, monomials and involution}

\subsubsection{Presentation with bicharacters}\label{presbica} The definition of the algebra $\hat{\Yim}_t^s$ with the Heisenberg algebra $\hat{\mathcal{H}}$ is a ``concrete'' construction. It can also be defined ``abstractly'' with bicharacters in the same spirit as Nakajima \cite{Nab} did for the simply laced case :

\noindent We define $\pi_+$\label{piplus} as the algebra homomorphism:
$$\pi_+:\hat{\Yim}_t^s\rightarrow \hat{\Yim}^s=\ZZ[Y_{i,l},A_{i,l}^{-1}]_{i\in I,l\in\ZZ/s\ZZ}$$\label{yil}
such that $\pi_+(\tilde{Y}_{i,l}^{\pm})=Y_{i,l}^{\pm}$, $\pi_+(\tilde{A}_{i,l}^{\pm})=A_{i,l}^{\pm}$ and $\pi_+(t)=1$ ($\hat{\Yim}^s$ is commutative).

\noindent We say that $m\in \hat{\Yim}_t^s$ is a $\hat{\Yim}_t^s$-monomial if it is a product of the $\tilde{A}_{i,l}^{-1},\tilde{Y}_{i,l},t^{\pm}$. For $m$ a $\hat{\Yim}_t^s$-monomial, $i\in I$, $l\in\ZZ/s\ZZ$ we define $y_{i,l}(m), v_{i,l}(m)\geq 0$ such that $\pi_+(m)=\underset{i\in I,l\in\ZZ/s\ZZ}{\prod}Y_{i,l}^{y_{i,l}}A_{i,l}^{-v_{i,l}}$. In order to simplify the formulas for a Laurent polynomial let $P(z)=\underset{k\in\ZZ}{\sum}P_kz^k\in\ZZ[z^{\pm}]$\label{ope} ($i\in I, l\in\ZZ/s\ZZ$):
$$(P(z))_{op}\mathcal{V}_{i,l}(m)=\underset{k\in\ZZ}{\sum}P_kv_{i,l+[k]}(m)$$
We define $u_{i,l}(m)\in\ZZ$ by\label{uil} :
$$u_{i,l}(m)=y_{i,l}(m)-\underset{j\in I}{\sum}(C_{i,j}(z))_{op}\mathcal{V}_{j,l}(m)$$
In particular if $C_{i,j}=0$ we have $u_{i,l}(\tilde{A}_{j,k}^{-1})=0$ and if $C_{i,j}<0$:
$$u_{i,l}(\tilde{A}_{j,k}^{-1})=-([C_{i,j}]_z)_{op}\mathcal{V}_{j,l}(\tilde{A}_{j,k}^{-1})=\underset{r=C_{i,j}+1... -C_{i,j}-1}{\sum}\delta_{l+r,k}$$

\noindent In the $ADE$-case the coefficients of $C$ are $-1,0$ or $2$, and we have the expression:
$$u_{i,l}(m)=y_{i,l}(m)-[z+z^{-1}]_{op}\mathcal{V}_{i,l}(m)+\underset{j\in I/C_{i,j}=-1}{\sum}v_{j,l}(m)$$
$$=y_{i,l}(m)-v_{i,l+1}(m)-v_{i,l-1}(m)+\underset{j\in I/C_{i,j}=-1}{\sum}v_{j,l}(m)$$
which is the formula used in \cite{Nab}.

\begin{defi}\label{defidd} For $m_1,m_2$ $\hat{\Yim}_t^s$-monomials we define:
$$d_1(m_1,m_2)=\underset{i\in I,l\in\ZZ/s\ZZ}{\sum}v_{i,l+r_i}(m_1)u_{i,l}(m_2)+y_{i,l+r_i}(m_1)v_{i,l}(m_2)$$
$$d_2(m_1,m_2)=\underset{i\in I,l\in\ZZ/s\ZZ}{\sum}u_{i,l+r_i}(m_1)v_{i,l}(m_2)+v_{i,l+r_i}(m_1)y_{i,l}(m_2)$$
\end{defi}\label{d}

\noindent For $m$ a $\hat{\Yim}_t^s$-monomial we have always $d_1(m,m)=d_2(m,m)$ (see section \ref{cartancompl}). In the $ADE$-case we have $d_1=d_2$ and it is the bicharacter of Nakajima \cite{Nab}.

\begin{prop} For $m_1,m_2$ $\hat{\Yim}_t^s$-monomials, we have in $\hat{\Yim}_t^s$:
$$m_1m_2=t^{2d_1(m_1,m_2)-2d_2(m_2,m_1)}m_2m_1=t^{2d_2(m_1,m_2)-2d_1(m_2,m_1)}m_2m_1$$
\end{prop}

\demo First we check that $m_1m_2=t^{2d_1(m_1,m_2)-2d_2(m_2,m_1)}m_2m_1$ on generators:
$$2d_1(\tilde{A}^{-1}_{i,l},\tilde{A}^{-1}_{i,k})-2d_2(\tilde{A}^{-1}_{i,k},\tilde{A}^{-1}_{i,l})=2u_{i,l-r_i}(\tilde{A}^{-1}_{i,k})-2u_{i,l+r_i}(\tilde{A}^{-1}_{i,k})=2(\delta_{l-k,-2r_i}-\delta_{l-k,2r_i})=\alpha(i,l,i,k)$$
$$2d_1(\tilde{A}^{-1}_{i,l},\tilde{A}^{-1}_{j,k})-2d_2(\tilde{A}^{-1}_{j,k},\tilde{A}^{-1}_{i,l})=2u_{i,l-r_i}(\tilde{A}^{-1}_{j,k})-2u_{i,l+r_i}(\tilde{A}^{-1}_{j,k})=\alpha(i,l,j,k)$$
$$2d_1(\tilde{A}^{-1}_{i,l},\tilde{Y}_{j,k})-2d_2(\tilde{Y}_{j,k},\tilde{A}^{-1}_{i,l})=2u_{i,l-r_i}(\tilde{Y}_{j,k})-2u_{i,l+r_i}(\tilde{Y}_{j,k})=-\beta(i,l,j,k)$$
$$2d_1(\tilde{Y}_{i,l},\tilde{A}^{-1}_{j,k})-2d_2(\tilde{A}^{-1}_{j,k},\tilde{Y}_{i,l})=2v_{i,l-r_i}(\tilde{A}^{-1}_{j,k})-2v_{i,l+r_i}(\tilde{A}^{-1}_{j,k})=-\beta(i,l,j,k)=\beta_{j,k,i,l}$$
The other equality $m_1m_2=t^{2d_2(m_1,m_2)-2d_1(m_2,m_1)}m_2m_1$ is checked in the same way.\qed

\noindent If $B(z)$ is not symmetric, the product is defined in section \ref{bznonsym}.

\subsubsection{Involution}\label{invol} We consider the $\ZZ[t^{\pm}]$-antilinear antimultiplicative involution of $\hat{\Yim}_t^s$ such that $\overline{\tilde{Y}_{i,l}}=\tilde{Y}_{i,l}$, $\overline{\tilde{A}_{i,l}^{-1}}=\tilde{A}_{i,l}^{-1}$, $\overline{t}=t^{-1}$.

\noindent In \cite{her} we gave a ``concrete'' construction of this involution for the generic case: in $\hat{\Yim}_u$ the involution is defined by $c_m\rightarrow -c_m$.

\begin{lem} There is a $\ZZ[t^{\pm}]$-basis $\overline{A}^s$\label{linea} of $\hat{\Yim}_t^s$ such that all $m\in \overline{A}^s$ is a $\hat{\Yim}_t^s$-monomial and:
$$\overline{m}=t^{2d_1(m,m)}m=t^{2d_2(m,m)}m$$
Moreover for $m_1,m_2\in\overline{A}^s$ we have $m_1m_2t^{-d_1(m_1,m_2)-d_2(m_1,m_2)}\in\overline{A}^s$.
\end{lem}

\demo For the first point it suffices to show that for $m$ a $\hat{\Yim}_t^s$-monomial there is a unique $\alpha\in\ZZ$ such that $\overline{t^{\alpha}m}=t^{2d_1(m,m)+\alpha}m$, that is to say for $m$ a $\hat{\Yim}_t^s$-monomial we have $m\overline{m}^{-1}\in t^{2\ZZ}$. This is proved as in lemma 6.12 of \cite{her}.

\noindent For the second point we compute:

\noindent $\overline{t^{-d_1(m_1,m_2)-d_2(m_1,m_2)}m_1m_2}
=t^{d_1(m_1,m_2)+d_2(m_1,m_2)}\overline{m_2}\overline{m_1}
\\=t^{2d_1(m_2,m_2)+2d_1(m_1,m_1)+d_1(m_1,m_2)+d_2(m_1,m_2)}m_2m_1
\\=t^{2d_1(m_2,m_2)+2d_1(m_1,m_1)+2d_1(m_2,m_1)+d_1(m_1,m_2)-d_2(m_1,m_2)}m_1m_2
\\=t^{2d_1(m_1m_2,m_1m_2)}(t^{-d_1(m_1,m_2)-d_2(m_1,m_2)}m_1m_2)$\qed

\noindent For example we have $\tilde{Y}_{i,l}\in\overline{A}^s$ (because $d_1(\tilde{Y}_{i,l},\tilde{Y}_{i,l})=0$) and if $s=0$ or $s>2r_i$ we have $t\tilde{A}_{i,l}^{-1}\in\overline{A}^s$ (because $d_1(\tilde{A}_{i,l}^{-1},\tilde{A}_{i,l}^{-1})=-1$).

\noindent For $m_1,m_2\in\overline{A}^s$ we set $m_1.m_2=m_1m_2t^{-d_1(m_1,m_2)-d_2(m_1,m_2)}\in\overline{A}^s$. We have $m_1.m_2=m_2.m_1$. The non commutative multiplication can be defined from $.$ by setting ($m_1,m_2\in\overline{A}^s$):
$$m_1m_2=t^{d_1(m_1,m_2)+d_2(m_1,m_2)}m_1.m_2$$
In the $ADE$-case it is the point of view adopted in \cite{Nab}. In particular if $s=0$ or $s>2r_i$, $\tilde{Y}_{i,l}$ (resp. $\tilde{A}_{i,l}^{-1}$) is denoted by $W_{i,l}$ (resp. $t^{-1}V_{i,l}$) in \cite{Nab}.

\noindent Let $\overline{A}=\overline{A}^0$ and for $s\geq 0$ there is a surjective map $p_s:\overline{A}\rightarrow \overline{A}^s$\label{tau} such that for $m\in\overline{A}$, $p_s(m)$ is the unique element of $\overline{A}^s$ such that for $i\in I,l\in\ZZ/s\ZZ$: 
$$y_{i,l}(p_s(m))=\underset{l'\in\ZZ/[l']=l}{\sum}y_{i,l}(m)\text{ , }v_{i,l}(p_s(m))=\underset{l'\in\ZZ/[l']=l}{\sum}v_{i,l}(m)$$
In particular it gives a $\ZZ[t^{\pm}]$-linear map $p_s:\hat{\Yim}_t\rightarrow \hat{\Yim}_t^s$.

\subsubsection{Notations and technical complements}\label{hatpi} A $\hat{\Yim}_t^s$-monomial is said to be $i$-dominant (resp. 
\\$i$-antidominant) if $\forall l\in\ZZ/s\ZZ$, $u_{i,l}(m)\geq 0$ (resp. $u_{i,l}(m)\leq 0$). We denote by $\overline{B}_i^s$ the set of $i$-dominant monomials $m$ such that $m\in\overline{A}^s$. 

\noindent A $\hat{\Yim}_t^s$-monomial is said to be dominant (resp. antidominant) if $\forall l\in\ZZ/s\ZZ,\forall i\in I$, $u_{i,l}(m)\geq 0$ (resp. $u_{i,l}(m)\leq 0$). We denote by $\overline{B}^s$ the set of dominant monomials $m$ such that $m\in\overline{A}^s$. In the generic case let $\overline{A}=\overline{A}^0$\label{hata}, $\overline{B}_i=\overline{B}_i^0$, $\overline{B}=\overline{B}^0$.\label{hatas}

\noindent We denote by $\hat{A}^s=\{m=\underset{i\in I,l\in\ZZ/s\ZZ}{\prod}Y_{i,l}^{y_{i,l}}A_{i,l}^{-v_{i,l}}/u_{i,l},v_{i,l}\geq 0\}\subset\hat{\Yim}^s$ the set of $\hat{\Yim}^s$-monomials. It is a $\ZZ$-basis of $\hat{\Yim}^s$ and $\pi_+(\overline{A}^s)=\hat{A}^s$. Let $\hat{B}_i^s=\{m\in A^s/\forall l\in\ZZ/s\ZZ,u_{i,l}(m)\geq 0\}=\pi_+(\overline{B}_i^s)$, $\hat{B}^s=\underset{i\in I}{\bigcap} \hat{B}_i^s=\pi_+(\overline{B}^s)$.

\noindent We define $\hat{\Pi}:\hat{\Yim}_t^s\rightarrow \Yim^s=\ZZ[Y_{i,l}^{\pm}]_{i\in I,l\in\ZZ/s\ZZ}$ as the ring morphism such that for $m$ a $\hat{\Yim}_t$-monomial $\hat{\Pi}(m)=\underset{i\in I,l\in\ZZ/s\ZZ}{\prod}Y_{i,l}^{u_{i,l}(m)}$ ($\Yim^s$ is commutative\label{pig}).

\noindent In particular for $i\in I$, $l\in\ZZ/s\ZZ$, we have:
$$\hat{\Pi}(\tilde{A}_{i,l}^{-1})=Y_{i,l-r_i}^{-1}Y_{i,l+r_i}^{-1}\underset{j/C_{j,i}<0}{\prod}\underset{k=C_{j,i}+1,C_{j,i}+3,...,-C_{j,i}-1}{\prod}Y_{j,l+2k}$$\label{ail}
\noindent and we denote this term by $A_{i,l}^{-1}=\hat{\Pi}(\tilde{A}_{i,l}^{-1})$. Let $A^s=\{m=\underset{i\in I,l\in\ZZ/s\ZZ}{\prod}Y_{i,l}^{u_{i,l}(m)}/u_{i,l}(m)\in\ZZ\}=\hat{\Pi}(\overline{A}^s)$\label{as} the set of $\Yim^s$-monomials, $B_i^s=\{m\in A^s/\forall l\in\ZZ/s\ZZ,u_{i,l}(m)\geq 0\}=\hat{\Pi}(\overline{B}_i^s)$, $B^s=\underset{i\in I}{\bigcap} B_i^s=\hat{\Pi}(\overline{B}^s)$.

\noindent If $q$ is generic then for $M\in\overline{A}$ and $m\in A$ there at most one $m'\in\overline{A}^s$ of the form $m'=t^{\alpha}M\tilde{A}_{i_1,l_1}^{-1}...\tilde{A}_{i_K,l_K}^{-1}$ such that $\hat{\Pi}(m')=m$ (the $A_{i,l}^{-1}$ are algebraically independant because we have supposed $\text{det}(C(z))\neq 0$ , see \cite{her}).

\noindent If $q$ is a root of unity the situation can be different: for example we suppose that $C$ is of type $A_2^{(1)}$ and $s=3$ (so $\text{det}(C(q))=q^{-3}(q^3-1)^2=0$). Then for all $L\geq 0$, we have:
$$\hat{\Pi}(\tilde{Y}_{1,0}\tilde{A}_{1,1}^{-L}\tilde{A}_{2,2}^{-L}\tilde{A}_{3,3}^{-L})=Y_{1,0}$$
and $\hat{\Pi}^{-1}(Y_{1,0})$ is infinite.

\noindent If $C$ is finite the situation is better. We have a generalization of lemma 3.14 of \cite{her} at roots of unity:

\begin{lem}\label{dommons} We suppose that $C$ is finite and that $s\geq 1$. Let $M$ be in $\overline{A}^s$. Then:

i) There is at most a finite number of $m'\in\overline{A}^s$ of the form $m'=t^{\alpha}M\tilde{A}_{i_1,l_1}^{-1}...\tilde{A}_{i_K,l_K}^{-1}$ such that $m'$ is dominant.

ii) For $m\in A^s$ there is at most a finite number of $m'\in\overline{A}^s$ of the form $m'=t^{\alpha}M\tilde{A}_{i_1,l_1}^{-1}...\tilde{A}_{i_K,l_K}^{-1}$ such that $\hat{\Pi}(m')=m$. \end{lem}

\demo First let us show (i) : let $m'$ be in $\overline{A}^s$ with $m'=t^{\alpha}M\underset{i\in I,l\in\ZZ/s\ZZ}{\prod}A_{i,l}^{-v_{i,l}}$ and the $v_{i,l}\geq 0$. It suffices to show that the condition $m'$ dominant implies that the $v_i=\underset{l'\in\ZZ/s\ZZ}{\sum}v_{i,l}$ are bounded (because $\ZZ/s\ZZ$ is finite). This condition implies :
$$u_i(m')=-2v_i+\underset{j\neq i}{\sum}(-C_{i,j}v_j) +u_i(M)\geq 0$$
Let $U$ be the column vector with coefficients $(u_1(M),...,u_n(M))$ and $V$ the column vector with coefficients $(v_1,...,v_n)$. So we have $U-CV\geq 0$. As $C$ is finite, the theorem 4.3 of \cite{kac} implies that $C^{-1}U-V\geq 0$ and so the $v_i$ are bounded.

\noindent For the $(ii)$ we use the same proof with the condition :
$$u_i(m')=-2v_i+\underset{j\neq i}{\sum}(-C_{i,j}v_j) +u_i(M)=u_i(m)$$\qed

\noindent In some cases we have another result. For $i\in I$ let $L_i=(C_{i,1},...,C_{i,n})$. 

\begin{lem}\label{affroot} We suppose that $s\geq 1$ and that there are $(\alpha_i)_{i\in I}\in\ZZ^I$ such that $\alpha_i>0$ and:
$$\underset{j\in I}{\sum}\alpha_j L_j=0$$
Then for $M\in A^s$ there are at most a finite number of dominant monomials $m\in B^s$ of the form $m=MA_{i_1,l_1}^{-1}A_{i_2,l_3}^{-1}...A_{i_k,l_k}^{-1}$.\end{lem}

\noindent In particular an affine Cartan matrix verifies the property of the lemma (see \cite{kac} for the coefficients $\alpha_j$).

\demo Consider $m'=\underset{i\in I,l\in\ZZ/s\ZZ}{\prod}A_{i,l}^{-v_{i,l}}$ and $m=Mm'$. For $i\in I$ let $v_i=\underset{l\in\ZZ/s\ZZ}{\sum}v_{i,l}\geq 0$. We have:
$$\underset{i\in I}{\sum}\alpha_i u_i(m')=\underset{i\in I}{\sum}\alpha_i\underset{j\in I}{\sum}(-C_{i,j})v_j=-\underset{j\in I}{\sum}v_j(\underset{i\in I}{\sum}\alpha_iC_{i,j})=0$$
We suppose that $m$ is dominant, in particular $u_{i,l}(m')\geq -u_{i,l}(M)$. So:
$$u_{i,l}(m')=u_i(m')-\underset{l'\in\ZZ/s\ZZ, l'\neq l}{\sum}u_{i,l'}(m')
\leq u_i(m')+\underset{l'\in\ZZ/s\ZZ, l'\neq l}{\sum}u_{i,l'}(M)$$
$$\leq \frac{1}{\alpha_i}\underset{j\neq i}{\sum}\alpha_{j}(-u_j(m')) +\underset{l'\in\ZZ/s\ZZ, l'\neq l}{\sum}u_{i,l'}(M)
\leq \frac{1}{\alpha_i}\underset{j\neq i}{\sum}\alpha_{j}u_j(M) +\underset{l'\in\ZZ/s\ZZ, l'\neq l}{\sum}u_{i,l'}(M)$$

\noindent So the $u_{i,l}(m')$ ($i\in I,l\in\ZZ/s\ZZ$) are bounded and there is at most a finite number of $m'$ such that $m$ is dominant.\qed

\section{$q,t$-characters in the generic case}\label{quatre}

In \cite{her} we defined $q,t$-characters for all finite Cartan matrices in the generic case. In this section we define $q$ and $q,t$-characters for all symmetrizable (non necessarily finite) Cartan matrix such that $i\neq j\Rightarrow C_{i,j}C_{j,i}\leq 3$, in particular for Cartan matrices of affine type (except $A_1^{(1)}$, $A_2^{(2)}$). We suppose $s=0$, that is to say $q$ is generic. The root of unity case will be studied in section \ref{rofcase}.

\subsection{Deformed screening operators}

Classical screening operators were introduced in \cite{Fre} and $t$-deformed screening operators were introduced in \cite{her01} for $C$ finite. We define and study deformed screening operators in the general case:

\begin{defi} $\hat{\Yim}_{i,u}$ is the $\hat{\Yim}_u$-bimodule defined by generators $\tilde{S}_{i,l}$ ($i\in I,l\in\ZZ$) and relations \label{tsil}:
$$\tilde{S}_{i,l}\tilde{A}_{j,k}^{-1}=t_{-C_{i,j}(z)(z^{(k-l)}+z^{(l-k)})}\tilde{A}_{j,k}^{-1}\tilde{S}_{i,l}$$
$$\tilde{S}_{i,l}\tilde{Y}_{j,k}=t_{\delta_{i,j}(z^{(k-l)}+z^{(l-k)})}\tilde{Y}_{j,k}\tilde{S}_{i,l}\text{ , }\tilde{S}_{i,l}t=t\tilde{S}_{i,l}$$
$$\tilde{S}_{i,l-r_i}-t_{-q^{-2r_i}-1}\tilde{A}_{i,l}^{-1}\tilde{S}_{i,l+r_i}=0$$
\end{defi}

\noindent In \cite{her} we made a concrete construction of $\hat{\Yim}_{i,u}$ by realizing it in $\hat{\mathcal{H}}_h$. Note that $\hat{\Yim}_t$ is a $\hat{\Yim}_u$-bimodule using the projection $\hat{\Yim}_u\rightarrow\hat{\Yim}_t$.

\begin{defi}\label{yit} $\hat{\Yim}_{i,t}$ is the $\hat{\Yim}_t$-bimodule $\hat{\Yim}_t\otimes_{\hat{\Yim}_u}\hat{\Yim}_{i,u}\otimes_{\hat{\Yim}_u}\hat{\Yim}_t$.\end{defi}

\noindent For $l\in\ZZ$ we denote by $\tilde{S}_{i,l}$ the image of $\tilde{S}_{i,l}$ in $\hat{\Yim}_{i,t}$. The $\hat{\Yim}_t$-module $\hat{\Yim}_{i,t}$ is torsion free. 

\noindent For $m$ a $\hat{\Yim}_t^s$-monomial only a finite number of $[\tilde{S}_{i,l},m]=(t^2-1)t^{u_{i,l}(m)-1}[u_{i,l}(m)]_t\tilde{S}_{i,l}\in (t^2-1)\hat{\Yim}_{i,t}$ are not equal to $0$, so we can define:

\begin{defi} The $i^{th}$-deformed screening operator is the map $S_{i,t}:\hat{\Yim}_t\rightarrow \hat{\Yim}_{i,t}$\label{tsit} defined by ($\lambda\in\hat{\Yim}_t$):
$$S_{i,t}(\lambda)=\frac{1}{t^2-1}\underset{l\in\ZZ}{\sum}[\tilde{S}_{i,l},\lambda]\in \hat{\Yim}_{i,t}$$
\end{defi}

\noindent Let $\hat{\mathfrak{K}}_{i,t}=\text{Ker}(S_{i,t})$\label{kit}. As $S_{i,t}$ is a derivation, $\hat{\mathfrak{K}}_{i,t}$ is a subalgebra of $\hat{\Yim}_t$.

\noindent At $t=1$ we define $S_i:\Yim\rightarrow \Yim_i=\underset{l\in\ZZ}{\bigoplus}\Yim S_{i,l}/\underset{l\in\ZZ}{\sum}\Yim.(S_{i,l-r_i}-A_{i,l}^{-1}S_{i,l+r_i})$\label{si} such for $m\in A$, $S_i(m)=m\underset{l\in\ZZ}{\sum}u_{i,l}(m)S_{i,l}$. It is the classical screening operator (see \cite{Fre}). For $m\in\hat{\Yim}_t$ we have $S_i(\hat{\Pi}(m))=\hat{\Pi}(S_{i,t}(m))$ where $\hat{\Pi}:\hat{\Yim}_{i,t}\rightarrow\Yim_i$ is defined by $\hat{\Pi}(m \tilde{S}_{i,l})=\hat{\Pi}(m)S_{i,l}$.

\noindent We set $\mathfrak{K}_i=\text{Ker}(S_i)$\label{ki} and $\mathfrak{K}=\underset{i\in I}{\bigcap}\mathfrak{K}_i$.

\noindent In the following a product $\overset{\rightarrow}{\underset{l\in\ZZ}{\prod}}M_l$ (resp. $\overset{\leftarrow}{\underset{l\in\ZZ}{\prod}}M_l$) is the ordered product $...M_{-2}M_{-1}M_0M_1...$ 
\\(resp. $...M_2M_1M_0M_{-1}...$).

\begin{defi}\label{teitm} For $M\in\hat{\Yim}_t$ a $i$-dominant monomial we define:
$$\overset{\leftarrow}{E}_{i,t}(M)=M(\underset{l\in\ZZ}{\prod}\tilde{Y}_{i,l}^{u_{i,l}(M)})^{-1}\overset{\leftarrow}{\underset{l\in\ZZ}{\prod}}(\tilde{Y}_{i,l}(1+t\tilde{A}_{i,l+r_i}^{-1}))^{u_{i,l}(M)}\in\hat{\Yim}_t$$
\end{defi}

\noindent For example we have $\overset{\leftarrow}{E}_{i,t}(\tilde{Y}_{i,l})=\tilde{Y}_{i,l}(1+t\tilde{A}_{i,l+r_i}^{-1})$, $\overset{\leftarrow}{E}_{i,t}(\tilde{A}_{i,l}^{-1}\tilde{Y}_{i,l+r_i}\tilde{Y}_{i,l-r_i})=\tilde{A}_{i,l}^{-1}\tilde{Y}_{i,l+r_i}\tilde{Y}_{i,l-r_i}$ and for $j\neq i$: $\overset{\leftarrow}{E}_{i,t}(\tilde{Y}_{j,l})=\tilde{Y}_{j,l}$.

\begin{thm}(\cite{her01})\label{kernelsoi} For all Cartan matrix $C$, the kernel $\hat{\mathfrak{K}}_{i,t}$ of $S_{i,t}$ is the $\ZZ[t^{\pm}]$-subalgebra of $\hat{\Yim}_t$ generated by the ($l\in\ZZ,j\neq i$):
$$\tilde{Y}_{i,l}(1+t\tilde{A}_{i,l+r_i}^{-1})\text{ , }\tilde{A}_{i,l}^{-1}\tilde{Y}_{i,l+r_i}\tilde{Y}_{i,l-r_i}\text{ , }\tilde{Y}_{j,l}\text{ , }\overset{\leftarrow}{E}_{i,t}(\tilde{A}_{j,l}^{-1})$$
For $M$ a $i$-dominant monomial we have $\overset{\leftarrow}{E}_{i,t}(M)\in \hat{\mathfrak{K}}_{i,t}$, and:
$$\hat{\mathfrak{K}}_{i,t}=\underset{M\in \overline{B}_i}{\bigoplus}\ZZ[t^{\pm}] \overset{\leftarrow}{E}_{i,t}(M)$$\end{thm}

\noindent Note that the proof of \cite{her01} works also if $C$ is not finite : the point of this proof is that an element $\chi\in\hat{\mathfrak{K}}_{i,t}-\{0\}$ has at least one $i$-dominant monomial, which is shown as in the $sl_2$-case. 

\noindent At $t=1$ it is a classical result of \cite{Fre}.

\noindent Note that in the $ADE$-case the identification (see section \ref{invol}) between the $t\tilde{A}_{i,l}^{-1}$ and the $V_{i,l}$ shows that the notation $\hat{\mathfrak{K}}_{i,t}$ coincides with the notation of \cite{Nab}.

\subsection{Reminder on the algorithm of Frenkel-Mukhin and on the deformed algorithm} 

\subsubsection{Completed algebras}\label{cpltdalg}

Let $\hat{\mathfrak{K}}_t=\underset{i\in I}{\bigcap}\text{Ker}(S_{i,t})\subset\hat{\Yim}_t$. It is a subalgebra of $\hat{\Yim}_t$.

\noindent We recall that a partial ordering is defined on the $\hat{\Yim}_t$-monomials by $m\in t^{\ZZ}\tilde{A}_{i,l}^{-1}m'\Leftrightarrow m< m'$. 

\noindent We define a $\NN$-graduation of $\hat{\Yim}_t$ by putting $\text{deg}(\tilde{A}_{i,l}^{-1})=1$, $\text{deg}(\tilde{Y}_{i,l})=0$. Note that $m<m'\Rightarrow \text{deg}(m)>\text{deg}(m')$.

\noindent We define the algebra $\hat{\Yim}_t^{\infty}\supset\hat{\Yim}_t$\label{tytinf} as the completion for this gradation. In particular the elements of $\hat{\Yim}_t^{\infty}$ are (infinite) sums $\underset{k\geq 0}{\sum}\lambda_k$ such that $\lambda_k$ is homogeneous of degree $k$.

\noindent In the same way we define $\hat{\mathfrak{K}}_{i,t}^{\infty}$\label{kitinf} such that $\hat{\Yim}_t^{\infty}\supset\hat{\mathfrak{K}}_{i,t}^{\infty}\supset\hat{\mathfrak{K}}_{i,t}$, that is to say $\chi\in\hat{\Yim}_t^{\infty}$ is in $\hat{\mathfrak{K}}_{i,t}^{\infty}$ if and only if it is of the form $\chi=\underset{k\geq 0}{\sum}\chi_k$ where:
$$\chi_k\in\underset{M\in \overline{B}_i/\text{deg}(M)=k}{\bigoplus}\ZZ[t^{\pm}]\overset{\leftarrow}{E}_{i,t}(M)$$

\noindent Let $\hat{\mathfrak{K}}_t^{\infty}=\underset{i\in I}{\bigcap}\hat{\mathfrak{K}}_{i,t}^{\infty}$. 

\noindent In the same way for $t=1$ we define $\hat{\Yim}^{\infty}\supset\hat{\Yim}$, $\Yim^{\infty}\supset\Yim$, $\hat{\mathfrak{K}}^{\infty}\supset\hat{\mathfrak{K}}$, $\mathfrak{K}^{\infty}\supset\mathfrak{K}$. They are well defined because in $\hat{\Yim}$ and in $\Yim$ the $A_{i,l}^{-1}$ are algebraically independent (see section \ref{hatpi}) and $\pi_+$ preserves the degree. In particular the maps $\pi_+$ and $\hat{\Pi}$ can be extended to maps $\pi_+:\hat{\Yim}_t^{\infty}\rightarrow \hat{\Yim}^{\infty}$ and $\hat{\Pi}:\hat{\Yim}_t^{\infty}\rightarrow\Yim^{\infty}$.

\noindent For $m$ a $\hat{\Yim}_t$-monomial let $u(m)=\text{max}\{l\in\ZZ/\forall k<l,\forall i\in I,u_{i,k}(m)=0\}$. We define the subset $C(m)\subset \overline{A}$\label{cpm}
$$C(m)=\{t^{\ZZ}m\tilde{A}_{i_1,l_1}^{-1}...\tilde{A}_{i_N,l_N}^{-1}/N\geq 0, l_1,...,l_N\geq u(m)\}\cap\overline{A}$$
We define the $\ZZ[t^{\pm}]$-submodule of $\hat{\Yim}_t^{\infty}$:
$$\tilde{C}(m)=\{\chi\in\hat{\Yim}_t^{\infty}/\chi=\underset{m'\in C(m)}{\sum}\lambda_{m'}(t)m'\}$$

\begin{lem}\label{ad} An element of $\hat{\mathfrak{K}}_t^{\infty}-\{0\}$ has at least one dominant monomial. An element of $\hat{\mathfrak{K}}_t-\{0\}$ has at least one dominant monomial and one antidominant monomial.\end{lem}

\demo For $\chi\in \hat{\mathfrak{K}}_t^{\infty}$ let $M$ be a maximal monomial of $\chi$. Then in the decomposition $\chi=\underset{k\geq 0}{\sum}\chi_k^{(i)}$ where $\chi_k^{(i)}\in\underset{M\in\overline{B}_i/\text{deg}(M)=k}{\bigoplus}\ZZ[t^{\pm}]\overset{\leftarrow}{E}_{i,t}(M)$ we see that $M$ is $i$-dominant.

\noindent For $\chi\in\hat{\mathfrak{K}}_t$, we can consider a maximal and a minimal monomial, and so we have a dominant monomial and an antidominant monomial in $\chi$.\qed

\subsubsection{Algorithms}\label{rapalgo}

In \cite{her} we defined a deformed algorithm to compute $q,t$-characters for $C$ finite. We had to show that this algorithm is well defined, that is to say that at each step the different ways to compute each term give the same result.

\noindent The formulas of \cite{her} gives also a (non necessarily well defined) deformed algorithm for all Cartan matrices, that is to say:

Let $m\in \overline{B}$. If the deformed algorithm beginning with $m$ is well defined, it gives an element $\hat{F}_t(m)\in\hat{\mathfrak{K}}_t^{\infty}$\label{tftm} such that $m$ is the unique dominant monomial of $\hat{F}_t(m)$.

\noindent An algorithm was also used by Nakajima in the $ADE$-case in \cite{Naa}. If we set $t=1$ and apply $\hat{\Pi}$ (where $\hat{\Pi}$ is defined in section \ref{hatpi}) we get a classical algorithm (it is analogous to the algorithm constructed by Frenkel and Mukhin in \cite{Fre2}). So:

Let $m\in B$. If the classical algorithm beginning with $m$ is well defined, it gives an element $F(m)\in\mathfrak{K}^{\infty}$\label{fm} such that $m$ is the unique dominant monomial of $F(m)$.

\noindent We say that the classical algorithm (resp. the deformed algorithm) is well defined if for all $m\in B$ (resp. all $m\in \overline{B}$) the classical algorithm (resp. deformed algorithm) beginning with $m$ is well defined.

\begin{lem}\label{clawd} If the deformed algorithm is well defined then the classical algorithm is well defined.\end{lem}

\demo If the deformed algorithm beginning with $m$ is well defined then the classical algorithm beginning with $\hat{\Pi}(m)$ is well defined and $F(\hat{\Pi}(m))=\hat{\Pi}(\hat{F}_t(m))$.\qed

\noindent The following results are known:

If $C$ is finite then the classical algorithm is well defined (\cite{Fre2}).

If  $C$ is finite and symmetric then the deformed algorithm is well defined (\cite{Nab}).

If $C$ is finite then the deformed algorithm is well defined (\cite{her}).

\noindent In this section (theorem \ref{geneth}) we show that the classical and the deformed algorithms are well defined for a (non necessarily finite) Cartan matrix such that $i\neq j\Rightarrow C_{i,j}C_{j,i}\leq 3$.

\subsection{Morphism of $q,t$-characters} The construction of \cite{her} is based on the fact that we can compute explicitly $q,t$-characters for the submatrices of format $2$ of the Cartan matrix. So:

\subsubsection{The case $n=2$}

\begin{prop}\label{deuxgene} We suppose that $C$ is a Cartan matrix of rank $2$. The following properties are equivalent:

i) For all $m\in B$, $F(m)\in\mathfrak{K}$

ii) $C$ is finite

iii) $C_{1,2}C_{2,1}\leq 3$

iv) For $i=1$ or $2$, $\hat{\mathfrak{K}}_t\cap \tilde{C}(\tilde{Y}_{i,0})\neq \{0\}$ 

v) For $i=1$ or $2$, $C(Y_{i,0})$ has an antidominant monomial
\end{prop}

\demo The Cartan matrices of rank $2$ such that $C_{1,2}C_{2,1}\leq 3$ are matrices of type $A_1\times A_1$, $A_2$, $B_2$, $C_2$, $G_2$ or $G_2^{t}$. Those are finite Cartan matrices of rank 2, so $(ii)\Leftrightarrow (iii)$. Moreover if $C$ is finite, the classical theory of $q$-characters shows $(ii)\Rightarrow (i)$.

\noindent We have seen in \cite{her} that $(ii)\Rightarrow (iv)$. It follows from lemma \ref{ad} that $(iv)\Rightarrow (v)$ and $(i)\Rightarrow (v)$.  

\noindent So it suffices to show that $(v)\Rightarrow (iii)$. We suppose there is an antidominant monomial $m\in C(Y_{1,0})$. We can suppose $C_{1,2}<0$ and $C_{2,1}<0$. $m$ verifies $\hat{\Pi}(m)=Y_{1,0}A_{1,l_1}^{-1}...A_{1,l_L}^{-1}A_{2,l_1}^{-1}...A_{2,l_M}^{-1}$ where $L,M\geq 0$. In particular we have:
$$u_1(m)=1-2L-MC_{1,2}\text{ and }u_2(m)=-2M-LC_{2,1}$$
As $m$ is antidominant, we have $u_1(m),u_2(m)\leq 0$.

if $M=0$, we have $u_2(m)=-LC_{2,1}\leq 0\Rightarrow L=0$ and $u_1(m)=1>0$, impossible.

if $M>0$, we have:
$$\frac{L}{M}>\frac{-C_{1,2}}{2}\text{ and }\frac{L}{M}\leq\frac{2}{-C_{2,1}}$$
$$\frac{-C_{1,2}}{2}<\frac{2}{-C_{2,1}}\Rightarrow C_{1,2}C_{2,1}\leq 3$$
\qed

\subsubsection{General case}

\begin{thm}\label{geneth} If $i\neq j\Rightarrow C_{i,j}C_{j,i}\leq 3$,  then the classical and the deformed algorithms are well defined.\end{thm}

\demo It suffices to show that the deformed algorithm is well defined (lemma \ref{clawd}).We follow the proof of theorem 5.13 of \cite{her}: it suffices to construct $\hat{F}_t(m)$ for $m=\tilde{Y}_{i,0}$ ($i\in I$) and it suffices to see the property for the matrices $\begin{pmatrix}2&C_{i,j}\\C_{j,i} & 2\end{pmatrix}$. If $r_i\wedge r_j=1$ this follows from proposition \ref{deuxgene}. If $r_i\wedge r_j>1$ it suffices to replace $r_i,r_j$ with $\frac{r_i}{r_i\wedge r_j},\frac{r_j}{r_i\wedge r_j}$ (in fact it means that we replace $q$ by $q^{r_i\wedge r_j}$).\qed

\noindent In the following we suppose that $i\neq j\Rightarrow C_{i,j}C_{j,i}\leq 3$. For example $C$ could be of finite or affine type (except $A_1^{(1)}$, $A_2^{(2)}$).

\noindent We conjecture that for $C$ of type $A_1^{(1)}$ (with $r_1=r_2=2$) and of type $A_2^{(2)}$ the algorithms are well defined. This conjecture is motivated by the remarks of the introduction about representation theory of quantum affinization algebras (note that for $C$ of type $A_1^{(1)}$ and $r_1=r_2=1$ the classical algorithm is not well defined).

\subsubsection{Definition of $\chi_{q,t}$} We verify as in \cite{her} that $\hat{F}_t(\tilde{Y}_{i,l})\hat{F}_t(\tilde{Y}_{j,l})=\hat{F}_t(\tilde{Y}_{j,l})\hat{F}_t(\tilde{Y}_{i,l})$. Let $\text{Rep}=\ZZ[X_{i,l}]_{i\in I,l\in\ZZ}$ as in section \ref{qaa} and a $\text{Rep}$-monomials is a product of the $X_{i,l}$.

\begin{defi}\label{chiqt} The morphism of $q,t$-characters $\chi_{q,t}:\text{Rep}\rightarrow\hat{\mathfrak{K}}^{\infty}_t$ is the $\ZZ$-linear map such that:
$$\chi_{q,t}(\underset{i\in I,l\in\ZZ}{\prod}X_{i,l}^{x_{i,l}})=\overset{\rightarrow}{\underset{l\in\ZZ}{\prod}}\underset{i\in I}{\prod}\hat{F}_t(\tilde{Y}_{i,l})^{x_{i,l}}$$
The morphism of $q$-characters $\chi_q:\text{Rep}\rightarrow \mathfrak{K}^{\infty}$ is defined by $\chi_q=\hat{\Pi}\circ\hat{\chi}_{q,t}$.\end{defi}

\begin{thm}\label{axiomes}(\cite{her}) The $\ZZ$-linear map $\chi_{q,t}:\text{Rep}\rightarrow\hat{\Yim}_t^{\infty}$ is injective and is characterized by the three following properties:

1) For $M$ a $\text{Rep}$-monomial define $m=\underset{i\in I,l\in\ZZ}{\prod}\tilde{Y}_{i,l}^{x_{i,l}(M)}\in \overline{B}$. Then we have :
$$\chi_{q,t}(M)=m+\underset{m'<m}{\sum}a_{m'}(t)m'\text{ (where $a_{m'}(t)\in\ZZ[t^{\pm}]$)}$$

2) The image $\text{Im}(\chi_{q,t})$ is contained in $\hat{\mathfrak{K}}_t^{\infty}$.

3) Let $M_1,M_2$ be $\text{Rep}$-monomials such that $\text{max}\{l/\underset{i\in I}{\sum}x_{i,l}(M_1)>0\}\leq \text{min}\{l/\underset{i\in I}{\sum}x_{i,l}(M_2)>0\}$. We have :
$$\chi_{q,t}(M_1M_2)=\chi_{q,t}(M_1)\chi_{q,t}(M_2)$$
\end{thm}

\noindent Those properties are generalizations of Nakajima's axioms \cite{Nab} for $q$ generic, so:

\begin{cor} If $C$ is finite then we have $\pi_+(\text{Im}(\chi_{q,t}))\subset\hat{\Yim}$ and $\chi_q:\text{Rep}\rightarrow \Yim$ is the classical morphism of $q$-characters and $\chi_{q,t}$ is the morphism of \cite{her}. In particular if $C$ is of type $ADE$ then $\chi_{q,t}$ is the morphism of $q,t$-characters of \cite{Nab}.\end{cor}

\section{$\epsilon,t$-characters in the root of unity case}\label{rofcase}

In this section we define and study $\epsilon,t$-characters at roots of unity: let $\epsilon\in\CC^*$ be a $s^{th}$-primitive root of unity. We suppose that $s>2r^{\vee}$.

\noindent The case $t=1$ was study in \cite{Fre3} (but classical screening operators in the root of unity case were not defined). The $t$-deformations were studied in the $ADE$-case by Nakajima in \cite{Nab} using quiver varieties. In this section we suppose that $i\neq j\Rightarrow C_{i,j}C_{j,i}\leq 3$ and $B(z)$ is symmetric. In particular $C$ can be of finite type or of affine type (except $A_1^{(1)}$, $A_{2l}^{(2)}$, $l\geq 1$, see section \ref{andsym}). The deformed algorithm is well defined and $\chi_{q,t}$ exists (theorem \ref{geneth}).

\subsection{Reminder: classical $\epsilon$-characters at roots of unity} We define $\tau_s:\Yim\rightarrow\Yim_s$\label{taus} as the ring homomorphism such that $\tau_s(Y_{i,l})=Y_{i,[l]}$ where for $l\in\ZZ$ we denote by $[l]$ its image in $\ZZ/s\ZZ$.

\noindent If $C$ is finite the morphism of $\epsilon$-characters $\chi_{\epsilon}:\text{Rep}^s\rightarrow \Yim^s$ is defined by Frenkel and Mukhin (see section \ref{qaa}). We have the following characterization:

\begin{thm}\label{tunroot}(\cite{Fre3}) If $C$ is finite, the morphism of $\epsilon$-characters $\chi_{\epsilon}:\text{Rep}^s\rightarrow\Yim^s$ verifies ($l_0\in\ZZ$):
$$\chi_{\epsilon}(\underset{i\in I,l\in\ZZ/s\ZZ}{\prod}X_{i,l}^{x_{i,l}})=\tau_s(\chi_q(\underset{i\in I, l_0\leq l\leq l_0+s-1}{\prod}X_{i,l}^{x_{i,[l]}}))$$
\end{thm}

\noindent Note that this formula suffices to characterize the $\ZZ$-linear map $\chi_{\epsilon}$.

\noindent If $C$ is not finite, we can consider $\hat{\Yim}^s=\ZZ[Y_{i,l},A_{i,l}^{-1}]_{i\in I,l\in\ZZ/s\ZZ}$ and the completion $\hat{\Yim}_t^{s,\infty}$ as in the generic case. We define $\hat{\chi}_{\epsilon}:\text{Rep}^s\rightarrow\hat{\Yim}^{s, \infty}$ with the formula of the theorem \ref{tunroot}. The map $\hat{\chi}_{\epsilon}$ is also an injective ring homomorphism. 

\noindent In the following we give an analogous construction in the deformed case $t\neq 1$. 

\subsection{Construction of $\chi_{\epsilon,t}$} The point for the $t$-deformation is that we can not define a natural $t$-analog of $\tau_s$ which is a ring homomorphism. In this section we construct an analog $\tau_{s,t}$ of $\tau_s$ which is not a ring homomorphism but has nice properties.

\subsubsection{Definition of $\tau_{s,t}$}\label{deficonst} First let us briefly explain how $\tau_{s,t}$ is constructed. The main property is a compatibility with some ordered products: suppose that $l_1> l_2$ ($l_1,l_2\in\ZZ$), that $m_1\in\hat{\Yim}_t$ involves only the $\tilde{Y}_{i,l_1}, \tilde{A}_{i,l_1}^{-1}$ and that $m_2$ involves only the $\tilde{Y}_{i,l_2}, \tilde{A}_{i,l_2}^{-1}$. Then $\tau_{s,t}$ is defined such that $\tau_{s,t}(m_1m_2)=\tau_{s,t}(m_1)\tau_{s,t}(m_2)$. Let us now write it in a formal way:

\noindent For $m$ a $\hat{\Yim}_t$-monomial and $l\in\ZZ$, let :
$$\pi_l(m)=(\underset{i\in I}{\prod}\tilde{Y}_{i,l}^{y_{i,l}(m)})(\underset{i\in I}{\prod}\tilde{A}_{i,l}^{-v_{i,l}(m)})$$
It is well defined because for $i,j\in I$ and $l\in\ZZ$ we have $\tilde{Y}_{i,l}\tilde{Y}_{j,l}=\tilde{Y}_{j,l}\tilde{Y}_{i,l}$, $\tilde{A}_{i,l}^{-1}\tilde{A}_{j,l}^{-1}=\tilde{A}_{j,l}^{-1}\tilde{A}_{i,l}^{-1}$ and for $i\neq j$, $\tilde{A}_{i,l}^{-1}\tilde{Y}_{j,l}=\tilde{Y}_{j,l}\tilde{A}_{i,l}^{-1}$ (theorem \ref{dessus}). 

\noindent Let $\overset{\rightarrow}{m}={\underset{l\in\ZZ}{\overset{\rightarrow}{\prod}}}\pi_l(m)$\label{tm}, $\overset{\leftarrow}{m}={\underset{l\in\ZZ}{\overset{\leftarrow}{\prod}}}\pi_l(m)$\label{marrow}, and \label{aarrow}:
$$\overset{\rightarrow}{A}=\{\overset{\rightarrow}{m}/\text{ $m$ $\hat{\Yim}_t$-monomial}\}\text{ and }\overset{\leftarrow}{A}=\{\overset{\leftarrow}{m}/\text{ $m$ $\hat{\Yim}_t$-monomial}\}$$
It follows from theorem \ref{dessus} that $\overset{\rightarrow}{A}$ and $\overset{\leftarrow}{A}$ are $\ZZ[t^{\pm}]$-basis of $\hat{\Yim}_t$.

\begin{defi}\label{taust} We define $\tau_{s,t}:\hat{\Yim}_t\rightarrow\hat{\Yim}_t^s$ as the $\ZZ[t^{\pm}]$-linear map such that for $m\in\overset{\leftarrow}{A}$:
$$\tau_{s,t}(m)={\underset{l\in\ZZ}{\overset{\leftarrow}{\prod}}}(\underset{j\in I}{\prod}\tilde{A}_{j,[l]}^{-v_{j,l}(m)})(\underset{j\in I}{\prod}\tilde{Y}_{j,[l]}^{y_{j,l}(m)})$$\end{defi}

\noindent Note that $\tau_{s,t}$ is not a ring homomorphism and is not injective.

\subsubsection{Definition of $\chi_{\epsilon,t}$}

We define a $\NN$-gradation of $\hat{\Yim}_t^s$, the completed algebra $\hat{\Yim}_t^{s,\infty}$\label{tytinfs} in the same way as we did for the generic case (section \ref{cpltdalg}). In particular $\tau_{s,t}$ is compatible with the gradations of $\hat{\Yim}_t$ and $\hat{\Yim}_t^s$ and is extended to a map $\tau_{s,t}:\hat{\Yim}_t^{\infty}\rightarrow\hat{\Yim}_t^{s,\infty}$.

\begin{defi}\label{epst} The morphism of $q,t$-characters at the $s^{th}$-primitive roots of unity $\chi_{\epsilon,t}:\text{Rep}^s\rightarrow\hat{\Yim}_t^{s,\infty}$ is the $\ZZ$-linear map such that:
$$\chi_{\epsilon,t}(\underset{i\in I,l\in\ZZ/s\ZZ}{\prod}X_{i,l}^{x_{i,l}})=\tau_{s,t}(\chi_{q,t}(\underset{i\in I,0\leq l\leq s-1}{\prod}X_{i,l}^{x_{i,[l]}}))$$
\end{defi}

\begin{prop}\label{axioms} The morphism $\chi_{\epsilon,t}$ verifies the following properties:

1)  The following diagram is commutative:
$$\begin{array}{rcccl}
\text{Rep}^s&\stackrel{\chi_{\epsilon,t}}{\longrightarrow}&\text{Im}(\chi_{\epsilon,t})\\
\text{id}\downarrow &&\downarrow&\pi_+\\
\text{Rep}^s&\stackrel{\hat{\chi}_{\epsilon}}{\longrightarrow}&\hat{\Yim}^{s,\infty}
\end{array}$$ 

2) If $C$ is finite we have $\pi_+(\text{Im}(\chi_{\epsilon,t}))\subset\hat{\Yim}^s$ and $\hat{\Pi}\circ\chi_{\epsilon,t}=\chi_{\epsilon}$.

3) The map $\chi_{\epsilon,t}$ is injective.

4) For a $\text{Rep}$-monomial $M$ define $m=\underset{i\in I,l\in\ZZ/s\ZZ}{\prod}\tilde{Y}_{i,l}^{x_{i,l}(M)}\in \overline{B}^s$. Then we have :
$$\chi_{q,t}(M)=m+\underset{m'<m}{\sum}a_{m'}(t)m'\text{ (where $a_{m'}(t)\in\ZZ[t^{\pm}]$)}$$
\end{prop}

\demo 

1) Consequence of the definition and of $(\tau_{s,t})_{t=1}=\tau_s$.

2) Consequence of (1) and of theorem \ref{tunroot}.

3) Consequence of (1) and of the injectivity of $\hat{\chi}_{\epsilon}$ (see section \ref{qaa}).

4) Consequence of the analogous property of $\chi_{q,t}$ (1. of theorem \ref{axiomes}).\qed

\noindent Note that 2) means that in the finite case we get at $t=1$ the map of \cite{Fre2}.

\noindent In the following we show other fundamental properties of $\chi_{\epsilon,t}$ (theorem \ref{qtint} and theorem \ref{axquat}).

\subsection{Classical and deformed screening operators at roots of unity} We define classical and deformed screening operators at roots of unity in order to have an analog of the property 2 of theorem \ref{axiomes} at roots of unity.

\subsubsection{Deformed bimodules}

\begin{defi} $\hat{\Yim}_{i,u}^s$ is the $\hat{\Yim}_u^s$-bimodule defined by generators $\tilde{S}_{i,l}$ ($i\in I,l\in\ZZ/s\ZZ$) and relations :
$$\tilde{S}_{i,l}\tilde{A}_{j,k}^{-1}=t_{-C_{i,j}(z)(z^{(k-l)}+z^{(l-k)})}\tilde{A}_{j,k}^{-1}\tilde{S}_{i,l}\text{ , }\tilde{S}_{i,l}\tilde{Y}_{j,k}=t_{\delta_{i,j}(z^{(k-l)}+z^{(l-k)})}\tilde{Y}_{j,k}\tilde{S}_{i,l}$$
$$\tilde{S}_{i,l}t=t\tilde{S}_{i,l}\text{ , }\tilde{S}_{i,l-r_i}-t_{-q^{-2r_i}-1}\tilde{A}_{i,l}^{-1}\tilde{S}_{i,l+r_i}\text{ , }\tilde{S}_{i,l+s}-\tilde{S}_{i,l}$$
\end{defi}

\noindent Note that this structure is well-defined: if $s\geq 1$, for example we have $t_{-C_{i,j}(z)(z^{(k+s-l)}+z^{(l-k-s)})}=t_{-C_{i,j}(z)(z^{(k-l)}+z^{(l-k)})}$.

\noindent Note that $\hat{\Yim}_t^s$ is a $\hat{\Yim}_u^s$-bimodule using the projection $\hat{\Yim}_u^s\rightarrow\hat{\Yim}_t^s$.

\begin{defi}\label{yits} $\hat{\Yim}_{i,t}^s$ is the $\hat{\Yim}_t^s$-bimodule $\hat{\Yim}_t^s\otimes_{\hat{\Yim}_u^s}\hat{\Yim}_{i,u}^s\otimes_{\hat{\Yim}_u^s}\hat{\Yim}_t^s$.\end{defi}

\noindent For $l\in\ZZ/s\ZZ$ we denote by $\tilde{S}_{i,l}$ the image of $\tilde{S}_{i,l}$ in $\hat{\Yim}_{i,t}^s$. If $s\geq 1$, the $\hat{\Yim}_t^s$-module $\hat{\Yim}_{i,t}^s$ has torsion: 
$$\tilde{S}_0=t^{\alpha}\tilde{A}_{r_i}^{-1}\tilde{A}_{3r_i}^{-1}...\tilde{A}_{(2s-1)r_i}^{-1}\tilde{S}_{2r_is}=t^{\alpha}\tilde{A}_{r_i}^{-1}\tilde{A}_{3r_i}^{-1}...\tilde{A}_{(2s-1)r_i}^{-1}\tilde{S}_{0}$$
where $\alpha=-2s$ if $s|2r_i$ and $\alpha=-s$ otherwise.

\subsubsection{Deformed screening operators} As in the generic case, we can define:

\begin{defi} The $i^{th}$-deformed screening operator is the map $S_{i,t}^s:\hat{\Yim}_t^s\rightarrow \hat{\Yim}_{i,t}^s$\label{tsits} defined by ($\lambda\in\hat{\Yim}_t^s$):
$$S_{i,t}^s(\lambda)=\frac{1}{t^2-1}\underset{l\in\ZZ/s\ZZ}{\sum}[\tilde{S}_{i,l},\lambda]\in \hat{\Yim}_{i,t}^s$$
\end{defi}

\noindent We define $\hat{\mathfrak{K}}_{i,t}^s=\text{Ker}(S_{i,t}^s)$\label{kits} and we complete this algebra $\hat{\mathfrak{K}}_{i,t}^{s,\infty}\supset\hat{\mathfrak{K}}_{i,t}^s$\label{kitinfs}.

\subsubsection{Classical screening operators at roots of unity} We suppose in this section that $t=1$. 

\noindent The classical screening operators at roots of unity are
$$S_i^s:\Yim^s\rightarrow \Yim_i^s=\underset{l\in\ZZ/s\ZZ}{\bigoplus}\Yim^s S_{i,l}/\underset{l\in\ZZ/s\ZZ}{\sum}\Yim^s.(S_{i,l-r_i}-A_{i,l}^{-1}S_{i,l+r_i})$$\label{sis} 
such that for $m\in A^s$, $S_i^s(m)=m\underset{l\in\ZZ/s\ZZ}{\sum}u_{i,l}(m)S_{i,l}$. 

\noindent For $\lambda\in\hat{\Yim}_t^s$ we have $S_i(\hat{\Pi}(\lambda))=\hat{\Pi}(S_{i,t}(\lambda))$ where $\hat{\Pi}:\hat{\Yim}_{i,t}^s\rightarrow\Yim_i^s$ is defined by $\hat{\Pi}(m \tilde{S}_{i,l})=\hat{\Pi}(m)S_{i,l}$.

\noindent The map $\tau_s:\Yim\rightarrow \Yim^s$ is a ring homomorphism. In particular we can define a $\ZZ$-linear map $\tau_s:\Yim_i\rightarrow\Yim_i^s$ such that:
$$\tau_s(mS_{i,l})=\Pi_{s}(m)S_{i,[l]}$$
Indeed it suffices to see it agrees with the defining relations of $\Yim_i$:
$$\tau_s(mA_{i,l+r_i}^{-1}S_{i,l+2r_i})=\tau_s(mA_{i,l+r_i}^{-1})S_{i,[l+2r_i]}=\tau_s(m)A_{i,[l+r_i]}^{-1}S_{i,[l+2r_i]}=\tau_s(m)S_{i,[l]}=\tau_s(m S_{i,l})$$
\noindent Note that the crucial point is that $\tau_s$ is a ring homomorphism.

\begin{lem}\label{tisone} We have $\tau_s\circ S_i=S_i^s\circ\tau_s$.\end{lem}

\demo It suffices to see for $m$ a $\Yim$-monomial:
$$\tau_s(S_i(m))=\underset{l\in\ZZ}{\sum}u_{i,l}(m)\tau_s(mS_{i,l})
=\tau_s(m)\underset{l\in\ZZ}{\sum}u_{i,l}(m)S_{i,[l]}
=\tau_s(m)\underset{0\leq l\leq s-1}{\sum}(\underset{r\in\ZZ}{\sum}u_{i,l+rs}(m))S_{i,[l]}$$
$$=\tau_s(m)\underset{0\leq l\leq s-1}{\sum}u_{i,[l]}(\tau_s(m))S_{i,[l]}
=S_i^s(\tau_s(m))$$
\qed

\noindent For $m\in B_i^s$, we set $E_i(m)=m\underset{l\in\ZZ/s\ZZ}{\prod}(1+A_{i,l+r_i}^{-1})^{u_{i,l}(m)}$\label{eim}. Let $\mathfrak{K}_i^s=\text{Ker}(S_i^s)$\label{kis}. 

\begin{prop}\label{dommon} $\tau_s(\mathfrak{K}_i)$ is a subalgebra of $\mathfrak{K}_i^s$. Moreover:
$$\tau_s(\mathfrak{K}_i)=\underset{m\in B_s}{\bigoplus}E_i(m)$$
In particular if $\chi\in\tau_s(\mathfrak{K}_i)$ has no $i$-dominant monomial then $\chi=0$.\end{prop}

\demo The lemma \ref{tisone} gives $\tau_s(\mathfrak{K}_i)\subset\mathfrak{K}_i^s$ and $\tau_s$ is an algebra homomorphism. 

\noindent For $m\in\overline{B}$ we have $\tau_s(E_i(m))=E_i(\tau_s(m))$ and so it follows from theorem \ref{kernelsoi} that $\tau_s(\mathfrak{K}_i)=\underset{m\in B_s}{\bigoplus}E_i(m)$.\qed

\subsection{The image of $\chi_{\epsilon,t}$} In this section we show an analog of the property 2 of theorem \ref{axiomes} at roots of unity.

\begin{thm}\label{qtint} The image of $\chi_{\epsilon,t}$ is contained in $\hat{\mathfrak{K}}_t^{s, \infty}$.\end{thm}

\noindent With the help of theorem \ref{axiomes} it suffices to show that $\tau_{s,t}(\hat{\mathfrak{K}}_{i,t})\subset \hat{\mathfrak{K}}_{i,t}^s$ which will be done in proposition \ref{kernels}.

\subsubsection{The bicharacters $D_1$, $D_2$} For $m$ a $\hat{\Yim}_t$-monomial and $k\in\ZZ$ let:
$$m[k]=m(\overset{\leftarrow}{m})^{-1}{\underset{l\in\ZZ}{\overset{\leftarrow}{\prod}}}(\underset{j\in I}{\prod}\tilde{Y}_{j,l}^{y_{j,l+k}(m)})(\underset{j\in I}{\prod}\tilde{A}_{j,l}^{-v_{j,l+k}(m)})$$

\noindent Note that $\tau_{s,t}(m[ks])=\tau_{s,t}(m)$ and for $m\in\overset{\leftarrow}{A}$, $k\in\ZZ$ we have $m[k]\in\overset{\leftarrow}{A}$.

\noindent For $m_1,m_2$ $\hat{\Yim}_t$-monomials, and $k\in\ZZ$ we have :
$$d_1(m_1,m_2[k])=d_1(m_1[-k],m_2)\text{ and }d_2(m_1,m_2[k])=d_2(m_1[-k],m_2)$$
Moreover there is only a finite number of $k\in\ZZ$ such that $d_1(m_1,m_2[k])\neq 0$ or $d_2(m_1,m_2[k])\neq 0$. So we can define:

\begin{defi}\label{gd} For $m_1,m_2$ $\hat{\Yim}_t$-monomials we define:
$$D_1(m_1,m_2)=\underset{r\in\ZZ}{\sum}d_1(m_1,m_2[rs])=\underset{r\in\ZZ}{\sum}d_1(m_1[rs],m_2)$$
$$D_2(m_1,m_2)=\underset{r\in\ZZ}{\sum}d_2(m_1,m_2[rs])=\underset{r\in\ZZ}{\sum}d_2(m_1[rs],m_2)$$\end{defi}

\begin{lem}\label{dspis} For $m_1,m_2$ $\hat{\Yim}_t$-monomials we have:
$$D_1(m_1,m_2)=d_1(\tau_{s,t}(m_1),\tau_{s,t}(m_2))\text{ , }D_2(m_1,m_2)=d_2(\tau_{s,t}(m_1),\tau_{s,t}(m_2))$$
\end{lem}
\noindent In particular we have in $\hat{\Yim}_t^s$:
$$\tau_{s,t}(m_1)\tau_{s,t}(m_2)=t^{D_1(m_1,m_2)-D_2(m_2,m_1)}\tau_{s,t}(m_2)\tau_{s,t}(m_1)$$

\demo For example for $d_1$ we compute:

\noindent $d_1(\tau_{s,t}(m_1),\tau_{s,t}(m_2))
\\=\underset{i\in I,l\in\ZZ/s\ZZ}{\sum}v_{i,l+r_i}(\tau_{s,t}(m_1))u_{i,l}(\tau_{s,t}(m_2))+w_{i,l+r_i}(\tau_{s,t}(m_1))v_{i,l}(\tau_{s,t}(m_2))
\\=\underset{i\in I,0\leq l\leq s-1,r\in\ZZ,r'\in\ZZ}{\sum}v_{i,l+r_i+rs}(m_1)u_{i,l+r's}(m_2)+w_{i,l+r_i+rs}(m_1)v_{i,l+r's}(m_2)
\\=\underset{i\in I,l\in\ZZ,r\in\ZZ}{\sum}v_{i,l+r_i}(m_1)u_{i,l+rs}(m_2)+w_{i,l+r_i}(m_1)v_{i,l+rs}(m_2)
\\=\underset{r\in\ZZ}{\sum}d_1(m_1,m_2[rs])$\qed

\subsubsection{Technical lemmas}

\begin{lem}\label{mappis} Let $m$ be a $\hat{\Yim}_t$-monomial of the form $m=Z_1Z_2...Z_K$ where $Z_k=\tilde{Y}_{i_k,l_k}$ or $Z_k=\tilde{A}_{i_k,l_k}^{-1}$. We suppose that $k>k'$ implies $l_k\leq l_{k'}+r^{\vee}$ and $(Z_k,Z_{k'})\notin\{(\tilde{A}_{i,l}^{-1},\tilde{A}_{i,l'}^{-1})/i\in I,l'<l\}$. Then we have:
$$\tau_{s,t}(m)=\tau_{s,t}(Z_1)\tau_{s,t}(Z_2)...\tau_{s,t}(Z_K)$$
\end{lem}

\demo First we order the factors of $m$:
$$m=t^{2\underset{k<k'/l_{k}<l_{k'}}{\sum}d_1(Z_k,Z_{k'})-d_2(Z_{k'},Z_k)}\overset{\leftarrow}{m}$$
So we can apply $\tau_{s,t}$:
$$\tau_{s,t}(m)=t^{2\underset{k<k'/l_{k}<l_{k'}}{\sum}d_1(Z_k,Z_{k'})-d_2(Z_{k'},Z_k)}\tau_{s,t}(\overset{\leftarrow}{m})$$
where:
$$\tau_{s,t}(\overset{\leftarrow}{m})={\underset{l\in\ZZ}{\overset{\leftarrow}{\prod}}}(\underset{j\in I}{\prod}\tilde{Y}_{j,[l]}^{y_{j,l}(m)})(\underset{j\in I}{\prod}\tilde{A}_{j,[l]}^{-v_{j,l}(m)})$$
If we order the factors of $\tau_{s,t}(Z_1)\tau_{s,t}(Z_2)...\tau_{s,t}(Z_k)$, we get:
$$\tau_{s,t}(Z_1)\tau_{s,t}(Z_2)...\tau_{s,t}(Z_k)=t^{2\underset{k<k'/l_{k}<l_{k'}}{\sum}(D_1(Z_k,Z_{k'})-D_2(Z_{k'},Z_k))}\tau_{s,t}(\overset{\leftarrow}{m})$$
So it suffices to show that $k<k'$ and $l_{k}<l_{k'}$ implies $d_1(Z_k,Z_{k'})-d_2(Z_{k'},Z_k)=D_1(Z_k,Z_{k'})-D_2(Z_{k'},Z_k)$. But we have $0<l_{k'}-l_k\leq r^{\vee}$ and $s>2r^{\vee}$. So for $p\in\ZZ$ such that $p\neq 0$ we have $|l_k-l_{k'}+ps|>r^{\vee}$. But in general for $k_1,k_2$, we have:
$$[Z_{k_1},Z_{k_2}]\neq 0\Rightarrow (Z_{k_1},Z_{k_2})=(\tilde{A}_{i_k,l_k}^{-1},\tilde{A}_{i_k,l_k\pm 2r_{i_k}}^{-1})\text{ or }|l_{k_1}-l_{k_2}|\leq r^{\vee}$$
So in our situation we have $d_1(Z_k,Z_{k'}[ps])=d_2(Z_{k'},Z_k[ps])=0$. In particular:
$$D_1(Z_k,Z_{k'})-D_2(Z_{k'},Z_k)$$
$$=d_1(Z_k,Z_{k'})-d_2(Z_{k'},Z_k)+\underset{p\neq 0}{\sum}(d_1(Z_k,Z_{k'}[ps])-d_2(Z_{k'},Z_k[ps]))=d_1(Z_k,Z_{k'})-d_2(Z_{k'},Z_k)$$
\qed

\begin{lem}\label{detail} Let $m$ be a $\hat{\Yim}_t^s$-monomial and $l,l'\in\ZZ$. 
$$l'\geq l+s-r_i\Rightarrow u_{i,l'}(\pi_l(m))=0$$
$$l'\leq l+r_i-s+1\Rightarrow u_{i,l'}(\pi_l(m)\pi_{l-1}(m)...)=u_{i,l'}(m)$$
\end{lem}

\demo First notice that for $l,l'\in\ZZ$, we have:
$$u_{i,l'}(\tilde{Y}_{i,l})\neq 0\Rightarrow l'=l\text{ , }u_{i,l'}(\tilde{A}_{i,l})\neq 0\Rightarrow l'=l\pm r_i$$
$$i\neq j\text{ , }u_{i,l'}(\tilde{A}_{j,l})\neq 0\Rightarrow |l'-l|\leq -C_{j,i}-1\leq r^{\vee}-1$$
As $r_i\leq r^{\vee}$ we have: $u_{i,l'}(\pi_l(m))\neq 0\Rightarrow l-r^{\vee}\leq l'\leq l+r^{\vee}$.

\noindent If we suppose $l'\geq l+s-r_i\geq l+2r^{\vee}+1-r_i\geq l+r^{\vee}+1$ we have $u_{i,l'}(\pi_l(m))=0$ and this gives the first point.

\noindent We suppose that $l'\leq l+r_i-s+1$. If $k\geq l+1\geq l'+s-r_i\geq l'+r^{\vee}+1$ we have $u_{i,l'}(\pi_k(m))=0$. So: $u_{i,l'}(\pi_l(m)\pi_{l-1}(m)...)=u_{i,l'}(m)-\underset{k>l}{\sum}u_{i,l'}(\pi_k(m))=u_{i,l'}(m)$.\qed

\subsubsection{Elements of $\hat{\mathfrak{K}}_{i,t}^s$}

\begin{prop}\label{kernels} We have $\tau_{s,t}(\hat{\mathfrak{K}}_{i,t})\subset \hat{\mathfrak{K}}_{i,t}^s$. Moreover for $m$ a $i$-dominant monomial:
$$\tau_{s,t}(\overset{\leftarrow}{E}_{i,t}(m))=\tau_{s,t}(m)\tau_{s,t}(\hat{m}^i)^{-1}\overset{\leftarrow}{\underset{l\in\ZZ}{\prod}}(\tilde{Y}_{i,[l]}(1+t\tilde{A}_{i,[l+r_i]}^{-1}))^{u_{i,l}(m)}$$
where $\hat{m}^i=\underset{l\in\ZZ}{\prod}\tilde{Y}_{i,l}^{u_{i,l}(m)}\in\overline{B}_i$.\end{prop}

\demo We have to show that for $m$ a $i$-dominant monomial, $\tau_{s,t}(\overset{\leftarrow}{E}_{i,t}(m))\in \hat{\mathfrak{K}}_{i,t}^s$. The proof has three steps:

1) First we suppose that $m=\tilde{Y}_{i,l}$ where $l\in\ZZ$. We have $\overset{\leftarrow}{E}_{i,t}(\tilde{Y}_{i,l})=\tilde{Y}_{i,l}(1+t\tilde{A}_{i,l+r_i}^{-1})$, and:
$$\tau_{s,t}(\overset{\leftarrow}{E}_{i,t}(\tilde{Y}_{i,l}))=\tau_{s,t}((1+t^{-1}\tilde{A}_{i,l+r_i}^{-1})\tilde{Y}_{i,l})=1+t^{-1}\tilde{A}_{i,[l+r_i]}^{-1})\tilde{Y}_{i,[l]}=\tilde{Y}_{i,[l]}(1+t\tilde{A}_{i,[l+r_i]}^{-1})$$
and so:
$$S_{i,t}^s(\tau_{s,t}(\overset{\leftarrow}{E}_{i,t}(\tilde{Y}_{i,l}))=\tilde{Y}_{i,[l]}S_{i,[l]}-t^{-2}t\tilde{Y}_{i,[l]}\tilde{A}_{i,[l+r_i]}^{-1}S_{i,[l+2r_i]}=\tilde{Y}_{i,[l]}(S_{i,[l]}-t^{-1}\tilde{A}_{i,[l+r_i]}^{-1}S_{i,[l+2r_i]})=0$$

2) Next we suppose that $m=\underset{l\in\ZZ}{\prod}\tilde{Y}_{i,l}^{u_{i,l}}$. We have $\overset{\leftarrow}{E}_{i,t}(m)=\overset{\leftarrow}{\underset{l\in\ZZ}{\prod}}(\overset{\leftarrow}{E}_{i,t}(\tilde{Y}_{i,l}))^{u_{i,l}}$. But $r_i\leq r^{\vee}$, and in $(\overset{\leftarrow}{E}_{i,t}(\tilde{Y}_{i,l}))^{u_{i,l}}$ there are only $\tilde{Y}_{i,l}$ and $\tilde{A}_{i,l+r_i}^{-1}$. So we are in the situation of the lemma \ref{mappis}, and:
$$\tau_{s,t}(\overset{\leftarrow}{E}_{i,t}(m))=\overset{\leftarrow}{\underset{l\in\ZZ}{\prod}}(\tau_{s,t}(\overset{\leftarrow}{E}_{i,t}(\tilde{Y}_{i,l})))^{u_{i,l}}$$
As $\hat{\mathfrak{K}}_{i,t}^s$ is a subalgebra of $\hat{\Yim}_t^s$, it follows from the first step that $\tau_{s,t}(\overset{\leftarrow}{E}_{i,t}(m))\in \hat{\mathfrak{K}}_{i,t}^s$.

3) Finally let $m\in\overline{B}_i$ be an $i$-dominant monomial. As for all $l\in\ZZ$, $u_{i,l}(m)=u_{i,l}(\hat{m}^i)$, we have:
$$(\tau_{s,t}(m))^{-1}S_{i,t}^s(\tau_{s,t}(m))=(\tau_{s,t}(\hat{m}^i))^{-1}S_{i,t}^s(\tau_{s,t}(\hat{m}^i))$$
It follows from the second point that $\tau_{s,t}(\overset{\leftarrow}{E}_{i,t}(\hat{m}^i))\in\hat{\mathfrak{K}}_{i,t}^s$. Let $\hat{\chi}$ be in $\hat{\Yim}_t^s$ defined by: 
$$\hat{\chi}=\tau_{s,t}(\hat{m}^i)^{-1}\tau_{s,t}(\overset{\leftarrow}{E}_{i,t}(\hat{m}^i))$$
We have $\tau_{s,t}(m)\hat{\chi}\in\hat{\mathfrak{K}}_{i,t}^s$, because:

\noindent $S_{i,t}^s(\tau_{s,t}(m)\hat{\chi})
=S_{i,t}^s(\tau_{s,t}(m))\hat{\chi}+\tau_{s,t}(m)S_{i,t}^s(\hat{\chi})
\\=\tau_{s,t}(m)(\tau_{s,t}(\hat{m}^i))^{-1}(S_{i,t}^s(\tau_{s,t}(\hat{m}^i))\hat{\chi}+\tau_{s,t}(\hat{m}^i)S_{i,t}^s(\hat{\chi}))$$
\\=\tau_{s,t}(m)(\tau_{s,t}(\hat{m}^i))^{-1}S_{i,t}^s(\tau_{s,t}(\hat{m}^i)\hat{\chi})
\\=S_{i,t}^s(\tau_{s,t}(\overset{\leftarrow}{E}_{i,t}(m)))=0$

\noindent So it suffices to show that $\tau_{s,t}(\overset{\leftarrow}{E}_{i,t}(m))=\tau_{s,t}(m)\hat{\chi}$.

\noindent Let $\chi$ be in $\hat{\Yim}_t$ defined by: 
$$\chi=(\hat{m}^i)^{-1}\overset{\leftarrow}{E}_{i,t}(\hat{m}^i)$$
By definition of $\overset{\leftarrow}{E}_{i,t}(m)$, we have in $\hat{\Yim}_t$:
$$\overset{\leftarrow}{E}_{i,t}(m)=m\chi$$
In particular we want to show that $\tau_{s,t}(m\chi)=\tau_{s,t}(m)\tau_{s,t}(\hat{m}^i)^{-1}\tau_{s,t}(\hat{m}^i\chi)$. Let $\lambda_{m'}(t)$ be in $\ZZ[t^{\pm}]$ such that:
$$\chi=\underset{m'\in\overline{A}}{\sum}\lambda_{m'}(t)m'$$
If $\lambda_{m'}(t)\neq 0$ then $m'$ is of the form $m'=\tilde{A}_{i,l_1}^{-1}...\tilde{A}_{i,l_k}^{-1}$. As $\tau_{s,t}$ is $\ZZ[t^{\pm}]$-linear, it suffices to show that for all $m'$ of this form, we have:
$$\tau_{s,t}(m)\tau_{s,t}(\hat{m}^i)^{-1}\tau_{s,t}(\hat{m}^im')=\tau_{s,t}(mm')$$
That is to say $\alpha=\beta$ where $\alpha,\beta\in\ZZ$ are defined by: 
$$\tau_{s,t}(mm')=t^{\alpha}\tau_{s,t}(m)\tau_{s,t}(m')\text{ and }\tau_{s,t}(\hat{m}^im')=t^{\beta}\tau_{s,t}(\hat{m}^i)\tau_{s,t}(m')$$
We can suppose without loss of generality that $m\in\overset{\leftarrow}{A}$ and $m'\in\overset{\leftarrow}{A}$ (because $\tau_{s,t}$ is $\ZZ[t^{\pm}]$-linear). Let us compute $\alpha$. First we have in $\hat{\Yim}_t$:
$$mm'=t^{2\underset{l'>l}{\sum}d_2(\pi_l(m),\pi_{l'}(m'))-d_1(\pi_{l'}(m'),\pi_l(m))}\overset{\leftarrow}{\underset{l\in\ZZ}{\prod}}\pi_l(m)\pi_l(m')$$
We are in the situation of lemma \ref{mappis}, so:
$$\tau_{s,t}(mm')=t^{2\underset{l'>l}{\sum}d_2(\pi_l(m),\pi_{l'}(m'))-d_1(\pi_{l'}(m'),\pi_l(m))}\overset{\leftarrow}{\underset{l\in\ZZ}{\prod}}\tau_{s,t}(\pi_l(m))\tau_{s,t}(\pi_l(m'))$$
But we have in $\hat{\Yim}_t^s$ (lemma \ref{dspis}):
$$\tau_{s,t}(m)\tau_{s,t}(m')=t^{2\underset{l'>l}{\sum}D_2(\pi_l(m),\pi_{l'}(m'))-D_1(\pi_{l'}(m'),\pi_l(m))}\overset{\leftarrow}{\underset{l\in\ZZ}{\prod}}\tau_{s,t}(\pi_l(m))\tau_{s,t}(\pi_l(m'))$$
And we get:
$$\alpha=2\underset{l'>l}{\sum}d_2(\pi_l(m),\pi_{l'}(m'))-d_1(\pi_{l'}(m'),\pi_l(m))-2\underset{l'>l}{\sum}D_2(\pi_l(m),\pi_{l'}(m'))-D_1(\pi_{l'}(m'),\pi_l(m))$$
And so we have from lemma \ref{dspis}:

\noindent $\alpha=2\underset{l'>l}{\sum}(d_2(\pi_l(m),\pi_{l'}(m'))-d_1(\pi_{l'}(m'),\pi_l(m)))-2\underset{l'>l,r\in\ZZ}{\sum}(d_2(\pi_l(m)[rs],\pi_{l'}(m'))-d_1(\pi_{l'}(m'),\pi_l(m)[rs]))
\\=-2\underset{l'>l,r\neq 0}{\sum}(d_2(\pi_l(m)[rs],\pi_{l'}(m'))-d_1(\pi_{l'}(m'),\pi_l(m)[rs]))$

\noindent But we have $\pi_{l'}(m')$ of the form $\tilde{A}_{i,l'}^{-v_{i,l'}}$, and so:

\noindent $\alpha=-2\underset{l'>l,r\neq 0}{\sum}v_{i,l'}(m')(u_{i,l'+r_i}(\pi_l(m)[rs])-u_{i,l'-r_i}(\pi_l(m)[rs]))
\\=-2\underset{l'\in\ZZ}{\sum}v_{i,l'}(m')\underset{l<l',r\neq 0}{\sum}(u_{i,l'+r_i-rs}-u_{i,l'-r_i-rs})(\pi_l(m))
\\=-2\underset{l'\in\ZZ}{\sum}v_{i,l'}(m')\underset{r\neq 0}{\sum}(u_{i,l'+r_i-rs}-u_{i,l'-r_i-rs})(\pi_{l'-1}(m)\pi_{l'-2}(m)...)$

\noindent We use lemma \ref{detail}:

\noindent $\alpha=-2\underset{l'\in\ZZ}{\sum}v_{i,l'}(m')\underset{r>0}{\sum}(u_{i,l'+r_i-rs}-u_{i,l'-r_i-rs})(\pi_{l'-1}(m)\pi_{l'-2}(m)...)
\\=-2\underset{l'\in\ZZ}{\sum}v_{i,l'}(m')\underset{r>0}{\sum}(u_{i,l'+r_i-rs}-u_{i,l'-r_i-rs})(m)$

\noindent It depends only of the $u_{i,l}(m)$, so with the same computation we get:
$$\beta=-2\underset{l'\in\ZZ}{\sum}v_{i,l'}(m')\underset{r>0}{\sum}(u_{i,l'+r_i-rs}-u_{i,l'-r_i-rs})(\hat{m}^i)$$
and we can conclude $\alpha=\beta$ because for all $l\in\ZZ$, $u_{i,l}(m)=u_{i,l}(\hat{m}^i)$.\qed

Note that there is another more direct proof if $C$ is symmetric (in particular if $C$ is of type $ADE$):

\demo Let $m$ be an $i$-dominant monomial. 
$$\overset{\leftarrow}{m}=\overset{\leftarrow}{\underset{l\in\ZZ}{\prod}}\underset{j\in I}{\prod}\tilde{A}_{j,l+1}^{-v_{j,l+1}}\underset{j\in I}{\prod}\tilde{Y}_{j,l}^{y_{j,l}}=\overset{\leftarrow}{\underset{l\in\ZZ}{\prod}}\tilde{A}_{i,l+1}^{-v_{i,l+1}}\tilde{Y}_{i,l}^{y_{i,l}}(\underset{j\neq i}{\prod}\tilde{Y}_{j,l}^{y_{j,l}}\underset{j\neq i}{\prod}\tilde{A}_{j,l}^{-v_{j,l}})$$
For $l\in\ZZ$, let $M_l=\underset{j\neq i}{\prod}\tilde{Y}_{j,l}^{y_{j,l}}\underset{j/ C_{i,j}=0}{\prod}\tilde{A}_{j,l}^{-v_{j,l}}$. We have $\overset{\leftarrow}{E}_{i,t}(M_l)=M_l$. The $\tilde{Y}_{i,l}$ and the $\tilde{A}_{j,l}^{-1}$ with $C_{i,j}=-1$ have the same relations with the $\tilde{A}_{i,l}^{-1}$, so we use indifferently the notation $Z_{i,l}$ for $\tilde{Y}_{i,l}$ or $\tilde{A}_{j,l}^{-1}$. The power of $Z_{i,l}$ is: 
$$z_{i,l}=y_{i,l}+\underset{j/C_{j,i}=-1}{\sum} v_{j,l+1}+v_{j,l-1}=u_{i,l}+v_{i,l-1}+v_{i,l+1}$$
In particular we have:
$$\overset{\leftarrow}{E}_{i,t}(m)=\overset{\leftarrow}{\underset{l\in\ZZ}{\prod}}(Z_{i,l}\tilde{A}_{i,l+1}^{-1})^{v_{i,l+1}}\overset{\leftarrow}{E}_{i,t}(Z_{i,l}^{u_{i,l}})M_l Z_{i,l}^{v_{i,l-1}}$$
and it follows from lemma \ref{mappis} that:
$$\tau_{s,t}(\overset{\leftarrow}{E}_{i,t}(m))=\overset{\leftarrow}{\underset{l\in\ZZ}{\prod}}(\tau_{s,t}(Z_{i,l}\tilde{A}_{i,l+1}^{-1}))^{v_{i,l+1}}\tau_{s,t}(\overset{\leftarrow}{E}_{i,t}(Z_{i,l}))^{u_{i,l}}\tau_{s,t}(M_l) \tau_{s,t}(Z_{i,l}^{v_{i,l-1}})\in\hat{\mathfrak{K}}_{i,t}^s$$\qed

\subsection{Description of $\chi_{\epsilon,t}$} In this section we prove the following theorem (the map $p_s$ is defined in section \ref{invol}):

\begin{thm}\label{axquat} If $\chi_{q,t}(\underset{i\in I,0\leq l\leq s-1}{\prod}X_{i,l}^{x_{i,[l]}})=\underset{m\in\overline{A}}{\sum}\lambda_m(t)m$, then:
$$\chi_{\epsilon,t}(\underset{i\in I,0\leq l\leq s-1}{\prod}X_{i,[l]}^{x_{i,l}})=\underset{m\in\overline{A}}{\sum}\lambda_m(t)t^{D_1^-(m)+D_2^-(m)}p_s(m)$$
where for $m$ a $\hat{\Yim}_t$-monomial:
$$D_1^-(m)=\underset{k<0}{\sum}d_1(m,m[ks])\text{ , }D_2^-(m)=\underset{k<0}{\sum}d_2(m,m[ks])$$
\end{thm}

\noindent Note that this result is a generalization of the axiom 4 of \cite{Nab} to the non necessarily finite simply laced case. In particular our construction fits with \cite{Nab} in the $ADE$-case.

\subsubsection{Description of the basis $\overline{A}$}

\begin{lem}\label{descbasis} For $m$ a $\hat{\Yim}_t$-monomial we have $t^{\gamma}\overset{\rightarrow}{m}\in\overline{A}$ and $t^{-\gamma-2d_1(m,m)}\overset{\leftarrow}{m}\in\overline{A}$ where:
$$\gamma=\underset{l\in\ZZ}{\sum}(\underset{i\in I}{\sum}v_{i,l}^2(m)-\underset{i,j/C_{i,j}+r_i=-1}{\sum}v_{i,l}(m)v_{j,l}(m)-\underset{i,j/C_{i,j}=-3\text{ and }r_i=1}{\sum}v_{i,l}(m)(v_{j,l+1}(m)+v_{j,l-1}(m)))$$
\end{lem}

\demo We have $\overline{\overset{\rightarrow}{m}}=\overset{\leftarrow}{m}=t^{2\beta}\overset{\rightarrow}{m}$ where:
$$\beta=\underset{l>l'}{\sum}d_1(\pi_l(m),\pi_{l'}(m))-d_2(\pi_{l'}(m),\pi_l(m))$$
$$=d_1(m,m)-\underset{l\in\ZZ}{\sum}d_1(\pi_l(m),\pi_l(m))-\underset{l<l'}{\sum}d_1(\pi_l(m),\pi_{l'}(m))+d_2(\pi_l(m),\pi_{l'}(m))$$
So $\overline{t^{\gamma}\overset{\rightarrow}{m}}=t^{2d_1(m,m)}t^{\gamma}\overset{\rightarrow}{m}$ where :
$$\gamma=-\underset{l\in\ZZ}{\sum}d_1(\pi_l(m),\pi_l(m))-\underset{l<l'}{\sum}d_1(\pi_l(m),\pi_{l'}(m))+d_2(\pi_l(m),\pi_{l'}(m))$$
But for $l\in\ZZ$ we have 
$$d_1(\pi_l(m),\pi_l(m))=-\underset{i\in I}{\sum}v_{i,l}^2(m)+\underset{i,j/C_{i,j}=-2\text{ and }r_i=-1}{\sum}v_{i,l}(m)v_{j,l}(m)+\underset{i,j/C_{i,j}=-3\text{ and }r_i=2}{\sum}v_{i,l}(m)v_{j,l}(m)$$
$$=-\underset{i\in I}{\sum}v_{i,l}^2(m)+\underset{i,j/C_{i,j}+r_i=-1}{\sum}v_{i,l}(m)v_{j,l}(m)$$ 
For $l<l'$ we have: 
$$d_1(\pi_l(m),\pi_{l'}(m))=\delta_{l'=l+1}\underset{i,j/C_{i,j}=-3\text{ and }r_i=1}{\sum}v_{i,l}(m)v_{j,l+1}(m)$$
$$d_2(\pi_l(m),\pi_{l'}(m))=\delta_{l'=l+1}\underset{i,j/C_{i,j}=-3\text{ and }r_i=1}{\sum}v_{i,l+1}(m)v_{j,l}(m)$$
and we get for $\gamma$ the annonced value.

\noindent For the second point we show that $t^{-\gamma-2d_1(m,m)}\overset{\leftarrow}{m}\in\overline{A}$:
$$\overline{t^{-\gamma-2d_1(m,m)}\overset{\leftarrow}{m}}=t^{\gamma+2d_1(m,m)}\overline{\overset{\leftarrow}{m}}=t^{\gamma+2d_1(m,m)-2\beta}\overset{\leftarrow}{m}=t^{-\gamma}\overset{\leftarrow}{m}=t^{2d_1(m,m)}(t^{-\gamma-2d_1(m,m)}\overset{\leftarrow}{m})$$
\qed

\subsubsection{Description of $\tau_{s,t}$}

\begin{prop}\label{dmoins} For $m\in\overline{A}$ we have:
$$\tau_{s,t}(m)=t^{D_1^-(m)+D_2^-(m)}p_s(m)$$
\end{prop}

\demo Using lemma \ref{descbasis} we can write $m=t^{-\gamma-2d_1(m,m)}\overset{\leftarrow}{m}$. So we have:
$$\tau_{s,t}(m)=t^{-\gamma-2d_1(m,m)}\overset{\leftarrow}{\underset{l\in\ZZ}{\prod}}\tau_{s,t}(\pi_l)$$
where $\pi_l=\pi_l(m)$. So we have $\overline{\tau_{s,t}(m)}=t^{2\alpha}\tau_{s,t}(m)$ where:
$$\alpha=\gamma+2d_1(m,m)+\underset{l<l'}{\sum}d_1(\tau_{s,t}(\pi_l),\tau_{s,t}(\pi_{l'}))-d_2(\tau_{s,t}(\pi_{l'}),\tau_{s,t}(\pi_l))$$
$$=\gamma+2d_1(m,m)+\underset{l<l'}{\sum}D_1(\pi_l,\pi_{l'})-D_2(\pi_{l'},\pi_l)$$
So it suffices to show that $\alpha=-D_1^-(m)-D_2^-(m)+d_1(p_s(m),p_s(m))$. But we have:
$$d_1(p_s(m),p_s(m))=D_1(m,m)=\underset{l<l'}{\sum}D_1(\pi_l,\pi_{l'})+\underset{l\geq l'}{\sum}D_1(\pi_l,\pi_{l'})$$
So we want to show:
$$-D_2(\pi_{l'},\pi_l)=\underset{l\geq l'}{\sum}D_1(\pi_l,\pi_{l'})-(d_1(m,m)+D_1^-(m))-(d_2(m,m)+D_2^-(m))-\gamma$$
The second term is:
$$\underset{l\geq l',r\in\ZZ}{\sum}d_1(\pi_l,\pi_{l'}[rs])-\underset{l,l'\in\ZZ, r\leq 0}{\sum}(d_1(\pi_l,\pi_{l'}[rs])+d_2(\pi_l,\pi_{l'}[rs]))+\underset{l\in\ZZ}{\sum}d_2(\pi_l,\pi_l)+\underset{l<l'}{\sum}(d_1(\pi_l,\pi_{l'})+d_2(\pi_l,\pi_{l'}))$$
But for $l<l'$ and $r<0$ (resp. $l\geq l'$ and $r>0$) we have $d_1(\pi_l,\pi_{l'}[rs])=d_2(\pi_l,\pi_{l'}[rs])=0$. So this term is:
$$\underset{l\geq l',r\leq 0}{\sum}d_1(\pi_l,\pi_{l'}[rs])-\underset{l\geq l', r\leq 0}{\sum}(d_1(\pi_l,\pi_{l'}[rs])+d_2(\pi_l,\pi_{l'}[rs]))+\underset{l\in\ZZ}{\sum}d_2(\pi_l,\pi_l)$$
$$=-\underset{l>l',r\leq 0}{\sum}d_2(\pi_l,\pi_{l'}[rs])=-\underset{l>l',r\in\ZZ}{\sum}d_2(\pi_l,\pi_{l'}[rs])=-\underset{l>l'}{\sum}D_2(\pi_l,\pi_{l'})$$\qed

\section{Applications}\label{app} In this section we see how we can generalize at roots of unity results of \cite{her} about Kazhdan-Lusztig polynomials and quantization of the Grothendieck ring. We suppose that $i\neq j\Rightarrow C_{i,j}C_{j,i}\leq 3$.

\noindent Such constructions were made by Nakajima \cite{Nab} in the simply laced case.

\subsection{Reminder: Kazhdan-Lusztig polynomials in the generic case \cite{Nab}\cite{her}} In this section we suppose that $s=0$. The involution of $\hat{\Yim}_t$ is naturally extended to an involution of $\hat{\Yim}_t^{\infty}$.

\noindent For $m$ a dominant $\hat{\Yim}_t$-monomial we set\label{tetm}:
$$\overset{\rightarrow}{E}_t(m)=m(\overset{\rightarrow}{\underset{l\in\ZZ}{\prod}}\underset{i\in I}{\prod}\tilde{Y}_{i,l}^{u_{i,l}(m)})^{-1}\overset{\rightarrow}{\underset{l\in\ZZ}{\prod}}\underset{i\in I}{\prod}\hat{F}_t(\tilde{Y}_{i,l})^{u_{i,l}(m)}$$
We denote by $\hat{\mathfrak{K}}_t^{f,\infty}\subset\hat{\mathfrak{K}}_t^{\infty}$ the subset of elements with only a finite number of dominant monomials. 

\noindent We show as in \cite{her} that for $m\in\overline{B}$, $C(m)\cap \overline{B}$ is finite, $\overset{\rightarrow}{E}_t(m)\in \hat{\mathfrak{K}}_t^{f,\infty}$, and:

\begin{prop}(\cite{her}) $\hat{\mathfrak{K}}_t^{f,\infty}$ is a subalgebra of $\hat{\mathfrak{K}}_t^{\infty}$, and:
$$\hat{\mathfrak{K}}_t^{f,\infty}=\underset{m\in\overline{B}}{\bigoplus}\ZZ[t^{\pm}]\hat{F}_t(m)=\underset{m\in\overline{B}}{\bigoplus}\ZZ[t^{\pm}]\overset{\rightarrow}{E}_t(m)$$
Moreover $\hat{\mathfrak{K}}_t^{f,\infty}$ is stable by the involution.
\end{prop}

\noindent For $m$ a $\hat{\Yim}_t^s$-monomial there is a unique $\alpha(m)\in\ZZ$\label{alpham} such that $\overline{t^{\alpha(m)}m}=t^{\alpha(m)}m$ (see the proof of lemma 6.12 of \cite{her}).

\noindent Let $\hat{A}^{\text{inv}}=\{t^{\alpha(m)}m/m\in\overline{A}\}$\label{ainv} and $\hat{B}^{\text{inv}}=\{t^{\alpha(m)}m/m\in\overline{B}\}$. 

\noindent The following theorem was given in \cite{Nab} for the $ADE$-case and in \cite{her} for the general finite case:

\begin{thm}\label{klgen} For $m\in \hat{B}^{\text{inv}}$ there is a unique $\hat{L}_t(m)\in\hat{\mathfrak{K}}_t^{f, \infty}$ such that\label{tltm}:
$$\overline{\hat{L}_t(m)}=\hat{L}_t(m)$$
$$\overset{\rightarrow}{E}_t(m)=\hat{L}_t(m)+\underset{m'<m, m' \in \hat{B}^{\text{inv}}}{\sum}P_{m',m}(t)\hat{L}_t(m')$$
where $P_{m',m}(t)\in t^{-1}\ZZ[t^{-1}]$\label{pmm}.
\end{thm}

\subsection{Kazhdan-Lusztig polynomials at roots of unity} In this section we suppose that $s> 2r^{\vee}$. The involution of $\hat{\Yim}_t^s$ is extended to an involution of $\hat{\Yim}_t^{s, \infty}$.

\subsubsection{Construction of stable subalgebras}\label{sblsubalg} For $m\in \overline{B}_i^s$ a $i$-dominant $\hat{\Yim}_t^s$-monomial we set:
$$\overset{\leftarrow}{E}_{i,t}(m)=m(\underset{i\in I,l=0.. s-1}{\prod}\tilde{Y}_{i,[l]}^{u_{i,[l]}(m)})^{-1}\overset{\leftarrow}{\underset{i\in I,l=0.. s-1}{\prod}}(\tilde{Y}_{i,[l]}(1+t\tilde{A}_{i,[l+r_i]}^{-1}))^{u_{i,[l]}(m)}$$
In particular the formula of proposition \ref{kernels} implies:
$$\overset{\leftarrow}{E}_{i,t}(m)=m(\tau_{s,t}(M))^{-1}\tau_{s,t}(\overset{\leftarrow}{E}_{i,t}(M))\text{ where }M=\underset{l=0... s-1}{\prod}\tilde{Y}_{i,l}^{u_{i,[l]}(m)}$$
We define:
$$\tilde{\mathfrak{K}}_{i,t}^s=\underset{m\in\overline{B}_i^s}{\bigoplus}\ZZ[t^{\pm}]\overset{\leftarrow}{E}_{i,t}(m)$$\label{tkits}
In particular if $\chi\in \tilde{\mathfrak{K}}_{i,t}^s$ has no $i$-dominant monomial then $\chi=0$.

\begin{lem}\label{aidetrois} We have $\tau_{s,t}(\hat{\mathfrak{K}}_{i,t})\subset \tilde{\mathfrak{K}}_{i,t}^s\subset \hat{\mathfrak{K}}_{i,t}^s$. Moreover $\tilde{\mathfrak{K}}_{i,t}^s$ is a subalgebra of $\hat{\mathfrak{K}}_{i,t}^s$ and is stable by the involution.\end{lem}

\demo As $\hat{\mathfrak{K}}_{i,t}^s$ is a subalgebra of $\hat{\Yim}_t^s$ and $\tilde{Y}_{i,[l]}(1+t\tilde{A}_{i,[l+r_i]}^{-1})\in \hat{\mathfrak{K}}_{i,t}^s$, $m(\underset{i\in I,l=0.. s-1}{\prod}\tilde{Y}_{i,[l]}^{u_{i,[l]}(m)})^{-1}\in \hat{\mathfrak{K}}_{i,t}^s$ we have $\tilde{\mathfrak{K}}_{i,t}^s\subset \hat{\mathfrak{K}}_{i,t}^s$.

\noindent Let us show that $\underset{m\in\overline{B}_i^s}{\bigoplus}\ZZ[t^{\pm}]\overset{\leftarrow}{E}_{i,t}(m)$ is a subalgebra of $\tilde{\mathfrak{K}}_{i,t}^s$ (note that in the generic case $s=0$ this point needs no proof because $\tilde{\mathfrak{K}}_{i,t}=\hat{\mathfrak{K}}_{i,t}$). For this point our proof is analogous to theorem 3.8 of \cite{Nab}. It suffices to show that for $0\leq k\leq s-1$, $M=\underset{l\in\ZZ/s\ZZ}{\prod}\tilde{Y}_{i,l}^{u_{i,l}}$ we have $\overset{\leftarrow}{E}_{i,t}(M)\overset{\leftarrow}{E}_{i,t}(\tilde{Y}_{i,k})\in \underset{m\in\overline{B}_i^s}{\bigoplus}\ZZ[t^{\pm}]\overset{\leftarrow}{E}_{i,t}(m)$. We can suppose without loss of generality that we are in the $sl_2$-case and that $r_i=r_1=1$. The $\overset{\leftarrow}{E}_t(\tilde{Y}_k)$ do not commute with $\overset{\leftarrow}{E}_t(\tilde{Y}_{k-2}^{u_{k-2}})$ and $\overset{\leftarrow}{E}_t(\tilde{Y}_{k+2}^{u_{k+2}})$. So if $k\geq 2$ that fact that $s\neq 0$ do not change anything and the result follows from the generic case. If $k=0$, we have:
$$\overset{\leftarrow}{E}_t(m)\overset{\leftarrow}{E}_t(\tilde{Y}_{i,0})=\overset{\leftarrow}{E}_t(m\tilde{Y}_{i,0})
+\overset{\leftarrow}{E}_t(\tilde{Y}_{i,0}^{u_{i,0}}\tilde{Y}_{i,1}^{u_{i,1}})[\overset{\leftarrow}{E}_t(\tilde{Y}_{i,2}^{u_{i,2}}),\overset{\leftarrow}{E}_t(\tilde{Y}_{i,0})]\overset{\leftarrow}{E}_t(\tilde{Y}_{i,3}^{u_{i,3}}...\tilde{Y}_{i,s-1}^{u_{i,s-1}})$$
$$+\overset{\leftarrow}{E}_t(\tilde{Y}_{i,0}^{u_{i,0}}...\tilde{Y}_{i,s-3}^{u_{i,s-3}})[\overset{\leftarrow}{E}_t(\tilde{Y}_{i,s-2}^{u_{i,s-2}}),\overset{\leftarrow}{E}_t(\tilde{Y}_{i,0})]\overset{\leftarrow}{E}_t(\tilde{Y}_{i,s-1}^{u_{i,s-1}})$$

\noindent It follows from the study of the generic case that:
$$[\overset{\leftarrow}{E}_t(\tilde{Y}_{i,2}^{u_{i,2}}),\overset{\leftarrow}{E}_t(\tilde{Y}_{i,0})]\in\underset{0\leq r<u_{i,2}}{\bigoplus}\ZZ[t^{\pm}]\overset{\leftarrow}{E}_t(\tilde{Y}_{i,2})^r$$
$$[\overset{\leftarrow}{E}_t(\tilde{Y}_{i,s-2}^{u_{i,s-2}}),\overset{\leftarrow}{E}_t(\tilde{Y}_{i,0})]\in\underset{0\leq r<u_{i,s-2}}{\bigoplus}\ZZ[t^{\pm}]\overset{\leftarrow}{E}_t(\tilde{Y}_{i,s-2})^r$$
and we can conclude by induction. The case $k=1$ is studied in the same way.

\noindent Let us study the stability by the involution: we see that $\overline{\overset{\leftarrow}{E}_{i,t}(\tilde{Y}_{i,l})}=\overset{\leftarrow}{E}_{i,t}(\tilde{Y}_{i,l})$, and:
$$\overline{\overset{\leftarrow}{E}_{i,t}(m)}=\overset{\rightarrow}{\underset{i\in I,l=s-1, s-2,... 0}{\prod}}\overline{\overset{\leftarrow}{E}_{i,t}(\tilde{Y}_{i,[l]})}^{u_{i,[l]}(m)}\overset{\leftarrow}{E}_{i,t}(m(\overset{\rightarrow}{\underset{i\in I,l=0.. s-1}{\prod}}\tilde{Y}_{i,[l]}^{u_{i,[l]}(m)})^{-1})\in \tilde{\mathfrak{K}}_{i,t}$$

\noindent Let us show that $\tau_{s,t}(\hat{\mathfrak{K}}_{i,t})\subset \tilde{\mathfrak{K}}_{i,t}^s$: the formula of proposition \ref{kernels} implies that for $m\in\overline{B}_i$:
$$\tau_{s,t}(\overset{\leftarrow}{E}_{i,t}(m))=\tau_{s,t}(m)\tau_{s,t}(\hat{m}^i)^{-1}\overset{\leftarrow}{\underset{l\in\ZZ}{\prod}}(\tilde{Y}_{i,[l]}(1+t\tilde{A}_{i,[l+r_i]}^{-1}))^{u_{i,l}(m)}$$
$$=\overset{\leftarrow}{E}_{i,t}(\tau_{s,t}(m)\tau_{s,t}(\hat{m}^i)^{-1})\overset{\leftarrow}{\underset{l\in\ZZ}{\prod}}\overset{\leftarrow}{E}_{i,t}(\tilde{Y}_{i,[l]}^{u_{i,l}(m)})$$
and we can conclude because $\tilde{\mathfrak{K}}_{i,t}^s$ is an algebra.\qed

\noindent We define the completion $\tilde{\mathfrak{K}}_{i,t}^{s,\infty}\subset\hat{\mathfrak{K}}_{i,t}^{s,\infty}$\label{tkitinfs} (as in section \ref{cpltdalg}) and:
$$\tilde{\mathfrak{K}}_t^{s,\infty}=\underset{i\in I}{\bigcap}\tilde{\mathfrak{K}}_{i,t}^{s,\infty}$$

\noindent For $m\in\overline{B}^s$ we define $\overset{\rightarrow}{E}_t(m)=m(\tau_{s,t}(M))^{-1}\tau_{s,t}(\overset{\rightarrow}{E}_t(M))$ where $M=\underset{i\in I,l=0...s}{\prod}\tilde{Y}_{i,l}^{u_{i,[l]}(m)}$.

\subsubsection{Polynomials at roots of unity (finite case)} In this section we suppose that $C$ is finite. Note that it follows from the lemma \ref{dommons} that for $m\in\overline{B}^s$, the set $C(m)\cap \overline{B}^s$ is finite.

\noindent We denote by $\tilde{\mathfrak{K}}_t^{s,f, \infty}$ the set of elements of $\tilde{\mathfrak{K}}_t^{s,\infty}$ with only a finite number of dominant monomials.

\begin{lem}\label{algroot} $\tilde{\mathfrak{K}}_t^{s, f,\infty}$ is a subalgebra of $\hat{\Yim}_t^{s,\infty}$, is stable by the involution, and:
$$\tilde{\mathfrak{K}}_t^{s, f,\infty}=\underset{m\in \overline{B}^s}{\bigoplus}\ZZ[t^{\pm}]\overset{\rightarrow}{E}_t(m)$$\end{lem}

\demo It follows from lemma \ref{aidetrois} that $\tilde{\mathfrak{K}}_t^{s,\infty}$ is a subalgebra of $\hat{\Yim}_t^{s,\infty}$. Let $m$ be in $\overline{B}^s$. For all $i\in I$ we have $m(\tau_{s,t}(M))^{-1}=\overset{\leftarrow}{E}_{i,t}(m(\tau_{s,t}(M))^{-1})$ and so $m(\tau_{s,t}(M))^{-1}\in \tilde{\mathfrak{K}}_t^{s,\infty}$. But $\tau_{s,t}(\overset{\rightarrow}{E}_t(M))\in\tilde{\mathfrak{K}}_{i,t}^{s,\infty}$ for all $i\in I$. So $\overset{\rightarrow}{E}_t(m)\in \tilde{\mathfrak{K}}_t^{s,\infty}$. Moreover lemma \ref{dommons} shows that $\overset{\rightarrow}{E}_t(m)$ has only a finite number of dominant monomials, so $\overset{\rightarrow}{E}_t(m)\in \tilde{\mathfrak{K}}_t^{s,f,\infty}$. It follows from lemma \ref{aidetrois} that a maximal monomial of an element of $\tilde{\mathfrak{K}}_t^{s,f,\infty}$ is dominant, and so we have the other inclusion $\tilde{\mathfrak{K}}_t^{s, f,\infty}\subset\underset{m\in \overline{B}^s}{\bigoplus}\ZZ[t^{\pm}]\overset{\rightarrow}{E}_t(m)$.

\noindent It follows from lemma \ref{aidetrois} that $\tilde{\mathfrak{K}}_t^{s,\infty}$ is stable by the involution. But for $m$ a dominant monomial, $\overline{m}$ is a dominant monomial and so $\tilde{\mathfrak{K}}_t^{s,f,\infty}$ is stable by the involution.

\noindent As $\tilde{\mathfrak{K}}_t^{s,\infty}$ is an algebra, $\tilde{\mathfrak{K}}_t^{s,f,\infty}$ is an algebra if for $m,m'\in\overline{B}^s$, $\overset{\rightarrow}{E}_t(m)\overset{\rightarrow}{E}_t(m')$ has only a finite number of dominant monomials. But the monomials of $\overset{\rightarrow}{E}_t(m)\overset{\rightarrow}{E}_t(m')$ are in $C(mm')$ and we can conclude with the help of lemma \ref{dommons}.\qed

\noindent Let $\hat{A}^{s, \text{inv}}=\{t^{\alpha(m)}m/m\in\overline{A}^s\}$\label{ainvs} and $\hat{B}^{s, \text{inv}}=\{t^{\alpha(m)}m/m\in\overline{B}^s\}$ where $\alpha(m)$ is defined by $\overline{t^{\alpha}m}=t^{\alpha(m)}m$ (see the proof of lemma 6.12 of \cite{her}).

\begin{thm} For $m\in \hat{B}^{s,\text{inv}}$ there is a unique $\hat{L}_t^s(m)\in\hat{\mathfrak{K}}_t^{s, f, \infty}$ such that:\label{tltms}
$$\overline{\hat{L}_t^s(m)}=\hat{L}_t^s(m)$$
$$\overset{\rightarrow}{E}_t(m)=\hat{L}_t^s(m)+\underset{m'<m, m' \in \hat{B}^{s,\text{inv}}}{\sum}P_{m',m}^s(t)\hat{L}_t^s(m')$$
where $P_{m',m}^s(t)\in t^{-1}\ZZ[t^{-1}]$.\label{pmms}
\end{thm}

\noindent The proof is analogous to the proof of theorem \ref{klgen} with the help of lemma \ref{algroot}. The result was first given by Nakajima \cite{Nab} for the $ADE$-case.

\subsubsection{Example and conjecture (finite case)} In the following example we suppose that we are in the $sl_2$-case and we study the decomposition with $m=\tilde{Y}_0\tilde{Y}_1\tilde{Y}_2$.

If $s=0$, we have:
$$\overset{\rightarrow}{E}_t(\tilde{Y}_0\tilde{Y}_1\tilde{Y}_2)=\tilde{Y}_0(1+t\tilde{A}_1^{-1})\tilde{Y}_1(1+t\tilde{A}_2^{-1})\tilde{Y}_2(1+t\tilde{A}_3^{-1})$$
$$=\hat{L}_t(\tilde{Y}_0\tilde{Y}_1\tilde{Y}_2)+t^{-1}\hat{L}_t(t^2\tilde{Y}_0\tilde{A}_1^{-1}\tilde{Y}_1\tilde{Y}_2)$$
where:
$$\hat{L}_t(\tilde{Y}_0\tilde{Y}_1\tilde{Y}_2)=\tilde{Y}_0\tilde{Y}_1\tilde{Y}_2(1+t\tilde{A}_3^{-1}(1+t\tilde{A}_1^{-1}))(1+t\tilde{A}_2^{-1})$$
$$\hat{L}_t(t^2\tilde{Y}_0\tilde{A}_1^{-1}\tilde{Y}_1\tilde{Y}_2)=t^2\tilde{Y}_0\tilde{A}_1^{-1}\tilde{Y}_1\tilde{Y}_2(1+t\tilde{A}_2^{-1})$$

If $s=3$, we have:

\noindent $\overset{\rightarrow}{E}_t(\tilde{Y}_0\tilde{Y}_1\tilde{Y}_2)
=\tau_{s,t}(\tilde{Y}_0(1+t\tilde{A}_1^{-1})\tilde{Y}_1(1+t\tilde{A}_2^{-1})\tilde{Y}_2(1+t\tilde{A}_3^{-1}))
\\=\tilde{Y}_0\tilde{Y}_1\tilde{Y}_2+t\tilde{Y}_0\tilde{A}_1^{-1}\tilde{Y}_1\tilde{Y}_2+t^{-1}\tilde{Y}_0\tilde{Y}_1\tilde{A}_2^{-1}\tilde{Y}_2+t^{-1}\tilde{Y}_0\tilde{Y}_1\tilde{A}_2^{-1}\tilde{Y}_2
\\+t^2\tilde{Y}_0\tilde{A}_1^{-1}\tilde{Y}_1\tilde{Y}_2\tilde{A}_2^{-1}+\tilde{Y}_0\tilde{Y}_1\tilde{A}_2^{-1}\tilde{Y}_2\tilde{A}_3^{-1}+\tilde{Y}_0\tilde{Y}_1\tilde{A}_1^{-1}\tilde{A}_2^{-1}\tilde{Y}_2+t^{-3}\tilde{Y}_0\tilde{A}_1^{-1}\tilde{Y}_1\tilde{A}_2^{-1}\tilde{Y}_2\tilde{A}_3^{-1}$

\noindent and so:
$$\overset{\rightarrow}{E}_t(\tilde{Y}_0\tilde{Y}_1\tilde{Y}_2)=\hat{L}_t^s(\tilde{Y}_0\tilde{Y}_1\tilde{Y}_2)+t^{-1}\hat{L}_t^s(t^2\tilde{Y}_0\tilde{A}_1^{-1}\tilde{Y}_1\tilde{Y}_2)+t^{-1}\hat{L}_t^s(\tilde{Y}_0\tilde{Y}_1\tilde{A}_2^{-1}\tilde{Y}_2)+t^{-1}\hat{L}_t^s(\tilde{Y}_0\tilde{Y}_1\tilde{Y}_2\tilde{A}_3^{-1})$$
where:
$$\hat{L}_t^s(\tilde{Y}_0\tilde{Y}_1\tilde{Y}_2)=\tilde{Y}_0\tilde{Y}_1\tilde{Y}_2+t^{-3}\tilde{Y}_0\tilde{A}_1^{-1}\tilde{Y}_2\tilde{A}_3^{-1}\tilde{Y}_4\tilde{A}_5^{-1}$$
$$\hat{L}_t^s(t^2\tilde{Y}_0\tilde{A}_1^{-1}\tilde{Y}_1\tilde{Y}_2)=t^2\tilde{Y}_0\tilde{A}_1^{-1}\tilde{Y}_1\tilde{Y}_2(1+\tilde{A}_3^{-1})$$
$$\hat{L}_t^s(\tilde{Y}_0\tilde{Y}_1\tilde{A}_2^{-1}\tilde{Y}_2)=\tilde{Y}_0\tilde{Y}_1\tilde{A}_2^{-1}\tilde{Y}_2(1+t\tilde{A}_2^{-1})$$
$$\hat{L}_t^s(\tilde{Y}_0\tilde{Y}_1\tilde{Y}_2\tilde{A}_3^{-1})=\tilde{Y}_0\tilde{Y}_1\tilde{Y}_2\tilde{A}_3^{-1}(1+t\tilde{A}_2^{-1})$$

\noindent In particular we see in this example that the decomposition of $\overset{\rightarrow}{E}_t(m)$ in general is not necessarily the same if $s=0$ or $s\neq 0$.

\noindent We recall that irreducible representations  of $\U_q(\hat{\Glie})$ (resp. $\U_{\epsilon}^{\text{res}}(\hat{\Glie})$) are classified by dominant monomials of $\Yim$ (resp. $\Yim^s$) or by Drinfel'd polynomials (see \cite{Cha}, \cite{Cha3}, \cite{Fre}, \cite{Fre3}). 

\noindent For $m\in B$ (resp. $m\in B^s$) we denote by $V_m^0=V_m\in\text{Rep}(\U_q(\hat{\Glie}))$ (resp. $V_m^s\in\text{Rep}(\U_{\epsilon}^{\text{res}}(\hat{\Glie}))$) the irreducible module of highest weight $m$. In particular for $i\in I,l\in\ZZ/s\ZZ$ let $V_{i,l}^s=V_{Y_{i,l}}$. The simple modules $V_{i,l}^s$ are called fundamental representations. In the ring $\text{Rep}^s$ it is denoted by $X_{i,l}$.

\noindent For $m\in B$ (resp. $m\in B^s$) we denote by $M_m^s\in\text{Rep}(\U_q(\hat{\Glie}))$ (resp. $M_m^s\in\text{Rep}(\U_{\epsilon}^{\text{res}}(\hat{\Glie}))$) the module $M_m^s=\underset{i\in I,l\in\ZZ/s\ZZ}{\bigotimes}V_{i,l}^{s, \otimes u_{i,l}(m)}$. It is called a standard module and in $\text{Rep}^s$ it is denoted by $\underset{i\in I,l\in\ZZ/s\ZZ}{\prod}X_{i,l}^{u_{i,l}(m)}$.

\noindent The irreducible $\U_q(\hat{sl_2})$-representation with highest weight $m$ is $V_m=V_{Y_0Y_2}\otimes V_{Y_1}$ (see \cite{Cha} or \cite{Fre}). In particular $\text{dim}(V_m)=6$, that is to say the number of monomials of $\hat{L}_t(m)$.

\noindent For $\epsilon$ such that $s=3$, the irreducible $\U_{\epsilon}^{\text{res}}(\hat{\Glie})$-representation with highest weight $m$ is $V^s_{m}$ the pull back by the Frobenius morphism of the $\overline{\U}(\hat{sl_2})$-module $\overline{V}$ of Drinfel'd polynomial $(1-u)$ (see \cite{Cha3} or \cite{Fre3}). In particular $\text{dim}(V^s_m)=2$, that is to say the number of monomials of $\hat{L}_t^s(m)$.

\noindent Those observations would be explained by the following conjecture which is a generalization of the  conjecture 7.3 of \cite{her} to the root of unity case. We know from \cite{Nab} that the result is true in the simply laced case (in particular in the last example).

\noindent For $m=\underset{i\in I,l\in\ZZ/s\ZZ}{\prod}Y_{i,l}^{u_{i,l}}$ a dominant $\Yim^s$-monomial let $M=\underset{i\in I,l\in\ZZ/s\ZZ}{\prod}\tilde{Y}_{i,l}^{u_{i,l}}\in\hat{\Yim}_t^s$. We suppose that $C$ is finite.

\begin{conj} For $m$ a dominant $\Yim_s$-monomial, $\hat{\Pi}(\hat{L}_t^s(M))$ is the $\epsilon$-character of the irreducible 
\\$\U_{\epsilon}^{\text{res}}(\hat{\Glie})$-representation $V_m^s$ associated to $m$. In particular for $m'$ another dominant $\Yim_s$-monomials the multiplicity of $V_{m'}^s$ in the standard module $M_m^s$ associated to $m$ is:
$$\underset{m''\in \hat{B}^{s,\text{inv}}/\hat{\Pi}(m'')=m'}{\sum}P_{m'',M}^s(1)$$
\end{conj}

\noindent Let us look at an application of the conjecture in the non-simply laced case: we suppose that $C=\begin{pmatrix}2 & -2\\-1 & 2\end{pmatrix}$ and $m=\tilde{Y}_{1,0}\tilde{Y}_{1,1}$. We have for $l\in\ZZ/s\ZZ$: 
$$\tilde{A}_{1,l}^{-1}=:\tilde{Y}_{1,l-1}\tilde{Y}_{1,l+1}\tilde{Y}_{2,l}:\text{ , }\tilde{A}_{2,l}^{-1}=:\tilde{Y}_{2,l-2}^{-1}\tilde{Y}_{2,l+2}^{-1}\tilde{Y}_{1,l-1}\tilde{Y}_{1,l+1}:$$

First we suppose that $s=0$. The formulas for $\hat{F}_t(\tilde{Y}_{1,0})$ and $\hat{F}_t(\tilde{Y}_{1,1})$ are given in \cite{her}: 
$$\hat{F}_t(\tilde{Y}_{1,0})=\tilde{Y}_{1,0}(1+t\tilde{A}_{1,1}^{-1}(1+t\tilde{A}_{2,3}^{-1}(1+t\tilde{A}_{1,5}^{-1})))$$
$$\hat{F}_t(\tilde{Y}_{1,1})=\tilde{Y}_{1,1}(1+t\tilde{A}_{1,2}^{-1}(1+t\tilde{A}_{2,4}^{-1}(1+t\tilde{A}_{1,6}^{-1})))$$
The product $\hat{F}_t(\tilde{Y}_{1,0})\hat{F}_t(\tilde{Y}_{1,1})$ has a unique dominant monomial $\tilde{Y}_{1,0}\tilde{Y}_{2,0}$, so:
$$\overset{\rightarrow}{E}_t(\tilde{Y}_{1,0}\tilde{Y}_{2,0})=\hat{F}_t(\tilde{Y}_{1,0}\tilde{Y}_{2,0})=\hat{L}_t(\tilde{Y}_{1,0}\tilde{Y}_{2,0})=\hat{F}_t(\tilde{Y}_{1,0})\hat{F}_t(\tilde{Y}_{1,1})$$
In particular the $V_{1,0}\otimes V_{1,1}$ is irreducible. Note that it is not a consequence of the conjecture but of classical theory of $q$-characters.

We suppose now that $s=5>4=2r^{\vee}$. There are two dominant monomials in $\tau_{s,t}(\overset{\rightarrow}{E}_t(\tilde{Y}_{1,0}\tilde{Y}_{1,1}))$: 
$$\tau_{s,t}(\tilde{Y}_{1,0}\tilde{Y}_{1,1})=\tilde{Y}_{1,0}\tilde{Y}_{1,1}\text{ and }\tau_{s,t}(t^3\tilde{Y}_{1,0}\tilde{A}_{1,1}^{-1}\tilde{A}_{2,3}^{-1}\tilde{A}_{1,5}^{-1}\tilde{Y}_{1,1})=t^{-1}$$
And so we have:
$$\tau_{s,t}(\overset{\rightarrow}{E}_t(\tilde{Y}_{1,0}\tilde{Y}_{1,1}))=\hat{L}_t^s(\tilde{Y}_{1,0}\tilde{Y}_{1,1})+t^{-1}\hat{L}_t(1)$$
where $\hat{L}_t(1)=1$. So if the conjecture is true, at $s=5$ the $V_{1,0}^s\otimes V_{1,1}^s$ is not irreducible and contains the trivial representation with multiplicity one.

\subsubsection{Non finite cases}\label{nonfini} In this section we suppose that $B(z)$ is symmetric and $s>2r^{\vee}$. An important difference with the finite case is that an infinite number of dominant monomials can appear in the $q,t$-character : let us briefly explain it for the example of section \ref{hatpi}. We consider the case $C$ of type $A_2^{(1)}$ and $s=3$. We have the following subgraph in the $q$-character given by the classical algorithm:
$$Y_{1,0}\rightarrow Y_{1,2}^{-1}Y_{2,1}Y_{3,1}\rightarrow Y_{3,2}Y_{3,1}Y_{2,3}^{-1}\rightarrow Y_{3,4}^{-1}Y_{3,1}Y_{1,0}$$
But at $s=3$ we have $Y_{3,4}^{-1}Y_{3,1}Y_{1,0}\simeq Y_{1,0}$. So we have a periodic chain and an infinity of dominant monomials in $\tau_{s,t}(\hat{F}_t(\tilde{Y}_{1,0}))$.

\noindent However we propose a construction of analogs of Kazdhan-Lusztig polynomials. As there is an infinity of dominant monomials, we have to begin the induction from the highest weight monomial. Let us describe it in a more formal way:

\noindent For $m\in \hat{B}^{s,\text{inv}}$ and $k\geq 0$ we denote by $\hat{B}^s_k(m)\subset \hat{B}^{s,\text{inv}}$ the set of dominant monomials of the form $m'=t^{\alpha}m\tilde{A}_{i_1,l_1}^{-1}...\tilde{A}_{i_k,l_k}^{-1}$. We set also $B^s(m)=\underset{k\geq 0}{\bigcup}B_k^s(m)$.

\noindent For $m\in\hat{B}^{s,\text{inv}}$, $\overset{\rightarrow}{E}_t(m)\in \hat{\mathfrak{K}}_t^{s,\infty}$ is defined as in section \ref{sblsubalg}. It will be useful to construct the element $\hat{F}_t^s(m)\in \hat{\mathfrak{K}}_t^{s,\infty}$ with a unique dominant monomial $m$: we denote by $m_0=m>m_1>m_2>...$ the dominant monomials appearing in $\overset{\rightarrow}{E}_t(m)$ with a total ordering compatible with the partial ordering and the degree (the set is countable because there is a finite number of monomials of degree $k$). We define $\lambda_k(t)\in\ZZ[t^{\pm}]$ inductively as the coefficient of $m_k$ in $\overset{\rightarrow}{E}_t(m)-\underset{1\leq l\leq k-1}{\sum}\lambda_l(t)\overset{\rightarrow}{E}_t(m_l)$. We define : 
$$\hat{F}_t^s(m)=\overset{\rightarrow}{E}_t(m)-\underset{l\geq 1}{\sum}\lambda_l(t)\overset{\rightarrow}{E}_t(m_l)\in \tilde{\mathfrak{K}}_t^{s,\infty}$$
(this infinite sum is allowed in $\tilde{\mathfrak{K}}_t^{s,\infty}$). The unique dominant monomial of $\hat{F}_t^s(m)$ is $m$. In particular $\overline{\hat{F}_t^s(m)}=\hat{F}_t^s(m)$ (see section \ref{sblsubalg}). In the following theorem the infinite sums are well-defined in $\tilde{\mathfrak{K}}_t^{s,\infty}$:

\begin{thm} For $m\in \hat{B}^{s,\text{inv}}$ there is a unique $\hat{L}_t^s(m)\in\tilde{\mathfrak{K}}_t^{s, \infty}$ of the form $\hat{L}_t^s(m)=m+\underset{m'<m}{\sum} \mu_{m',m}(t) m'$ such that:
$$\overline{\hat{L}_t^s(m)}=\hat{L}_t^s(m)$$
$$\overset{\rightarrow}{E}_t(m)=\hat{L}_t^s(m)+\underset{m'\in \hat{B}^s_k(m), k\geq 1}{\sum}P_{m',m}^s(t)\hat{L}_t^s(m')$$
where $P_{m',m}^s(t)\in t^{-1}\ZZ[t^{-1}]$. Moreover we have:
$$\hat{\Pi}(m)=\hat{\Pi}(m')\Rightarrow m^{-1}\hat{L}_t^s(m)={m'}^{-1}\hat{L}_t^s(m')$$
\end{thm}

\demo We aim at defining the $\mu_{m'm}(t)\in\ZZ[t^{\pm}]$ such that:
$$\hat{L}_t^s(m)=\underset{m'\in \hat{B}^s(m)}{\sum}\mu_{m',m}(t)\hat{F}_t(m')$$
The condition $\overline{\hat{L}_t^s(m)}=\hat{L}_t^s(m)$ means that $\mu_{m',m}(t^{-1})=\mu_{m',m}(t)$.

\noindent We define by induction on $k\geq 0$, for $m'\in \hat{B}^s_k(m)$ the $P^s_{m',m}(t)$ and the $\mu_{m',m}(t)$ such that: 
$$\hat{E}_t(m)-\underset{k\geq l\geq 0,m'\in \hat{B}^s_l(m)}{\sum}P^s_{m',m}(t)\underset{k\geq r\geq 0,m''\in \hat{B}^s_{r}(m')}{\sum}\mu_{m'',m'}(t)\hat{F}_t(m'')$$
$$\in \underset{m'\in \hat{B}^s_{k+1}(m)}{\sum}(\mu_{m',m}(t)+P^s_{m',m}(t))\hat{F}_t(m')+ \underset{l>k+1,m'\in \hat{B}^s_l(m')}{\sum}\ZZ[t^{\pm}] \hat{F}_t(m')$$ 
For $k=0$ we have $P_{m,m}^s(t)=\mu_{m,m}(t)=1$. And the the equation determines uniquely $P^s_{m',m}(t)\in t^{-1}\ZZ[t^{-1}]$ and $\mu_{m',m}(t)\in\ZZ[t^{\pm}]$ such that $\mu_{m',m}(t)=\mu_{m',m}(t^{-1})$.

\noindent For the last point we see also by induction on $k$ that for $m_1,m_2\in \hat{B}^{s , \text{inv}}$ such that $\hat{\Pi}(m_1)=\hat{\Pi}(m_2)$ and $m_1'\in \hat{B}^s(m_1)$, $m_2'\in\hat{B}^s(m_2)$ such that $m_1^{-1}m_1'\in t^{\ZZ}m_2^{-1}m_2'$ we have:
$$\mu_{m_1',m_1}(t)=\mu_{m_2',m_2}(t)\text{ , } P^s_{m_1',m_1}(t)=P^s_{m_2',m_2}(t)$$
\qed

\noindent Let us look at an example: we suppose that $C$ is of type $A_2^{(1)}$. In the generic case, the classical algorithm gives the $q$-characters beginning with $Y_{1,0}$, and the first terms are:
$$\xymatrix{&Y_{1,0} \ar[d]^{1,1}&
\\&Y_{1,2}^{-1}Y_{3,1}Y_{2,1}\ar[ld]^{3,2}\ar[rd]^{2,2}&
\\Y_{2,0}^{-1}Y_{3,1}Y_{3,2}&&Y_{3,3}^{-1}Y_{2,1}Y_{2,2}}$$
The deformed algorithm gives:
$$\overset{\rightarrow}{E}_t(\tilde{Y}_{1,0})=\tilde{Y}_{1,0}(1+t\tilde{A}_{1,1}^{-1}(1+t\tilde{A}_{2,2}^{-1}+t\tilde{A}_{3,2}))+\text{terms of higher degree}$$
We suppose now that $s=3$. First let $m=\tilde{Y}_{1,0}\tilde{Y}_{1,2}$, $m'=t^2\tilde{Y}_{1,0}\tilde{A}_{1,1}^{-1}\tilde{Y}_{1,2}$. We have:
$$\overset{\rightarrow}{E}_t(m)=\hat{F}_t(m)+t^{-1}\hat{F}_t(m')+...$$
In particular $P^s_{m',m}(t)=t^{-1}$.

\noindent Let $m=\tilde{Y}_{2,1}\tilde{Y}_{3,1}$, $m'=t\tilde{Y}_{2,1}\tilde{Y}_{3,1}\tilde{A}_{3,2}^{-1}\tilde{A}_{2,3}^{-1}$, $m''=t\tilde{Y}_{2,1}\tilde{Y}_{3,1}\tilde{A}_{2,2}^{-1}\tilde{A}_{3,3}^{-1}$. We have:
$$\overset{\rightarrow}{E}_t(m)=\hat{F}_t(m)+t^{-1}\hat{F}_t(m')+t^{-1}\hat{F}_t(m'')+...$$
In particular $P^s_{m',m}(t)=t^{-1}$ and $P^s_{m'',m}(t)=t^{-1}$.

Let us go back to general case and we want to define $P^s_{m',m}(t)$ for $m,m'\in B^s$. We can not set as in the finite case $P^s_{m',m}(t)=\underset{M'\in \hat{B}^s(M)/\hat{\Pi}(M')=m'}{\sum}P_{M',M}(t)$ (where $M\in\hat{B}^{s,\text{inv}}$ verifies $\hat{\Pi}(M)=m$) because this sum is not finite in general. However we propose the following construction. For $m,m'\in B^s$, we define $k(m,m')\geq 0$ such that for $M\in\hat{\Pi}^{-1}(m)$ we have $k(m,m')=\text{min}\{k\geq 0/ \exists M'\in \hat{B}_k^s(M),\hat{\Pi}(M')=m'\}$. 
\begin{defi} For $m,m'\in B^s$ we define $P^s_{m',m}(t)\in\ZZ[t^{\pm}]$ by:
$$P^s_{m',m}(t)=\underset{M'\in \hat{B}^s(M)/\hat{\Pi}(M')=m'\text{ and }deg(M')=deg(M)+k(m,m')}{\sum}P_{M',M}(t)$$
where $M$ an element of $\hat{B}^{s,\text{inv}}\cap \hat{\Pi}^{-1}(m)$.\end{defi}

Note that if $C$ affine it follows from lemma \ref{affroot} that for each $m\in B^s$, there is a finite number of $m'\in B^s$ such that $P^s_{m',m}(t)\neq 0$. In particular in this situation the proof of the theorem gives an algorithm to compute the polynomials with a finite number of steps (although there could be an infinite number of monomials in the $\epsilon,t$-character).

\noindent For example if $C$ is of type $A_2^{(1)}$ and $s=3$ we have:
$$P_{Y_{3,1}Y_{2,1},Y_{1,0}Y_{1,2}}(t)=t^{-1}\text{ , }P_{Y_{1,0},Y_{1,2},Y_{3,1}Y_{2,1}}=2t^{-1}$$

\subsection{Quantization of the Grothendieck ring}

\subsubsection{General quantization}

\noindent We set $\text{Rep}_t^s=\text{Rep}^s\otimes\ZZ[t^{\pm}]=\ZZ[X_{i,l},t^{\pm}]_{i\in I,l\in\ZZ/s\ZZ}$\label{repts} and we extend $\chi_{\epsilon,t}$ to a $\ZZ[t^{\pm}]$-linear injective map $\chi_{\epsilon,t}:\text{Rep}_t^s\rightarrow \tilde{\mathfrak{K}}_t^{\infty,s}$. We set $B^s=\{m=\underset{i\in I,l\in\ZZ/s\ZZ}{\prod}\tilde{Y}_{i,l}^{u_{i,l}(m)}\}\subset \overline{B}^s$. We have a map $\pi:\overline{B}^s\rightarrow B^s$ defined by $\pi(m)=\underset{i\in I,l\in\ZZ/s\ZZ}{\prod}\tilde{Y}_{i,l}^{u_{i,l}(m)}$.

\noindent We have:
$$\text{Im}(\chi_{\epsilon,t})=\underset{m\in B^s}{\bigoplus}\ZZ[t^{\pm}]\overset{\rightarrow}{E}_t(m)\subset \tilde{\mathfrak{K}}_t^{s,\infty}$$
But in general $\text{Im}(\chi_{\epsilon,t})$ is not a subalgebra of $\tilde{\mathfrak{K}}_t^{s,\infty}$. 

\noindent If $s=0$ or $C$ is finite we have $\text{Im}(\chi_{\epsilon,t})\subset \tilde{\mathfrak{K}}_t^{s,f,\infty}=\underset{m\in \overline{B}^s}{\bigoplus}\ZZ[t^{\pm}]\overset{\rightarrow}{E}_t(m)$ and we have a $\ZZ[t^{\pm}]$-linear map $\pi:\tilde{\mathfrak{K}}_t^{s,f,\infty}\rightarrow\text{Im}(\chi_{q,t})$ such that for $m\in\overline{B}^s$:
$$\pi(\overset{\rightarrow}{E}_t(m))=\overset{\rightarrow}{E}_t(\pi(m))$$

\noindent If $s>2r^{\vee}$ and $C$ verifies the property of lemma \ref{affroot} (for example $C$ is affine) then there is a $\ZZ[t^{\pm}]$-linear map $\pi:\tilde{\mathfrak{K}}_t^{s,\infty}\rightarrow\text{Im}(\chi_{q,t})$ such that for $m\in\overline{B}^s$ of the form $m=Mm'$ where $M\in B^s$ and $m'=t^{\alpha}\tilde{A}_{i_1,l_1}^{-1}...\tilde{A}^{-1}_{i_k,l_k}$ (see the definition of $k(m_1,m_2)\in\ZZ$ in section \ref{nonfini}):
$$\pi(\overset{\rightarrow}{E}_t(m))=\overset{\rightarrow}{E}_t(\pi(m))\text{ if $k=k(\hat{\Pi}(m),\hat{\Pi}(M))$}$$
$$\pi(\overset{\rightarrow}{E}_t(m))=0 \text{ if $k>k(\hat{\Pi}(m),\hat{\Pi}(M))$}$$

\noindent In both cases, as $\chi_{\epsilon,t}$ is injective, we can define a $\ZZ[t^{\pm}]$-bilinear map $*$ such that for $\alpha,\beta\in\text{Rep}_t^s$\label{star}:
$$\alpha*\beta=\chi_{q,t}^{-1}(\pi(\chi_{q,t}(\alpha)\chi_{q,t}(\beta)))$$
This is a deformed multiplication on $\text{Rep}_t^s$. But in general this multiplication is not associative.

\subsubsection{Associative quantization} In some cases it is possible to define an associative quantization (see \cite{Vas}, \cite{Nab}, \cite{her}). The point is to use a $t$-deformed algebra $\Yim_t=\ZZ[\tilde{Y}_{i,l}^{\pm},t^{\pm}]_{i\in I,l\in\ZZ}$ instead of $\hat{\Yim}_t$: in this case $\text{Im}(\chi_{q,t})$ is an algebra and we have an associative quantization of the Grothendieck ring (see \cite{her} for details). In this section we see how this construction can be generalized to other Cartan matrices. We suppose that $s=0$ and that $q$ is transcendental.

\begin{lem}\label{condcinv} Let $C$ be a Cartan matrix such that:
$$C_{i,j}<-1\Rightarrow -C_{j,i}\leq r_i$$
Then:
$$\text{det}(C(z))=z^{-R}+\alpha_{-R+1}z^{-R+1}+...+\alpha_{R-1}z^{R-1}+z^{R}$$ 
where $R=\underset{i=1...n}{\sum}r_i$ and $\alpha(-l)=\alpha(l)\in\ZZ$.\end{lem}

\noindent In particular finite and affine Cartan matrices ($A_1^{(1)}$ with $r_1=r_2=2$) verify the property of lemma \ref{condcinv}. Note that the condition $C_{i,j}<0\Rightarrow C_{i,j}=-1\text{ or }C_{j,i}=-1$ is sufficient; in particular Cartan matrices such that $i\neq j\Rightarrow C_{i,j}C_{j,i}\leq 3$ verify the property.

\demo For $\sigma\in S_n$ let us look at the term $\text{det}_{\sigma}=\underset{i\in I}{\prod}C_{i,\sigma (i)}(z)$ of $\text{det}(C(z))$. If $\sigma=\text{Id}$ then the degree $\text{deg}(\text{det}_{\text{Id}})$ is $\underset{i\in I}{\sum}r_i$. So it suffices to show that for $\sigma\neq \text{Id}$ we have $\text{deg}(\text{det}_{\sigma})<\underset{i\in I}{\sum}r_i$. If $i\neq \sigma(i)$, we have the following cases:

if $C_{i,\sigma(i)}=0$ or $-1$, $\text{deg}([C_{i,\sigma(i)}]_{z})\leq 0<r_{\sigma(i)}$

if $C_{i,\sigma(i)}<-1$, we have $C_{\sigma(i),i}=-1$ and so $r_iC_{i,\sigma(i)}=-r_{\sigma(i)}$ and so
$$\text{deg}([C_{i,\sigma(i)}]_{z})=-C_{i,\sigma(i)}-1=-\frac{r_{\sigma(i)}C_{\sigma (i),i}}{r_i}-1\leq r_{\sigma (i)}-1<r_{\sigma(i)}$$
So if $\sigma\neq \text{Id}$ we have:
$$\text{deg}(\text{det}_{\sigma})=\underset{i\in I/i=\sigma(i)}{\sum}r_i+\underset{i\in I/i\neq \sigma(i)}{\sum}\text{deg}([C_{i,\sigma(i)}]_{z_i})<\underset{i\in I/i=\sigma(i)}{\sum}r_i+\underset{i\in I/i\neq \sigma(i)}{\sum}r_{\sigma(i)}=\underset{i\in I}{\sum}r_i$$
For the last point $\text{det}(C(z))$ is symmetric polynomial because the coefficients of $C(z)$ are symmetric.
\qed

\noindent We suppose in this section that $C$ verifies the property of lemma \ref{condcinv}.

\noindent In particular $\text{det}(C(z))\neq 0$ and $C(z)$ has an inverse $\tilde{C}(z)$\label{invcar} with coefficients of the form $\frac{P(z)}{Q(z^{-1})}$ where $P(z)\in\ZZ[z^{\pm}]$, $Q(z)\in\ZZ[z]$, $Q(0)=\pm 1$ and the dominant coefficient of $Q$ is $\pm 1$. We denote by $\mathfrak{V}\subset\ZZ((z^{-1}))$ the set of rational fractions of this form. Note that $\mathfrak{V}$ is a subring of $\QQ(z)$, and for $R(z)\in\mathfrak{V},m\in\ZZ$ we have $R(z^m)\in\mathfrak{V}$. 
\noindent In particular for $m\in\ZZ-\{0\}$, $\tilde{C}(q^m)$ makes sense.

\noindent We denote by $\ZZ((z^{-1}))$ the ring of series of the form $P=\underset{r\leq R_P}{\sum}P_rz^r$ where $R_P\in\ZZ$ and the coefficients $P_r\in\ZZ$. We have an embedding $\mathfrak{V}\subset\ZZ((z^{-1}))$ by expanding $\frac{1}{Q(z^{-1})}$ in $\ZZ[[z^{-1}]]$ for $Q(z)\in\ZZ[z]$ such that $Q(0)=1$. So we can introduce maps ($\pi_r$, $r\in\ZZ$)\label{pir}:
$$\pi_r:\mathfrak{V}\rightarrow \ZZ\text{ , }P=\underset{r\leq R_P}{\sum}P_rz^r\mapsto P_r$$

\noindent We denote by $\mathcal{H}$ the algebra with generators $a_i[m],y_i[m],c_r$, relations \ref{aa}, \ref{ay} (of definition \ref{hatcalh}) and ($j\in I, m\neq 0$):
\begin{equation}\label{yeqa}y_j[m]=\underset{i\in I}{\sum}\tilde{C}_{i,j}(q^m)a_i[m]\end{equation}
\noindent Note that the relations \ref{yeqa} are compatible with the relations \ref{ay}.

\noindent We define $\Yim_u$ as the subalgebra of $\mathcal{H}[[h]]$ generated by the $\tilde{Y}_{i,l}^{\pm}$, $\tilde{A}_{i,l}^{\pm}$ ($i\in I,l\in\ZZ$), $t_R$ ($R\in\mathfrak{V}$).

\noindent Let the algebra $\Yim_t$ be the quotient of $\Yim_u$ by relations 
$$t_R=t_{R'}\text{ if $\pi_0(R)=\pi_0(R')$}$$
We keep the notations $\tilde{Y}_{i,l}^{\pm},\tilde{A}_{i,l}^{\pm}$ for their image in $\Yim_t$. We denote by $t$ the image of $t_1=\text{exp}(\underset{m>0}{\sum}h^{2m}c_m)$ in $\Yim_t$.

\noindent The following theorem is a generalization of theorem 3.11 of \cite{her}:

\begin{thm}(\cite{her})\label{dessusinv} The algebra $\Yim_t$ is defined by generators $\tilde{Y}_{i,l}^{\pm}$, $(i\in I,l\in\ZZ)$ central elements $t^{\pm}$  and relations ($i,j\in I, k,l\in\ZZ$)\label{gamma}:
$$\tilde{Y}_{i,l}\tilde{Y}_{j,k}=t^{\gamma(i,l,j,k)}\tilde{Y}_{j,k}\tilde{Y}_{i,l}$$
where $\gamma: (I\times\ZZ)^2\rightarrow\ZZ$ is given by:
$$\gamma(i,l,j,k)=\underset{r\in\ZZ}{\sum}\pi_r(\tilde{C}_{j,i}(z))(-\delta_{l-k,-r_j-r}-\delta_{l-k,r-r_j}+\delta_{l-k,r_j-r}+\delta_{l-k,r_j+r})$$
\end{thm}

\section{Complements}\label{complements}

\subsection{Finiteness of algorithms}

\noindent In the construction of $q,t$ and $\epsilon,t$-character we deal with completed algebras $\hat{\Yim}_t^{s,\infty}$, so the algorithms can produce an infinite number of monomials. In some cases we can say when this number is finite:

\subsubsection{Finiteness of the classical and deformed algorithms}

\begin{defi} We say that the classical algorithm stops if the classical algorithm is well defined and for all $m\in B$, $F(m)\in\mathfrak{K}$.\end{defi}

\noindent It follows from the classical theory of $q$-characters that if $C$ is finite then the classical algorithm stops.

\noindent For $i\in I$ let $L_i=(C_{i,1},...,C_{i,n})$. 

\begin{prop}\label{qtfinchar} We suppose that there are $(\alpha_i)_{i\in I}\in\ZZ^I$ such that $\alpha_i>0$ and:
$$\underset{j\in I}{\sum}\alpha_j L_j=0$$
Then the classical algorithm does not stop.\end{prop}

\noindent In particular if $C$ is an affine Cartan matrix then the classical algorithm does not stop.

\demo It follows from lemma \ref{ad} at $t=1$ that it suffices to show that there is no antidominant monomial in $C(Y_{1,0})$. So let $m=Y_{1,0}\underset{i\in I,l\in\ZZ}{\prod}A_{i,l}^{-v_{i,l}}$ be in $C(Y_{1,0})$. We see as in lemma \ref{affroot} $u_{i}(Y_{1,0}^{-1}m)=0$. In particular $u_1(m)=1$ and $m$ is not antidominant.\qed

\noindent Note that in the $A_r^{(1)}$-case ($r\geq 2$) we have a more ``intuitive'' proof : for all $l\in\ZZ$, $i\in I$ we have $A_{i,l}^{-1}=Y_{i,l+1}^{-1}Y_{i,l-1}^{-1}Y_{i+1,l}Y_{i-1,l}$, and: 
$$u(A_{i,l}^{-1})=\underset{j\in I,k\in\ZZ}{\sum}u_{j,k}(A_{i,l}^{-1})=(-u_{i,l+1}-u_{i,l-1}+u_{i+1,l}+u_{i-1,l})(A_{i,l}^{-1})=0$$
where we set in $I$: $(1)-1=r+1$ and $(r+1)+1=1$. So for all $m\in C(Y_{1,0})$ we have $u(m)=1$ and $m$ is not antidominant.

\subsubsection{Finiteness of the deformed algorithm}

\begin{prop}\label{propqtf} The following properties are equivalent:

i) For all $i\in I$, $\hat{F}_t(\tilde{Y}_{i,0})\in\hat{\mathfrak{K}}_t$.

ii) For all $m\in\overline{B}$, $\hat{F}_t(m)\in\hat{\mathfrak{K}}_t$.

iii) $\text{Im}(\chi_{q,t})\subset \hat{\mathfrak{K}}_t$.
\end{prop}

\begin{defi} If the properties of the proposition \ref{propqtf} are verified we say that the deformed algorithm stops.\end{defi}

\noindent Let us give some examples:

-If $C$ is of type $ADE$ then the deformed algorithm stops: \cite{Nab} (geometric proof) and \cite{Nac} (algebraic proof in $AD$ cases)

-If $C$ is of rank $2$ ($A_1\times A_1$, $A_2$, $B_2$, $C_2$, $G_2$) then the deformed algorithm stops: \cite{her} (algebraic proof)

-In \cite{her} we give an alternative algebraic proof for Cartan matrices of type $A_n$ ($n\geq 1$) and we conjecture that for all finite Cartan matrices the deformed algorithm stops. The cases $F_4$, $B_n$, $C_n$ ($n\leq 10$) have been checked on a computer (with the help of T. Schedler).

\begin{lem}\label{oneimp} If the deformed algorithm stops then the classical algorithm stops.\end{lem}

\demo This is a consequence of the formula $F(\hat{\Pi}(m))=\hat{\Pi}(\hat{F}_t(m))$ (see section \ref{rapalgo}).\qed

\noindent In particular if $C$ is affine then the deformed algorithm does not stop.

\noindent Let $C$ be a Cartan matrix such that $i\neq j\Rightarrow C_{i,j}C_{j,i}\leq 3$. We conjecture that the deformed algorithm stops if and only if the classical algorithm stops.

\subsection{$q,t$-characters of affine type and quantum toroidal algebras}

We have seen in \cite{her} that if $C$ is finite then the defining relations of $\hat{\mathcal{H}}$:
$$[a_i[m],a_j[r]]=\delta_{m,-r}(q^m-q^{-m})B_{i,j}(q^m)c_{|m|}$$
appear in the $\CC$-subalgebra $\U_q(\hat{\Hlie})$ of $\U_q(\hat{\Glie})$ generated by the $h_{i,m},c^{\pm}$ ($i\in I, m\in\ZZ-\{0\}$): it suffices to send $a_i[m]$ to $(q-q^{-1})h_{i,m}$ and $c_r$ to $\frac{c^{r}-c^{-r}}{r}$.

\noindent In this section we see that in the affine case $A_n^{(1)}$ ($n\geq 2$) the relations of $\hat{\mathcal{H}}$ appear in the structure of the quantum toroidal algebra. In particular we hope that $q,t$-characters will play a role in representation theory of quantum toroidal algebras (see the introduction).

\subsubsection{Reminder on quantum toroidal algebras \cite{var}} Let be $d\in\CC^*$ and $n\geq 3$. In the quantum toroidal algebra of type $sl_n$ there is a subalgebra $\mathcal{Z}$ generated by the $k_i^{\pm}, h_{i,l}$ ($i\in \{1,..., n\}$, $l\in\ZZ-\{0\}$) with relations :
$$k_ik_i^{-1}=cc^{-1}=1\text{ , }[k_{\pm,i}(z),k_{\pm,j}(w)]=0$$
\begin{equation}\label{kikj}\theta_{-a_{i,j}}(c^2 d^{-m_{i,j}}wz^{-1})k_{+, i}(z)k_{-, j}(w)=\theta_{-a_{i,j}}(c^{-2}d^{-m_{i,j}}wz^{-1})k_{-,j}k_{+,i}(z)\end{equation}
where $k_{\pm,i}(z)\in\mathcal{Z}[[z]]$ is defined by:
$$k_{\pm,i}(z)=k_i^{\pm}\text{exp}(\pm (q-q^{-1})\underset{k>0}{\sum}h_{i,\pm k}z^{\mp k})$$
$\theta_m(z)\in\CC[[z]]$ is the expansion of $\frac{q^mz-1}{z-q^m}$, $A=(a_{i,j})_{0\leq i,j\leq n}$ is the affine Cartan matrix of type $A_{n-1}^{(1)}$:
$$A=\begin{pmatrix}2 & -1 & ... & 0 & -1
\\-1 & 2 & ... & 0 & 0
\\ & &\ddots & &
\\0 & 0& ... &2 & -1
\\-1 & 0& ... &-1& 2\end{pmatrix}$$
and $M=(m_{i,j})_{1\leq i,j\leq n}$ is given by:
$$M=\begin{pmatrix}0 & -1 & ... & 0 &1
\\1 & 0 & ... & 0 & 0
\\ & &\ddots & &
\\0 & 0& ... &0 & -1
\\-1 & 0& ... & 1& 0 \end{pmatrix}$$

\subsubsection{Relations of the Heisenberg algebra}

\begin{lem}\label{lemrel} The relation (\ref{kikj}) are consequences of:
$$[h_{i,l},h_{j,m}]=\delta_{l,-m}\frac{q^{l a_{i,j}}-q^{-la_{i,j}}}{(q-q^{-1})^2}d^{-|l|m_{i,j}}\frac{c^{2l}-c^{-2l}}{l}$$
\end{lem}

\demo First for $m\in\ZZ$, we have in $\CC[[z]]$:
$$\theta_m(z)=\frac{q^m z-1}{z-q^m}=q^{-m}\text{exp}(\text{ln}(1-q^mz)-\text{ln}(1-q^{-m}z))=q^{-m}\text{exp}(\underset{r\geq 1}{\sum}(-\frac{(q^mz)^r}{r}+\frac{(q^{-m}z)^r}{r}))$$
and so $k_{+,i}(z)k_{-,j}(z)k_{+,i}(z)^{-1}k_{-,j}(w)^{-1}=\theta_{-a_{i,j}}(c^{-2}d^{-m_{i,j}}wz^{-1})\theta_{-a_{i,j}}(c^2d^{-m_{i,j}}wz^{-1})^{-1}$ is given by:
$$\theta_{-a_{i,j}}(c^{-2}d^{-m_{i,j}}wz^{-1})\theta_{-a_{i,j}}(c^2d^{-m_{i,j}}wz^{-1})^{-1}=\text{exp}(\underset{r\geq 1}{\sum}(q^{-ra_{i,j}}-q^{ra_{i,j}})d^{-rm_{i,j}}(wz^{-1})^r\frac{c^{2r}-c^{-2r}}{r}))$$
But following the proof of lemma \ref{reltr} we see that the relation of lemma \ref{lemrel} give:
$$[-(q-q^{-1})\underset{r\geq 1}{\sum}h_{j,-r}w^r,-(q-q^{-1})\underset{l\geq 1}{\sum}h_{i,l}z^{-r}]=\underset{r\geq 1}{\sum}(q-q^{-1})^2\frac{q^{-ra_{i,j}}-q^{ra_{i,j}}}{(q-q^{-1})^2}d^{-rm_{i,j}}(wz^{-1})^r\frac{c^{2r}-c^{-2r}}{r}$$
\qed

\noindent In particular for $d=1$, $a_i[m]=\frac{h_{i,m}}{q-q^{-1}}$ and $c_m=\frac{c^{2m}-c^{-2m}}{m}$, we get the defining relation of the Heisenberg algebra $\hat{\mathcal{H}}$ of section \ref{defihcal} in the affine case $A_{n-1}^{(1)}$:
$$[a_i[m],a_j[r]]=\delta_{m,-r}(q^m-q^{-m})[B_{i,j}]_{q^m}c_{|m|}$$
In the case $d\neq 1$ we have to extend the former construction:

\subsubsection{Twisted multiplication with two variables} Let us study the case $d\neq 1$: in this section we suppose that $q,d$ are indeterminate and we construct a $t$-deformation of $\ZZ[\tilde{A}_{i,l,k}^{\pm}]_{i\in I,l,k\in\ZZ}$. 

\noindent We define the $\CC[q^{\pm},d^{\pm}]$-algebra $\hat{\mathcal{H}}_d$ by generators $a_i[m]$ ($i\in I=\{1,...,n\}, m\in\ZZ$) and relations:
$$[a_i[m],a_j[r]]=\delta_{m,-r}(q^m-q^{-m})[A_{i,j}]_{q^m}d^{-|m|m_{i,j}}c_{|m|}$$
For $i\in I,l,k\in\ZZ$ we define $\tilde{A}_{i,l,k}\in\hat{\mathcal{H}}_d[[h]]$ by:
$$\tilde{A}_{i,l,k}=\text{exp}(\underset{m>0}{\sum}h^m q^{lm} d^{km}a_{i}[m])\text{exp}(\underset{m>0}{\sum}h^m q^{-lm}d^{-km}a_i[-m])$$
and for $R(q,d)\in\ZZ[q^{\pm},d^{\pm}]$, $t_R\in\hat{\mathcal{H}}_d[[h]]$ by:
$$t_{R(q,d)}=\text{exp}(\underset{m>0}{\sum}h^m R(q^m,d^m)c_m)$$
A computation analogous to the proof of lemma \ref{reltr} gives: 
$$\tilde{A}_{i,l,p}\tilde{A}_{j,k,r}\tilde{A}_{i,l,p}^{-1}\tilde{A}_{j,k,r}^{-1}=t_{(q-q^{-1})[A_{i,j}]_q(-q^{l-k}d^{p-r}+q^{k-l}d^{r-p})d^{-m_{i,j}}}$$
In particular, in the quotient of $\hat{\mathcal{H}}_d[[h]]$ by relations $t_R=1$ if $R\neq 0$, we have:
$$\tilde{A}_{i,l,p}^{-1}\tilde{A}_{j,k,r}^{-1}=t^{\alpha(i,j,k,l,p,r)}\tilde{A}_{j,k,r}^{-1}\tilde{A}_{i,l,p}^{-1}$$
where $\alpha: (I\times\ZZ\times\ZZ)^2\rightarrow\ZZ$ is given by ($l,k\in\ZZ$, $i,j\in I$):
$$\alpha(i,i,l,k,p,r)=2(\delta_{l-k,2r_i}-\delta_{l-k,-2r_i})\delta_{r,p}$$
$$\alpha(i,j,l,k,p,r)=\underset{r=r_iC_{i,j}+1,r_iC_{i,j}+3,...,-r_iC_{i,j}-1}{\sum}(-\delta_{l-k,r+r_i}\delta_{p-r,m_{i,j}}+\delta_{l-k,r-r_i}\delta_{r-p,m_{i,j}})\text{ (if $i\neq j$)}$$
In particular this would lead to the construction of $q,t$-characters with variables $\tilde{Y}_{i,l,p},\tilde{A}_{i,l,p}^{-1}$ associated to quantum toroidal algebras. But we shall leave further discussion of this point to another place.

\subsection{Combinatorics of bicharacters and Cartan matrices}\label{cartancompl}

In this section $C=(C_{i,j})_{1\leq i,j\leq n}$ is an indecomposable generalized (non necessarily symmetrizable) Cartan matrix and $(r_1,...,r_n)$ are positive integers. Let $D=\text{diag}(r_1,...,r_n)$ and $B=DC$ (which is non necessarily symmetric).

\noindent We show that the quantization of $\hat{\Yim}^s\otimes\ZZ[t^{\pm}]=\ZZ[Y_{i,l},V_{i,l},t^{\pm}]_{i\in I,l\in\ZZ/s\ZZ}$ is linked to fundamental combinatorial properties of $C$ and $(r_1,...,r_n)$ (propositions \ref{facile}, \ref{faciledeux}, \ref{qsym} and theorem \ref{result}). Let us begin with some general background about twisted multiplication defined by bicharacters.

\subsubsection{Bicharacters and twisted multiplication} Let $\Lambda$ be a set, $Y$ be the commutative polynomial ring:
$$Y=\ZZ[X_{\alpha},t^{\pm}]_{\alpha\in\Lambda}$$
and $A$ the set of monomials of the form $m=\underset{\alpha\in\Lambda}{\prod}X_{\alpha}^{x_{\alpha}(m)}\in Y$. The usual commutative multiplication of $Y$ is denoted by $.$ in the following.

\begin{defi} A bicharacter on $A$ is a map $d:A\times A\rightarrow\ZZ$ such that ($m_1,m_2,m_3\in A$):
$$d(m_1.m_2,m_3)=d(m_1,m_3)+d(m_2,m_3)\text{ , }d(m_1,m_2.m_3)=d(m_1,m_2)+d(m_1,m_3)$$
\end{defi}

\noindent The symmetric bicharacter $\mathfrak{S}d$ and the antisymmetric bicharacter $\mathfrak{A}d$ of $d$ are defined by:
$$\mathfrak{S}d(m_1,m_2)=\frac{1}{2}(d(m_1,m_2)+d(m_2,m_1))\text{ , }\mathfrak{A}d(m_1,m_2)=\frac{1}{2}(d(m_1,m_2)-d(m_2,m_1))$$
and we have $d=\mathfrak{A}d+\mathfrak{S}d$.

\noindent Let be $d$ be a bicharacter on $A$. One can define a $\ZZ[t^{\pm}]$-bilinear map $*:Y\times Y\rightarrow Y$ such that: 
$$m_1*m_2=t^{d(m_1,m_2)}m_1.m_2$$
This map is associative\footnote{In fact it suffices that $-d(m_2,m_3)+d(m_1m_2,m_3)-d(m_1,m_2m_3)+d(m_1,m_2)=0$.} and we get a $\ZZ[t^{\pm}]$-algebra structure on $Y$. We say that the new multiplication is the twisted multiplication associated to the bicharacter $d$, and it is given by formulas:
$$m_1*m_2=t^{d(m_1,m_2)-d(m_2,m_1)}m_2*m_1=t^{2\mathfrak{A}d(m_1,m_2)}m_2*m_1$$

\begin{lem}\label{wddeux} Let $d_1,d_2$ be two bicharacters. One can define a multiplication on $Y$ such that ($m_1,m_2\in A$):
$$m_1*m_2=t^{2d_1(m_1,m_2)-2d_2(m_2,m_1)}m_2*m_1$$
if and only if $\mathfrak{S}d_1=\mathfrak{S}d_2$.

\noindent In this case, the multiplication is the twisted multiplication associated to the bicharacter $d=d_1+d_2$:
$$m_1*m_2=t^{d_1(m_1,m_2)+d_2(m_1,m_2)}m_1.m_2$$
\end{lem}

\demo It follows immediately from the definition of $*$:
$$m_1*m_2=t^{2d_1(m_1,m_2)-2d_2(m_2,m_1)}m_2*m_1=t^{4(\mathfrak{S}d_1-\mathfrak{S}d_2)(m_1,m_2)}m_1*m_2$$
If $\mathfrak{S}d_1=\mathfrak{S}d_2$, let $*$ be the twisted multiplication associated with the bicharacter $d=d_1+d_2$. We have:
$$m_1*m_2=t^{d_1(m_1,m_2)+d_2(m_1,m_2)-d_1(m_2,m_1)-d_2(m_2,m_1)}m_2*m_1=t^{2d_1(m_1,m_2)-2d_2(m_2,m_1)}m_2*m_1$$
\qed

\subsubsection{Definition of $d_1$ and $d_2$}\label{bicaddstu}

For $s\geq 0$ let $\Lambda_s=I\times (\ZZ/s\ZZ)$\label{lsg} and $\overline{A}^s$ be the set of monomials of $\hat{\Yim}^s$, that is to say elements of the form $m=\underset{(i,l)\in\Lambda_s}{\prod}Y_{i,l}^{y_{i,l}(m)}V_{i,l}^{v_{i,l}(m)}$. Let $D(z)=\text{diag}([r_1]_z,...,[r_n]_z)$.

\noindent For $\alpha\in\Lambda_s$, we define a character \footnote{ie. $u_{\alpha}(m_1.m_2)=u_{\alpha}(m_1)+u_{\alpha}(m_2)$} $u_{\alpha}$ on $\overline{A}^s$ as in section \ref{presbica}. In particular $u_{\alpha}(Y_{\beta})=\delta_{\alpha,\beta}$. 

\noindent We define $d_1,d_2$ the bicharacters on $\overline{A}^s$ as in section \ref{presbica}, that is to say ($m_1,m_2\in\overline{A}^s$):
$$d_1(m_1,m_2)=\underset{\alpha\in\Lambda_s}{\sum}v_{b(\alpha)}(m_1)u_{\alpha}(m_2)+y_{b(\alpha)}(m_1)v_{\alpha}(m_2)$$
$$d_2(m_1,m_2)=\underset{\alpha\in\Lambda_s}{\sum}u_{b(\alpha)}(m_1)v_{\alpha}(m_2)+v_{b(\alpha)}(m_1)y_{\alpha}(m_2)$$
where $b:\Lambda_s\rightarrow \Lambda_s$\label{ptb} is the bijection defined by $b(i,l)=(i,l+r_i)$.

\begin{prop}\label{facile} The following properties are equivalent:

i) For $s\geq 0$, $d_1=d_2$ 

ii) For $s\geq 0$, $\forall \alpha,\beta\in\Lambda_s, u_{\alpha}(V_{\beta})=u_{b(\beta)}(V_{b(\alpha)})$

iii) $C$ is symmetric and $\forall i,j\in I$, $r_i=r_j$.

\end{prop}

\demo We have always:
$$d_1(Y_{\alpha},Y_{\beta})=d_2(Y_{\alpha},Y_{\beta})=0$$
For $\alpha,\beta\in\Lambda_s$, we have $u_{\alpha}(Y_{\beta})=\delta_{\alpha,\beta}$. In particular:
$$d_1(Y_{\alpha},V_{\beta})=\delta_{b(\beta),\alpha}=u_{b(\beta)}(Y_{\alpha})=d_2(Y_{\alpha},V_{\beta})$$
$$d_1(V_{\beta},Y_{\alpha})=u_{b^{-1}(\beta)}(Y_{\alpha})=\delta_{b(\alpha),\beta}=d_2(V_{\beta},Y_{\alpha})$$
So the condition $d_1=d_2$ means $\forall \alpha,\beta\in\Lambda_s$, $d_1(V_{\alpha},V_{\beta})=d_2(V_{\alpha},V_{\beta})$. But the equation (ii)  means:
$$d_1(V_{\alpha},V_{\beta})=u_{b^{-1}(\alpha)}(V_{\beta})=u_{b(\beta)}(V_{\alpha})=d_2(V_{\alpha},V_{\beta})$$
In particular we have $(i)\Leftrightarrow (ii)$.

\noindent For $i,j\in I$ and $l,k\in\ZZ/s\ZZ$ we have:
$$u_{i,l}(V_{j,k})=\underset{r=C_{i,j}+1... -C_{i,j}-1}{\sum}\delta_{l+r,k}=\underset{r=C_{i,j}+1... -C_{i,j}-1}{\sum}\delta_{l-k,r}$$
$$u_{j,k+r_j}(V_{i,l+r_i})=\underset{r=C_{j,i}+1... -C_{j,i}-1}{\sum}\delta_{k+r+r_j,l+r_i}=\underset{r=C_{j,i}+1... -C_{j,i}-1}{\sum}\delta_{l-k,r_j-r_i+r}$$
If $s=0$, those terms are equal for all $l,k\in\ZZ$ if and only if $C_{i,j}\neq 0$ implies $C_{i,j}=C_{j,i}$ and $r_i=r_j$. So as $C$ is indecomposable we have $(ii)\Leftrightarrow (iii)$.

\noindent If $s\geq 0$ and $(iii)$ is verified we see in the same way that those terms are equal, so $(iii)\Rightarrow (ii)$.\qed

\noindent In particular if $C$ is of type $ADE$, we get the bicharacter of \cite{Nab} and $d_1=d_2$ is the equation (\cite{Nab}, 2.1). 

\subsubsection{Bicharacters and symmetrizable Cartan matrices}\label{andsym}

We have seen in lemma \ref{wddeux} that we can define a twisted multiplication if and only if $\mathfrak{S}d_1=\mathfrak{S}d_2$, so we investigate those cases:

\begin{thm}\label{result} The following properties are equivalent:

i) For $s\geq 0$, we have $\mathfrak{S}d_1=\mathfrak{S}d_2$

ii) For $s\geq 0$, $\forall\alpha,\beta\in\Lambda_s, u_{\alpha}(V_{b(\beta)})-u_{b^2(\alpha)}(V_{b(\beta)})=u_{b^2(\beta)}(V_{b(\alpha)})-u_{\beta}(V_{b(\alpha)})$

iii) For $s\geq 0$ and $m\in \overline{A}^s$, $d_1(m,m)=d_2(m,m)$

iv) $B(z)$ is symmetric

v) $B$ is symmetric and $C_{i,j}\neq C_{j,i}\Longrightarrow (r_i=-C_{j,i}\text{ and }r_j=-C_{i,j})$ \end{thm}

\demo 

First we show that $(i)\Leftrightarrow (ii)$. We have always:
$$\mathfrak{S}d_1(Y_{\alpha},Y_{\beta})=\mathfrak{S}d_2(Y_{\alpha},Y_{\beta})=0$$
and:
$$2\mathfrak{S}d_1(Y_{\alpha},V_{\beta})=u_{b(\beta)}(Y_{\alpha})-u_{b^{-1}(\beta)}(Y_{\alpha})=\delta_{b(\beta),\alpha}-\delta_{b(\alpha),\beta}=2\mathfrak{S}d_2(Y_{\alpha},V_{\beta})$$
But the equation (ii) means:
$$d_1(V_{\alpha},V_{\beta})-d_2(V_{\beta},V_{\alpha})=d_2(V_{\alpha},V_{\beta})-d_1(V_{\beta},V_{\alpha})$$
that is to say:
$$2\mathfrak{S}d_1(V_{\alpha},V_{\beta})=2\mathfrak{S}d_2(V_{\alpha},V_{\beta})$$
and we can conclude because $d_1,d_2$ are bicharacters.

Let us show that $(iv)\Leftrightarrow (v)$: the matrix $B(z)$ is symmetric if and only if for all $i\neq j$ we have:
$$(z^{r_i}-z^{-r_i})(z^{C_{i,j}}-z^{-C_{i,j}})=(z^{r_j}-z^{-r_j})(z^{C_{j,i}}-z^{-C_{j,i}})$$
If $C_{i,j}=C_{j,i}=0$ it is obvious. If $C_{i,j}=C_{j,i}\neq 0$, the equation means $r_i=r_j$. If $C_{i,j}\neq C_{j,i}$, the equality means $(r_i=-C_{j,i}\text{ and }r_j=-C_{i,j})$.

\noindent The equation $(ii)$ means:
$$\underset{r=C_{i,j}+1... -C_{i,j}-1}{\sum}\delta_{l-k,r_j-r}-\delta_{l-k,r_j-2r_i-r}=\underset{r=C_{j,i}+1... -C_{j,i}-1}{\sum}\delta_{l-k,2r_j+r-r_i}-\delta_{l-k,r-r_i}$$
At $s=0$, the formula holds for all $l,k\in\ZZ$, if and only the coefficients of Kronecker's functions are equal, that is to say in $\ZZ[X^{\pm}]$:
$$\underset{r=C_{i,j}+1... -C_{i,j}-1}{\sum}X^{r_j-r}-X^{r_j-2r_i-r}=\underset{r=C_{j,i}+1... -C_{j,i}-1}{\sum}X^{2r_j+r-r_i}-X^{r-r_i}$$
$$(X^{r_j}-X^{r_j-2r_i})X^{C_{i,j}+1}\frac{1-X^{-2C_{i,j}}}{1-X^{2}}=(X^{2r_j-r_i}-X^{-r_i})X^{C_{j,i}+1}\frac{1-X^{-2C_{j,i}}}{1-X^{2}}$$
$$X^{r_j-2r_i+C_{i,j}}(1-X^{-2C_{i,j}})(1-X^{2r_i})=X^{-r_i+C_{j,i}}(1-X^{-2C_{j,i}})(1-X^{2r_j})$$
$$\text{($C_{i,j}=C_{j,i}=0$) or ($r_i=r_j$ and $C_{i,j}=C_{j,i}\neq 0$) or ($r_j=-C_{i,j}$ and $r_i=-C_{j,i}$)}$$
and so $(ii)\Rightarrow (v)$. If we suppose that $(iv)$ is true, then the above equation is also verified in $\ZZ[X^{\pm}]/(X^s=1)$ and $(ii)$ is true.

\noindent To conclude it suffices to show that $(iii) \Leftrightarrow (i)$. If $(iii)$ is verified we have for $m,m'\in\overline{A}^s$:
$$d_1(m,m')+d_1(m',m)=d_1(mm',mm')-d_1(m,m)-d_1(m',m')=d_2(m,m')+d_2(m',m)$$
and $(i)$ is verified. If $(i)$ is verified we have for $m\in\overline{A}^s$: $2d_1(m,m)=2d_2(m,m)$.\qed

\subsubsection{Bicharacters and $q$-symmetrizable Cartan matrices}\label{bznonsym} There is a way to define a deformation multiplication if $B(z)$ is non necessarily symmetric. First we define the matrix $C_{i,j}'(z)=[C_{i,j}]_{z_i}$\label{czp} and the characters \label{uilp}:
$$u_{i,l}'(m)=y_{i,l}(m)-\underset{j\in I}{\sum}(C'_{i,j}(z))_{op}\mathcal{V}_{j,l}(m)$$

\noindent We define the bicharacters $d_1'$ and $d_2'$\label{dp} from $\tilde{u}_{i,l}$ in the same way $d_1$ and $d_2$ were defined from $u_{i,l}$ (section \ref{bicaddstu}). 

\noindent We also define $B_{i,j}'(z)=[B_{i,j}]_z$\label{bzp}. Note that we have always $B'(z)=D(z)C'(z)$. Indeed:
$$B_{i,j}'(z)=\frac{z^{r_iC_{i,j}}-z^{-r_iC_{i,j}}}{z-z^{-1}}=\frac{z_i-z_i^{-1}}{z-z^{-1}}\frac{z_i^{C_{i,j}}-z_i^{-C_{i,j}}}{z_i-z_i^{-1}}=D_{i,i}(z)C_{i,j}'(z)$$

\begin{prop}\label{faciledeux} The following properties are equivalent:

i) For $s\geq 0$, $\mathfrak{S}d_1'=\mathfrak{S}d_2'$

ii) For $s\geq 0$, $\forall\alpha,\beta\in\Lambda_s, u_{\alpha}'(V_{b(\beta)})-u_{b^2(\alpha)}'(V_{b(\beta)})=u_{b^2(\beta)}'(V_{b(\alpha)})-u_{\beta}'(V_{b(\alpha)})$

iii) $B$ is symmetric

iv) $B'(z)$ is symmetric \end{prop}

\noindent In particular if $C$ is symmetrizable we can define the deformed structure for all $s\geq 0$.

\demo First we have $(iii)\Leftrightarrow (iv)$ because $B_{i,j}'(z)=[B_{i,j}]_z$.

\noindent We show as in theorem \ref{result} that $(ii)\Leftrightarrow (i)$. 

\noindent Let us write the equation (ii):
$$u_{i,l}'(V_{j,k+r_j})-u_{i,l+2r_i}'(V_{j,k+r_j})=u_{j,k+2r_j}'(V_{i,l+r_i})-u_{j,k}'(V_{i,l+r_i})$$
If $i=j$, we are in the symmetric case, and it follows from proposition \ref{facile} that this equation is verified. In the case $i\neq j$, if $C_{i,j}=0$ then all is equal to $0$. In the cases $C_{i,j}<0$ the equation reads:
$$\underset{r=C_{i,j}+1... -C_{i,j}-1}{\sum}\delta_{l+r_ir,k+r_j}-\delta_{l+2r_i+rr_i,k+r_j}=\underset{l=C_{j,i}+1... -C_{j,i}-1}{\sum}\delta_{k+2r_j+lr_j,l+r_i}-\delta_{k+rr_j,l+r_i}$$
$$\underset{r=C_{i,j}+1... -C_{i,j}-1}{\sum}\delta_{l-k,r_j-rr_i}-\delta_{l-k,r_j-2r_i-r_ir}=\underset{r=C_{j,i}+1... -C_{j,i}-1}{\sum}\delta_{l-k,2r_j+rr_j-r_i}-\delta_{l-k,rr_j-r_i}$$
$$\delta_{l-k,r_j-r_i-r_iC_{i,j}}-\delta_{l-k,r_j-r_i+r_iC_{i,j}}=\delta_{l-k,r_j-r_i-r_jC_{j,i}}-\delta_{l-k,r_j-r_i+r_jC_{j,i}}$$
That is to say:
$$(2r_iC_{i,j}\in s\ZZ\text{ and }2r_jC_{j,i}\in s\ZZ)\text{ or }r_iC_{i,j}-r_jC_{j,i}\in s\ZZ$$
If $s=0$, the equation means $r_iC_{i,j}=r_jC_{j,i}$ that is to say $B=DC$ symmetric. So $(ii)\Leftrightarrow (iii)$.

\noindent If $s\geq 0$ and $B$ symmetric we have $r_iC_{i,j}-r_jC_{j,i}\in s\ZZ$. So $(iii)\Rightarrow (ii)$.\qed

\noindent In some situations the two constructions are the same:

\begin{prop}\label{qsym} The following properties are equivalent:

i) For $s\geq 0$, $u'=u$

ii) For $s\geq 0$, $d_1'=d_1$ 

iii) For $s\geq 0$, $d_2'=d_2$

iv) $C'(z)=C(z)$

v) $B'(z)=B(z)$

vi) $i\neq j\Rightarrow (r_i=1\text{ or }C_{i,j}=-1\text{ or }C_{i,j}=0)$

\end{prop}

\demo We have $(iv)\Leftrightarrow (v)$ because $B(z)=D(z)C(z)$, $B'(z)=D(z)C'(z)$ and $D(z)$ is invertible.

\noindent The $(i)\Rightarrow (ii)$ (resp. $(i)\Rightarrow (iii)$) is clear and we get $(ii)\Rightarrow (i)$ (resp. $(iii)\Rightarrow (i)$) by looking at $d_1(V_{i,l},V_{j,k})=d_1'(V_{i,l},V_{j,k})$ (resp. $d_2(V_{i,l},V_{j,k})=d_2'(V_{i,l},V_{j,k})$).

\noindent The $(iv)\Rightarrow (i)$ is clear. If $(i)$ is true we have for $i\neq j$ and all $l,k\in\ZZ$: 
$$u_{i,l}(V_{j,k})=\underset{r=C_{i,j}+1,C_{i,j}+3,...,-1C_{i,j}-1}{\sum}\delta_{l-k,r}=\underset{r=C_{i,j}+1,C_{i,j}+3,...,-1C_{i,j}-1}{\sum}\delta_{l-k,rr_i}=u_{i,l}'(V_{j,k})$$
and so $\frac{z^{C_{i,j}}-z^{-C_{i,j}}}{z-z^{-1}}=\frac{z_i^{C_{i,j}}-z_i^{-C_{i,j}}}{z_i-z_i^{-1}}$ that is to say $(iv)$.

\noindent So it suffices to show that $(v)\Leftrightarrow (vi)$. We have always:
$$B_{i,i}(z)=\frac{z^{r_i}-z^{-r_i}}{z-z^{-1}}(z^{r_i}+z^{-r_i})=\frac{z^{2r_i}-z^{-2r_i}}{z-z^{-1}}=[2r_i]_z=[B_{i,i}]_z$$
If $i\neq j$, the equality $B_{i,j}(z)=B'_{i,j}(z)$ means:
$$z^{r_i+C_{i,j}}+z^{-r_i-C_{i,j}}-z^{C_{i,j}-r_i}-z^{r_i-C_{i,j}}=z^{r_iC_{i,j}+1}+z^{-1-r_iC_{i,j}}-z^{r_iC_{i,j}-1}-z^{1-r_iC_{i,j}}$$
If $r_i=1$ or $C_{i,j}=-1$ or $C_{i,j}=0$ the equality is clear and so $(vi)\Rightarrow (v)$. Suppose that (v) is true and let be $i\neq j$. We have to study different cases:

$r_i+C_{i,j}=C_{i,j}-r_i\Rightarrow r_i=0$ (impossible)

$r_i+C_{i,j}=r_i-C_{i,j}\Rightarrow C_{i,j}=0$

$r_i+C_{i,j}=r_iC_{i,j}+1\text{ and }r_iC_{i,j}-1=C_{i,j}-r_i\Rightarrow r_i=1$

$r_i+C_{i,j}=r_iC_{i,j}+1\text{ and }-r_iC_{i,j}+1=C_{i,j}-r_i\Rightarrow C_{i,j}=1$ (impossible)

$r_i+C_{i,j}=-r_iC_{i,j}-1\text{ and }r_iC_{i,j}-1=C_{i,j}-r_i\Rightarrow C_{i,j}=-1$

$r_i+C_{i,j}=-r_iC_{i,j}-1\text{ and }-r_iC_{i,j}+1=C_{i,j}-r_i\Rightarrow r_i=-1$ (impossible)

\noindent and so we get (vi).\qed

\begin{lem}\label{defiqs} If the properties of the proposition \ref{qsym} are verified and $B=DC$ is symmetric then the properties of the proposition \ref{faciledeux} are verified.\end{lem}

\demo We verify the property (iv) of proposition \ref{faciledeux}: we suppose that $C_{i,j}\neq C_{j,i}$. So $C_{i,j}\neq 0$, $C_{j,i}\neq 0$ and we do not have $C_{i,j}=C_{j,i}=-1$. As $r_iC_{i,j}=r_jC_{j,i}$, we do not have $r_i=r_j=1$. So we have (property (vi) of proposition \ref{qsym}) $r_i=-C_{j,i}=1$ or $r_j=-C_{i,j}=1$. For example in the first case, $r_iC_{i,j}=r_jC_{j,i}$ leads to $C_{i,j}=-r_j$.\qed

\begin{defi} We say that $C$ is $q$-symmetrizable if $B=DC$ is symmetric and:
$$i\neq j\Rightarrow (r_i=1\text{ or }C_{i,j}=-1\text{ or }C_{i,j}=0)$$\end{defi}

\noindent In particular $C$ $q$-symmetrizable verifies the properties of proposition \ref{faciledeux}, \ref{qsym} and of theorem \ref{result}.

\subsubsection{Examples}

If $C$ is symmetric then for all $i\in I$ we have $r_i=1$ and so $C$ is $q$-symmetrizable.

\begin{lem} The Cartan matrices of finite or affine type (except $A_1^{(1)}$, $A_{2l}^{(2)}$ case, $l\geq 2$) are $q$-symmetrizable. The affine Cartan matrices $A_1^{(1)}$, $A_{2l}^{(2)}$ with $l\geq 2$ are not $q$-symmetrizable.
\end{lem}

\noindent In particular if $C$ is finite then $u=\tilde{u}$ and the presentation adopted in this paper fits with former articles, in particular in the non symmetric cases (\cite{Fre}, \cite{Fre2}, \cite{Fre3}, \cite{her}).

\demo As those matrices are symmetrizable, it suffices to check the property (vi) of proposition \ref{qsym}:

the finite Cartan matrices $A_l$ ($l\geq 1$), $D_l$ ($l\geq 4$), $E_6$, $E_7$, $E_8$ and the affine Cartan matrices $A_l^{(1)}$ ($l\geq 1$), $D_l^{(1)}$ ($l\geq 4$), $E_6^{(1)}$, $E_7^{(1)}$, $E_8^{(1)}$ are symmetric and so $q$-symmetrizable.

the finite Cartan matrices $B_l$ ($l\geq 2$), $G_2$ and the affine Cartan matrices $B_l^{(1)}$ ($l\geq 3$), $G_2^{(1)}$ verify  $r_n=1$ and for $i\neq j$: $i\leq n-1\Rightarrow C_{i,j}=-1\text{ or }0$.

the finite Cartan matrices $C_l$ ($l\geq 2$), the affine Cartan matrices $A_{2l-1}^{(2)}$ ($l\geq 3$), $D_4^{(3)}$ verify $r_1=...=r_{n-1}=1$, $C_{n,1}=...=C_{n,n-2}=0$ and $C_{n,n-1}=-1$.

the affine Cartan matrices $C_l^{(1)}$ ($l\geq 2$) verify $r_2=...=r_{n-1}=1$ and $C_{1,3}=...=C_{1,n}=0$, $C_{1,2}=-1$, $C_{n,1}=...=C_{n,n-2}=0$, $C_{n,n-1}=-1$.

the affine Cartan matrices $D_{l+1}^{(2)}$ ($l\geq 2$) verify $r_1=r_n=1$ and for $i\neq j$: $2\leq i\leq n-1\Rightarrow C_{i,j}=-1\text{ or }0$.

\noindent The other particular cases are studied one after one:

for the finite Cartan matrix $F_4=\begin{pmatrix}2 &-1&0&0\\-1&2&-1&0\\0&-2&2&-1\\0&0&-1&2\end{pmatrix}$ we have $(2,2,1,1)$

for the affine Cartan matrix $F_4^{(1)}=\begin{pmatrix}2&-1&0&0&0\\-1&2 &-1&0&0\\0&-1&2&-1&0\\0&0&-2&2&-1\\0&0&0&-1&2\end{pmatrix}$ we have $(2,2,2,1,1)$

for the affine Cartan matrix $A_2^{(2)}=\begin{pmatrix}2&-4\\-1&2\end{pmatrix}$ we have $(1,4)$

for the affine Cartan matrix $E_6^{(2)}=\begin{pmatrix}2&-1&0&0&0\\-1&2 &-1&0&0\\0&-1&2&-2&0\\0&0&-1&2&-1\\0&0&0&-1&2\end{pmatrix}$ we have $(1,1,1,2,2)$.

\noindent Finally the affine Cartan matrices $A_1^{(1)}$ and $A_{2l}^{(2)}$ ($l\geq 2$) are not $q$-symmetrizable because $C_{n-1,n}=-2$ and $r_{n-1}=2$.
\qed

\noindent One can understand ``intuitively'' the fact that $A_{2l}^{(2)}$ ($l\geq 2$) is not $q$-symmetrizable: in the Dynkin diagram there is an oriented path without loop with two arrows in the same direction.

\noindent There are $q$-symmetrizable Cartan matrices which are not finite and not affine: here is an example such that for all $i,j\in I$, $C_{i,j}\geq -2$:
$$C=\begin{pmatrix}2 &-2&-2&0\\-1&2&0&-1\\-1&0&2&-1\\0&-2&-2&2\end{pmatrix}$$
$$(r_1,r_2,r_3,r_4)=(1,2,2,1)$$

\newpage

\twocolumn

\section*{Notations}

\begin{tabular}{lll}

$A^s$ & set of $\Yim^s$-monomials & p \pageref{as}

\\$\overline{A}, \overline{A}^s$ & sets of $\hat{\Yim}_t, \hat{\Yim}_t^s$-monomials & p \pageref{linea}

\\$\hat{A}^{\text{inv}}, \hat{B}^{\text{inv}}$ & sets of $\hat{\Yim}_t$-monomials & p \pageref{ainv}

\\$\hat{A}^{s, \text{inv}}, \hat{B}^{s, \text{inv}}$ & sets of $\hat{\Yim}_t^s$-monomials & p \pageref{ainvs}

\\$\overset{\rightarrow}{A}, \overset{\leftarrow}{A}$ & sets of $\hat{\Yim}_t$-monomials & p \pageref{aarrow}

\\$\alpha$ & map $(I\times \ZZ/s\ZZ)^2\rightarrow \ZZ$ & p \pageref{alpha}

\\$\alpha(m)$ & character & p \pageref{alpham}

\\$a_i[m]$ & element of $\hat{\mathcal{H}}$ & p \pageref{aim}

\\$\tilde{A}_{i,l}, \tilde{A}_{i,l}^{-1}$ & elements of $\hat{\Yim}_u$ or $\hat{\Yim}_t$ & p \pageref{tail}

\\$A_{i,l}, A_{i,l}^{-1}$ & elements of $\hat{\Yim}$ & p \pageref{ail}

\\ $b$ & bijection of $\Lambda_s$ & p \pageref{ptb}

\\ $\overline{B}_i, \overline{B}$ & sets of $\hat{\Yim}_t$-monomials & p \pageref{hata}

\\ $\overline{B}_i^s, \overline{B}^s$ & sets of $\hat{\Yim}_t^s$-monomials & p \pageref{hatas}

\\$B^s, B_i^s$ & sets of $\Yim^s$-monomials & p \pageref{as}

\\$B=(B_{i,j})$ & symmetrized &

\\ &Cartan matrix & p \pageref{symcar}

\\$B(z)$ & deformation of $B$ & p \pageref{bz}

\\$B'(z)$ & deformation of $B$ & p \pageref{bzp}

\\$\beta$ & map $(I\times \ZZ/s\ZZ)^2\rightarrow \ZZ$ & p \pageref{beta}

\\$C=(C_{i,j})$ & Cartan matrix & p \pageref{carmat}

\\$C(z)$ & deformation of $B$ & p \pageref{cz}

\\$C'(z)$ & deformation of $B$ & p \pageref{czp}

\\$C(m)$ & set of monomials & p \pageref{cpm}

\\$(\tilde{C}_{i,j})$ & inverse of $C$ & p \pageref{invcar}

\\$c_r$ & central element of $\hat{\mathcal{H}}$ & p \pageref{cr}

\\$d_1,d_2$ & bicharacters & p \pageref{d}

\\$d_1',d_2'$ & bicharacters & p \pageref{dp}

\\$D_1,D_2$ & bicharacters & p \pageref{gd}

\\$\epsilon$ & root of unity & p \pageref{epsilon}

\\$E_i(m)$ & element of $\mathfrak{K}_i,\mathfrak{K}_i^s$ & p \pageref{eim}

\\$\overset{\leftarrow}{E}_{i,t}(m)$ & element of $\hat{\mathfrak{K}}_{i,t},\hat{\mathfrak{K}}_{i,t}^s $ & p \pageref{teitm}

\\$\overset{\rightarrow}{E}_t(m)$ & element of $\hat{\mathfrak{K}}_{t}^{\infty}$ & p \pageref{tetm}

\\$F(m)$ & element of $\hat{\mathfrak{K}}$ & p \pageref{fm}

\\$\hat{F}_t(m)$ & element of $\hat{\mathfrak{K}}_{t}^{\infty}$ & p \pageref{tftm}

\\$\gamma$ & map $(I\times \ZZ/s\ZZ)^2\rightarrow \ZZ$ & p \pageref{gamma}

\\$\hat{\mathcal{H}}$ & Heisenberg algebra & p \pageref{zq}

\\$\hat{\mathcal{H}}_h$ & formal series in $\hat{\mathcal{H}}$ & p \pageref{zqh}

\\$\mathfrak{K}_i, \mathfrak{K}$ & subrings of $\Yim$ & p \pageref{ki}

\\$\hat{\mathfrak{K}}_{i,t}, \hat{\mathfrak{K}}_t$ & subrings of $\hat{\Yim}_t$ & p \pageref{kit}

\\$\hat{\mathfrak{K}}_{i,t}^{\infty}, \hat{\mathfrak{K}}_t^{\infty}$ & subrings of $\hat{\Yim}_t^{\infty}$ & p \pageref{kitinf}

\\$\mathfrak{K}_i^s, \mathfrak{K}^s$ & subrings of $\Yim^s$ & p \pageref{kis}

\\$\hat{\mathfrak{K}}_{i,t}^s, \hat{\mathfrak{K}}_t^s$ & subrings of $\hat{\Yim}_t^s$ & p \pageref{kits}

\\$\hat{\mathfrak{K}}_{i,t}^{s,\infty}, \hat{\mathfrak{K}}_t^{\infty}$ & subrings of $\hat{\Yim}_t^{s,\infty}$ & p \pageref{kitinfs}

\\$\tilde{\mathfrak{K}}_{i,t}^{s}$ & subring of $\hat{\Yim}_t^s$ & p \pageref{tkits}

\\$\tilde{\mathfrak{K}}_{i,t}^{s, \infty},\tilde{\mathfrak{K}}_t^{s, \infty, f}$ & subring of $\hat{\Yim}_t^{s,\infty}$ & p \pageref{tkitinfs}

\\$\chi_{\epsilon}$ & morphism &

\\ & of $\epsilon$-characters &p \pageref{eps}

\end{tabular}

\begin{tabular}{lll}

 &&

\\ &&

\\$\chi_{\epsilon,t}$ & morphism &

\\ & of $\epsilon,t$-characters &p \pageref{epst}

\\$\chi_q$ & morphism &

\\ & of $q$-characters &p \pageref{chiqdefi}

\\$\chi_{q,t}$ & morphism &

\\ & of $q,t$-characters &p \pageref{chiqt}

\\$[l]$ & element of $\ZZ/s\ZZ$ & p \pageref{lc}

\\$\hat{L}_t(m)$ & element of $\hat{\mathfrak{K}}_t^{\infty}$ & p \pageref{tltm}

\\$\hat{L}_t^s(m)$ & element of $\hat{\mathfrak{K}}_t^{s,\infty}$ & p \pageref{tltms}

\\$\Lambda_s$ & set & p \pageref{lsg}

\\$\overset{\rightarrow}{m},\overset{\leftarrow}{m}$ & $\hat{\Yim}_t^s$-monomial& p \pageref{marrow}

\\$op$ & operator & p \pageref{ope}

\\$p_s$ & morphism & p \pageref{tau}

\\$\pi_r$ & map to $\ZZ$ & p \pageref{pir}

\\$\pi_+$ & ring homomorphism of & p \pageref{piplus}

\\$\hat{\Pi}$ & morphism & p \pageref{pig}

\\$P_{m',m}(t)$ & polynomial & p \pageref{pmm}

\\$P_{m',m}^s(t)$ & polynomial & p \pageref{pmms}

\\$q$ & complex number & p \pageref{qz}

\\$r^{\vee}$ & integer & p \pageref{rvee}

\\$r_i$ & integer & p \pageref{ri}

\\$\text{Rep}$ & Grothendieck ring & p \pageref{rep}

\\$\text{Rep}^s$ & Grothendieck ring & p \pageref{reps}

\\$\text{Rep}_t^s$ & deformed &

\\ &Grothendieck ring & p \pageref{repts}

\\ $s$ & integer & p \pageref{s}

\\$S_i$ & screening operator & p \pageref{si}

\\$S_i^s$ & screening operator & p \pageref{sis}

\\$\tilde{S}_{i,l}$ & screening current & p \pageref{tsil}

\\$S_{i,t}$ & $t$-screening operator & p \pageref{tsit}

\\$S_{i,t}^s$ & $t$-screening operator & p \pageref{tsits}

\\$t$ & central element of $\Yim_t$ & p \pageref{t}

\\$t_R$ & central element of $\Yim_u$ & p \pageref{tr}

\\$\tau_s$ & morphism & p \pageref{taus}

\\$\tau_{s,t}$ & morphism & p \pageref{taust}

\\$u_{i,l}$ & character & p \pageref{uil}

\\$u_{i,l}'$ & character & p \pageref{uilp}

\\$X_{i,l}$ & element of $\text{Rep},\text{Rep}^s$ & p \pageref{reps}

\\$y_i[m]$ & element of $\hat{\mathcal{H}}$ & p \pageref{yim}

\\$Y_{i,l}, Y_{i,l}^{-1}$ & elements of $\Yim$ & p \pageref{yil}

\\$\tilde{Y}_{i,l},\tilde{Y}_{i,l}^{-1}$ & elements of $\hat{\Yim}_u$ or $\hat{\Yim}_t$ & p \pageref{tail}

\\$\Yim$ & commutative algebra & p \pageref{y}

\\$\hat{\Yim}_t^s,\hat{\Yim}_t$ & quotient of $\hat{\Yim}_u^s,\hat{\Yim}_u$ & p \pageref{tyt}

\\$\hat{\Yim}_u^s,\hat{\Yim}_u$ & subalgebra of $\hat{\mathcal{H}}_h$ & p \pageref{yu}

\\$\hat{\Yim}_{i,t}$ & $\hat{\Yim}_t$-module & p \pageref{yit}

\\$\hat{\Yim}_{i,t}^s$ & $\hat{\Yim}_t^s$-module & p \pageref{yits}

\\$\hat{\Yim}_t^{\infty}$ & completion of $\hat{\Yim}_t$ & p \pageref{tytinf}

\\$\hat{\Yim}_t^{s, \infty}$ & completion of $\hat{\Yim}_t^s$ & p \pageref{tytinfs}

\\$z$ & indeterminate & p \pageref{qz}

\\$*$ & $t$-product & p \pageref{star}

\end{tabular}

\onecolumn

\end{document}